\documentclass[11pt]{amsart}

\usepackage{amsmath, amsthm, amssymb, mathrsfs}
\usepackage[hidelinks]{hyperref}
\usepackage[parfill=0pt]{parskip}    \parskip = 0.2cm    
\usepackage{etoolbox} 
\usepackage[margin=1.12in]{geometry} 
\usepackage{tikz}
\usepackage{tikz-cd}
\usetikzlibrary{matrix}
\usepackage{comment} 
\usepackage{subcaption} 
\DeclareCaptionLabelFormat{lc}{\MakeLowercase{#1}~#2}
\usepackage[style=alphabetic, firstinits=true, isbn=false, url=false, eprint=false]{biblatex} 
\usepackage{mathtools}

\graphicspath{ {images/} }
\usepackage{import}

\numberwithin{equation}{section}

\title[LR for polytopal windows]{A characterisation of linear repetitivity for cut and project sets with general polytopal windows}
\date{\today}

\author{James J.\ Walton} 
\address{School of Mathematical Sciences, Mathematical Sciences Building, University Park, Nottingham, NG7 2RD, United Kingdom}
\email{James.Walton@nottingham.ac.uk}
\urladdr{https://www.nottingham.ac.uk/mathematics/people/james.walton}

\theoremstyle{plain}
\newtheorem{theorem}{Theorem}[section]
\newtheorem{lemma}[theorem]{Lemma}

\newtheorem{proposition}[theorem]{Proposition}
\newtheorem{corollary}[theorem]{Corollary}

\newtheorem{Mainthm}{Theorem}

\theoremstyle{plain}
\newtheorem{definition}[theorem]{Definition}
\newtheorem{question}[theorem]{Question}

\theoremstyle{definition}
\newtheorem{remark}[theorem]{Remark}
\AtEndEnvironment{remark}{\null\hfill\exend}
\newtheorem{example}[theorem]{Example}
\AtEndEnvironment{example}{\null\hfill\exend}
\newtheorem{notation}[theorem]{Notation}
\AtEndEnvironment{notation}{\null\hfill\exend}

\newcommand{\R}{\mathbb R}
\newcommand{\C}{\mathbb C}
\newcommand{\Z}{\mathbb Z}
\newcommand{\Q}{\mathbb Q}

\newcommand{\N}{\mathbb N}

\newcommand{\sH}{\mathscr{H}}
\newcommand{\sA}{\mathscr{A}}
\newcommand{\sC}{\mathscr{C}}
\newcommand{\sF}{\mathscr{F}}

\newcommand{\Ver}{\mathrm{Ver}}

\newcommand{\E}{\mathbb E}

\newcommand{\LR}{{\bf LR}}
\newcommand{\Cpx}{{\bf C}}
\newcommand{\D}{{\bf D}}
\newcommand{\DF}{\bf D\(_F\)}
\newcommand{\vol}{\mathrm{vol}}
\newcommand{\Aut}{\mathrm{Aut}}

\renewcommand{\epsilon}{\varepsilon}

\newcommand{\exend}{\hfill \ensuremath{\Diamond}}

\newcommand{\iW}{\mathring{W}}      

\DeclareMathOperator{\rk}{rk}

\DeclareMathOperator{\Bad}{\mathrm{Bad}}

\DeclareMathOperator{\tot}{\E}      
\DeclareMathOperator{\phy}{\E_\vee} 
\DeclareMathOperator{\intl}{\E_<}   
\DeclareMathOperator{\cps}{\Lambda} 


\begin{document}

\begin{abstract}
The cut and project method is a central construction in the theory of Aperiodic Order for generating quasicrystals with pure point diffraction. Linear repetitivity (\LR) is a form of ideal regularity of aperiodic patterns. Recently, Koivusalo and the present author characterised \LR\ for cut and project sets with convex polytopal windows whose supporting hyperplanes are commensurate with the lattice, the weak homogeneity property. For such cut and project sets, we show that \LR\ is equivalent to two properties. One is a low complexity condition, which may be determined from the cut and project data by calculating the ranks of the intersections of the projection of the lattice to the internal space with the subspaces parallel to the supporting hyperplanes of the window. The second condition is that the projection of the lattice to the internal space is Diophantine (or `badly approximable'), which loosely speaking means that the lattice points in the total space stay far from the physical space, relative to their norm. We review then extend these results to non-convex and disconnected polytopal windows, as well as windows with polytopal partitions producing cut and project sets of labelled points. Moreover, we obtain a complete characterisation of {\LR} in the fully general case, where weak homogeneity is not assumed. Here, the Diophantine property must be replaced with an inhomogeneous analogue. We show that cut and project schemes with internal space isomorphic to \(\R^n \oplus G \oplus \Z^r\), for \(G\) finite Abelian, can, up to MLD equivalence, be reduced to ones with internal space \(\R^n\), so our results also cover cut and project sets of this form, such as the (generalised) Penrose tilings.
\end{abstract}


\maketitle

\centerline{\em Dedicated to Uwe Grimm.}

\section*{Introduction}

The cut and project method \cite{Moo00} is perhaps the most important general construction of quasicrystals \cite{Mey95,Lag96,GriKra19} with pure point diffraction \cite{BaaGri12}. Within the field of Aperiodic Order \cite{BGD16, AOI}, cut and project schemes naturally arise from lifting aperiodic structures \cite{BEG16}. Lagarias and Pleasants posed the problem of characterising linear repetitivity for cut and project sets, see \cite[Problem 8.3]{LP03}. See also \cite{DHS99}, which concerns the dynamical properties of linearly repetitive or `recurrent' subshifts in the dimension \(d=1\) case. Generally, a Delone set or tiling of \(\mathbb{R}^d\) is called linearly repetitive (\LR) if there exists a constant \(C > 0\) for which every patch of radius \(r \geq 1\) can be found within distance \(Cr\) of any point in the pattern. Linear repetitivity is a very strong property of an aperiodic pattern, signifying a high degree of rigid structural order. Amongst other properties, {\LR} implies rapid convergence of patch counts to ergodic averages \cite{LP03} and uniform subadditive ergodic theorems \cite{BBL} which are important for the study of lattice gas models and random Schr\"{o}dinger operators \cite{GeeHof91,DamLen01}. It is know that {\LR} Delone sets have a finite number of nonperiodic Delone system factors, see \cite{CDP10}. For further properties of \LR\ Delone sets, we refer the reader to \cite{ACCDP15}.

Characterising {\LR} for cut and project sets in full generality is difficult because of the degrees of freedom involved, such as in the choice of the window. However, in recent work \cite{KoiWalI,KoiWalII}, the problem was settled for Euclidean cut and project sets with convex polytopal windows, satisfying an additional condition of `weak homogeneity' (Definition \ref{def: homogeneous}), meaning that the supporting hyperplanes of the window are in some sense mutually commensurate with the projection of the lattice to internal space. This result greatly generalises the main result of \cite{HayKoiWal18}, on `cubical windows', and also covers the widely studied canonical cut and project sets, those with Euclidean total space, lattice \(\Z^k\) and window the projection of the unit hypercube \([0,1]^k\) to the internal space (provided the lattice projects densely to the internal space). For cut and project sets with Euclidean internal space, it seems that {\LR} (and even low complexity, see below) requires the window to either be polytopal or have fractal boundary (Rauzy fractals being important examples \cite{Pyt02,BaaGri20,BPfG21}), as one needs large overlaps of the boundary with itself under translation to avoid producing many acceptance domains and thus high complexity, concepts that are explained later in this article. So we consider the resolution of the polytopal window case as a significant milestone.

An obviously necessary (see comments below Definition \ref{def: low complexity}) condition for {\LR} is `low complexity', which we denote by property \Cpx, meaning that, for some \(C > 0\), there are at most \(Cr^d\) patches of any radius \(r \geq 1\). This is the minimal polynomial growth rate of the complexity function \(p(r)\), that is, for any nonperiodic polytopal cut and project set there is some \(c > 0\) so that \(p(r) \geq c r^d\) (Lemma \ref{lem: minimum complexity}). Despite the large amount of defining data, in the weakly homogeneous case the answer to which cut and project sets are {\LR} is strikingly neat: given {\Cpx} there is just one additional necessary and sufficient condition for {\LR}, a Diophantine condition {\D} on the positioning of the lattice relative to the physical space. This is the requirement that all lattice points are distant from the physical space, relative to their distance from the origin (in a particular quantitative sense of course, as given in Definition \ref{def: Diophantine scheme}, and in a sense that we only define here in the presence of property {\Cpx} in the atypical case that the window is decomposable or a `product'). For \(2\)-to-\(1\) canonical cut and project schemes, this is equivalent to the slope of the physical space being badly approximable relative to the lattice, so \cite[Theorem A]{KoiWalII} extends part of a result of Morse and Hedlund on {\LR} for Beatty sequences \cite{MH38,MH40}.

To demonstrate the degree of data reduction in the above result, consider Figures \ref{fig: cps}, \ref{fig: messy internal} and \ref{fig: easy internal}. Figure \ref{fig: cps} illustrates the full data required for a low dimensional (\(2\)-to-\(1\)) cut and project scheme. Given the standard restriction of the lattice projecting injectively to the physical space, all data required to determine {\LR} may be considered within the internal space, via the machinery of acceptance domains (reviewed in Section \ref{sec: basic definitions}). This data is illustrated in Figure \ref{fig: messy internal} (in this case for a \(4\)-to-\(2\) scheme). Here, we allow for a generalisation of the setting of \cite{KoiWalII}, by allowing the polytopal window to be nonconvex, disconnected or even equipped with a polytopal partition into regions that assign `colours' or `labels' to the points of the cut and project set. In the weakly homogeneous case, {\LR} is determined via two simple conditions on the data illustrated in Figure \ref{fig: easy internal}, consisting only of the set \(\sH_0\) of subspaces supporting the boundary of the window (or partition regions), and \(\Gamma_<\), the projection of the lattice to internal space, both of which can be described by a finite set of vectors in the internal space (the normals to the supporting hyperplanes and a choice of \(k\) generating vectors of the projected lattice).

Generally, even without weak homogeneity, the degree \(\alpha \in \N\) of the patch counting function which satisfies (using Vinogradov) \(p(r) \asymp r^\alpha\) (and thus also property \Cpx) is still determined just by the pair \((\sH_0,\Gamma_<)\), as we will see in Theorem \ref{thm: generalised complexity} (see also the earlier work of \cite{Jul10}). In particular, see Corollary \ref{cor: sum of ranks for low complexity}, which shows that {\Cpx} holds if and only if the sum of ranks of restrictions of \(\Gamma_<\) (the projection of the lattice to internal space) to subspaces parallel to the hyperplanes supporting the window over any flag (Definition \ref{def: flags}) are equal to \(k(n-1)\), where \(k\) is the `total dimension' and \(n\) is the `codimension' (in particular, \Cpx\ holds for every codimension \(n=1\) example but only holds on a measure 0 set of parameters in codimension \(n > 1\)). Property {\D} is determined just by \(\Gamma_<\), considered as a subgroup of the Euclidean internal space, at least in the standard `indecomposable' case, see \cite{KoiWalII} or Definition \ref{def: decomposable} (otherwise, one needs to consider the individual factors of a certain splitting of the lattice, given by \(\sH_0\), which is only guaranteed in the presence of the simpler condition {\Cpx}, see Remark \ref{rem: independence of C and D}). The set \(\Bad(m,n)\) of \(m\) badly approximable forms in \(n\) variables is known to be full Hausdorff dimension but a null set, see for example \cite{Sch66}, from which it follows that \D\ is also an atypical property in all codimensions. In this paper we show that \cite[Theorem A]{KoiWalII} may be generalised to arbitrary (i.e., labelled and not necessarily convex) polytopal windows:

\begin{Mainthm} \label{thm: main1}
For a weakly homogeneous polytopal cut and project scheme, the following are equivalent:
\begin{enumerate}
	\item \LR;
	\item {\Cpx} and \D.
\end{enumerate}
\end{Mainthm}

The definition of `polytopal' required for our proofs is given in Definition \ref{def: polytopal}, although this turns out to be equivalent to being a finite union of convex polytopes (Remark \ref{rem: nonconvex polytopes}).

Next we consider schemes that are not weakly homogeneous. To avoid small acceptance domains of \(r\)-patches (and thus infrequent patches), one needs to force vertices of acceptance domains to remain distant relative to \(r\). We consider the `generalised vertices' (see Figure \ref{fig: generalised vertices}) of the window, given by singleton intersection points of supporting hyperplanes of the window's boundary (or of its labelling partition, in the labelled case), and the set \(F\) of displacements between these. Then, at least when {\Cpx} holds, the set of displacements between vertices of acceptance domains (and the simpler `cut regions' that may be used to approximate them) are closely related to the set \(\Gamma_< + F\), where \(\Gamma_<\) is the projection of the lattice to internal space (Remark \ref{rem: C versus finite index vertices}). This motivates an inhomogeneous Diophantine condition, called {\DF}, requiring that projections of lattice points to internal space remain distant from \(F\) relative to their norm in the total space, in a particular sense made precise in Definition \ref{def: DF} (which, again, needs care in the decomposable case). We obtain the following:

\begin{Mainthm} \label{thm: main2}
For a polytopal cut and project scheme, if \LR\ holds then so does {\Cpx} and {\DF}. Conversely, there is some \(N \in \N\) (which may be taken as \(N=1\) in codimension \(n=1\)) so that if {\Cpx} holds and {\DF} still holds after replacing the lattice \(\Gamma\) with \(\frac{1}{N}\Gamma\), then \LR\ holds.
\end{Mainthm}

The value of \(N \in \N\) in the above may be determined from the data of \(\Gamma_<\) (the projection of the lattice to the internal space) and the set \(\sH_0\) of subspaces parallel to the support hyperplanes (that is, it does not depend on which translations are used to cover the windows boundary). One would expect the condition {\DF} to be typically quite difficult to verify or preclude but, in any case, this comes down to a question within Diophantine Approximation.

In fact, we are able to derive a single necessary and sufficient condition for {\LR}, in Theorem \ref{thm: iff for LR}. It is similar to the above, in that it requires {\Cpx} (which may be checked directly from the data of the cut and project scheme using Theorem \ref{thm: generalised complexity}) and a Diophantine condition, and it still only depends on \(\Gamma_<\) and the set \(\sH\) of supporting hyperplanes of the window. However, the Diophantine condition is now more complicated. Given {\Cpx} we may define certain groups \(\Gamma[f,f']\) (see Notation \ref{not: double flag groups} and Definition \ref{def: lifted vertex groups}) of the total space, given in terms of \(\Gamma_<\) and any pair of flags (Definition \ref{def: flags}) \(f\) and \(f'\), where \(\Gamma\) is finite index in each \(\Gamma[f,f']\). Then our characterisation of {\bf LR} is that {\bf C} holds and that each projection \(\Gamma[f,f']_<\) of \(\Gamma[f,f']\) to the internal space satisfies a certain inhomogeneous Diophantine condition of avoiding the relative displacement \(v(f')-v(f)\) between the flags' centres. Theorem \ref{thm: iff for LR} also establishes that {\LR} is equivalent (not just sufficient) to validity of a subadditive ergodic theorem \cite{BBL}. It is, in part, the introduction of the groups \(\Gamma[f,f']_<\) (which may be defined in terms of \(\Gamma_<\) and \(\sH_0\)) which allows us to extend the main result of \cite{KoiWalII} characterising {\LR}. Theorems \ref{thm: main1} and \ref{thm: main2} then quickly follow from the fully general Theorem \ref{thm: iff for LR}.

In the weakly homogeneous case, {\D} and {\DF} are equivalent over all lattices \(\frac{1}{N}\Gamma\) for fixed \(\Gamma\) here (see Proposition \ref{prop: F for weakly homogeneous} and Remarks \ref{rem: hom versus inhom Diophantine}, \ref{rem: DF => D}), so Theorem \ref{thm: main2} generalises Theorem \ref{thm: main1}. Indeed, passing between finite index subgroups of lattices does not affect property {\D} (Lemma \ref{lem: Diophantine under finite index}). Unfortunately, it is not clear that this is the case for the inhomogeneous condition {\DF} in general. Given a Diophantine lattice (Definition \ref{def: Diophantine}) one may ask if being inhomogeneously Diophantine with respect to some \(F\) (Definition \ref{def: inhomogeneous Diophantine}) is preserved under rational rescaling. This may be a subtle question in Diophantine Approximation and, to the author's knowledge, is not even known in the lowest possible dimensions, see Question \ref{q: question in DA}. If invariance under rational rescaling does hold, then the above result would show that {\LR} is equivalent to {\Cpx} and {\DF}, for all polytopal windows, i.e., one may remove the weakly homogeneous assumption in Theorem \ref{thm: main1} by replacing {\D} with {\DF}. In any case, Theorem \ref{thm: main2} gives two tightly related necessary and sufficient conditions that coincide for weakly homogeneous and codimension \(n=1\) schemes. If Question \ref{q: question in DA} answers in the negative, then the more complicated Diophantine property of Theorem \ref{thm: iff for LR} is essentially necessary for the complete characterisation of {\LR}.

Many interesting examples of cut and project sets have convex polytopal windows with multiple Euclidean components, such as the (generalised) Penrose tilings, which have internal space \(\Z^2 \oplus (\Z/5\Z)\). In Section \ref{sec: multiple components} we show that, up to MLD equivalence \cite[Section 5.2]{AOI}, cut and project schemes with internal spaces \(\R^n \oplus \Z^r \oplus G\) (for \(r \in \N \cup \{0\}\) and \(G\) finite Abelian) may be reduced to ones with internal space just \(\R^n\). The process easily determines the new projected lattice \(\Z_<\) and a new window, in particular the set \(\sH\) of supporting hyperplanes, which is the data needed for our main theorems characterising \LR. Other interesting examples of cut and project schemes have internal spaces with \(p\)-adic components, such as the cut and project scheme for the period doubling sequences, see, for example, \cite[Section 3]{BaaMoo04}. One would imagine that a few techniques from this paper may have analogues for \(p\)-adic internal spaces, or conjecture that \LR\ is equivalent to a low complexity condition coupled with a Diophantine condition. However, it is not clear to the author what a natural class of windows would be (to replace the polytopal windows here) where a characterisation seems attainable. Many structures and results in this paper depend on the Euclidean hyperplanes that support the window's boundary.

The paper is organised as follows. In Section \ref{sec: basic definitions} we cover the basic definitions for polytopal cut and project sets, their complexity, repetitivity, acceptance domains and decomposability. In Section \ref{sec: complexity} we review constructions from \cite{KoiWalI} (see also \cite{Jul10}) that determine the complexity of a general polytopal cut and project set. The main result here is Theorem \ref{thm: generalised complexity}, giving a formula for the complexity exponent \(\alpha\) of \(p(r) \asymp r^\alpha\). This is a mild extension of \cite[Theorem 6.1]{KoiWalI} with only a few details of the proof needing to be adjusted, but we cover the essential ideas of the proof for completeness. In Section \ref{sec: Diophantine schemes} we review the notion of (families of) numbers or vectors being badly approximable, and also a reference lattice free version of this: Diophantine densely embedded lattices. We introduce an inhomogeneous Diophantine condition in Definition \ref{def: inhomogeneous Diophantine}. These definitions are then applied to the projection of the lattice to internal space to define the homogeneous Diophantine property {\D} and inhomogeneous Diophantine property {\DF} of a polytopal cut and project scheme, where the latter may, a priori (subject to the answer to Question \ref{q: general question in DA}), depend on how the lattice is split in the decomposable case. In Section \ref{sec: main proof} we prove Theorem \ref{thm: iff for LR}, from which Theorems \ref{thm: main1} and \ref{thm: main2} easily follow. In Section \ref{sec: multiple components} we show how to reduce, up to MLD equivalence, a scheme with multiple Euclidean components to one with just a single Euclidean component. Finally, in Section \ref{sec: examples} we apply our main theorems to some illustrative examples.

\subsection*{Acknowledgements}

The author thanks Sam Chow, Franz G\"{a}hler, Alan Haynes, Johannes Hofscheier and Henna Koivusalo for useful discussions.

\section{Basic definitions} \label{sec: basic definitions}

We begin by reviewing the definition of a polytopal cut and project set, its complexity, repetitivity and acceptance domains. For ease of reference, we mostly use the notation of \cite{KoiWalI,KoiWalII}. Throughout, by a `vector space' we mean a finite-dimensional vector space over \(\R\).

\subsection{Labelled polytopal cut and project sets} \label{sec: cut and project sets}

\begin{figure}
	\def\svgwidth{\textwidth}
	\centering
	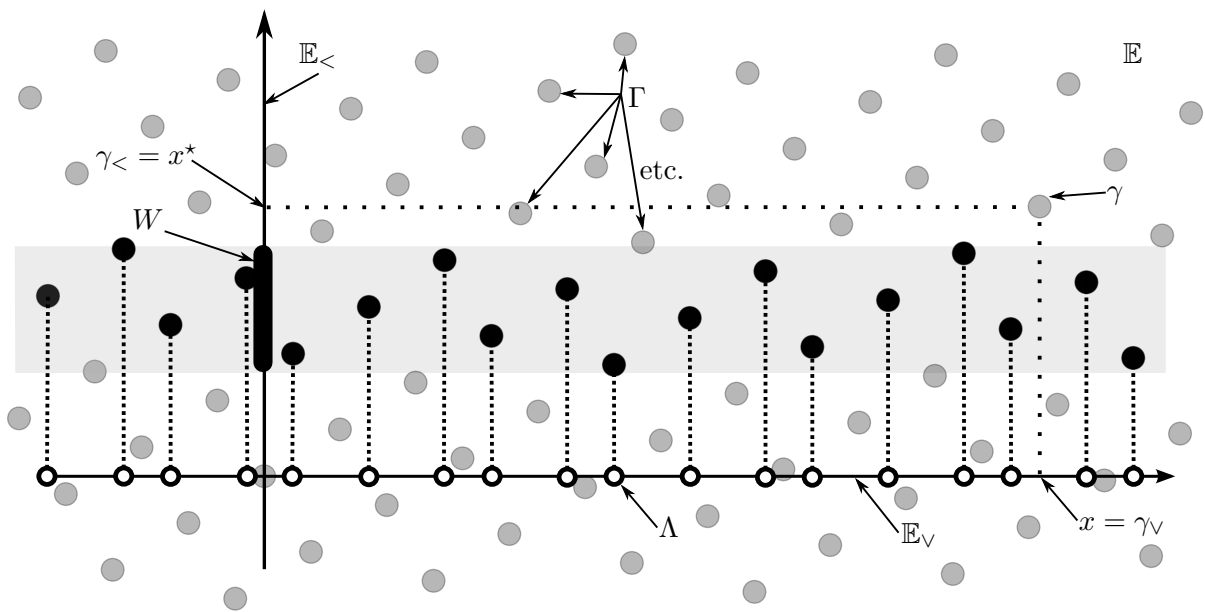
	\caption{The data of a cut and project set, here in the \(2\)-to-\(1\) case, with demonstration of the notation \(\gamma_<\), \(\gamma_\vee\) and \(x^\star\) for a particular \(\gamma \in \Gamma\) and \(x = \gamma_\vee\).}\label{fig: cps}
\end{figure}

A (Euclidean) {\bf \(k\)-to-\(d\) cut and project scheme} \(\mathcal{S}\), where \(d \in \N\) is the {\bf dimension} and \(n \coloneqq k-d \in \N\) is the {\bf codimension}, consists of the following data: a \(k\)-dimensional vector space \(\tot\) called the {\bf total space}, a \(d\)-dimensional subspace \(\phy < \tot\) called the {\bf physical space}, a complementary \(n\)-dimensional subspace \(\intl < \tot\) called the {\bf internal space}, a lattice (i.e., a discrete co-compact subgroup) \(\Gamma < \tot\) and a {\bf window} \(W \subset \intl\).

\begin{notation}
The decomposition \(\tot = \phy + \intl\) induces projections \(\pi_\vee \colon \tot \to \phy\) and \(\pi_< \colon \tot \to \intl\). For \(x \in \tot\) we write \(x_\vee \coloneqq \pi_\vee(x)\). We use analogous notation for the projections of subsets, and projections to the internal space. The notation is visual (see Figure \ref{fig: cps}): it is remembered by the direction of arrows, when in the first quadrant of a \(2\)-to-\(1\) scheme. For example, we think of \(\phy\) as the space we project `downwards' to (when in the first quadrant). The notation is similar to a widely used notation \(\tot_\|\) for the physical (or `parallel') space and \(\tot_\bot\) for the internal (or `complementary') space.
\end{notation}

We will always assume that our scheme \(\mathcal{S}\) is {\bf polytopal}, meaning that \(W\) is polytopal. Unlike for convex polytopal, there are multiple (inequivalent) possible definitions. The one below is general (see Remark \ref{rem: nonconvex polytopes}) and is directly convenient to our setting.

\begin{definition} \label{def: polytopal}
A {\bf hyperplane} \(H\) in a vector space \(X\) is a translate of a codimension \(1\) subspace. A (closed) {\bf half-space} \(H^+\) for \(H\) is then the union of \(H\) and one of the two connected components of \(X \setminus H\).

Given a subset \(W \subset X\), we call \(H \subset X\) {\bf properly supporting} if there exists some \(w_H \in H\) and open neighbourhood \(U_H \subset X\) of \(w_H\) for which \mbox{\(W \cap U_H = H^+ \cap U_H\)}, where \(H^+\) is a closed half-space for \(H\). A subset \(W \subset \R^n\) is called {\bf polytopal} if \(W\) is compact, equal to the closure of its interior and has boundary \(\partial W\) contained in the union of a finite set \(\sH\) of properly supporting hyperplanes. 
\end{definition}

\begin{notation}
For any polytope \(W\), every properly supporting hyperplane \(H\) must be a member of \(\sH\), otherwise it would be impossible to cover \(\partial W\) about \(w_H\). For brevity, we drop `properly' and refer to \(\sH\) as the set of {\bf supporting hyperplanes}. For a hyperplane \(H\), we let \(V(H) \coloneqq H-H = H-h\) for any \(h \in H\), i.e., the subspace parallel to \(H\). We denote the {\bf supporting subspaces} of a polytope by \(\sH_0 \coloneqq V(\sH) = \{V(H) \mid H \in \sH \}\).
\end{notation}

\begin{example} \label{ex: convex polytope}
Let \(W\) be a convex polytope, that is, a compact subset with non-empty interior that is an intersection of finitely many closed half-spaces, or equivalently the convex hull of a finite set of points with non-empty interior. Then \(W\) is polytopal and \(\sH\) is the set of `irredundant' supporting hyperplanes, i.e., the smallest set of hyperplanes covering the boundary.
\end{example}

\begin{remark} \label{rem: nonconvex polytopes}
Generally, if \(W\) is polytopal, \(\partial W \subset \bigcup_{H \in \sH} H\) for a finite set \(\sH\) of hyperplanes and thus \(\iW\) is a (necessarily finite) union of bounded connected components of \(X \setminus \bigcup_{H \in \sH} H\), ignoring any intersections of supporting hyperplanes with \(\iW\). Since \(W\) is the closure of its interior, we see that \(W\) is a union of closures of such regions (and thus a finite union of convex polytopes). Conversely, given a union \(A\) of bounded open subsets of the complement of a finite union of hyperplanes, the closure \(W = \overline{A}\) is polytopal. To see this, note that \(W\) is a finite union of convex polytopes (the closures of the open components defining \(A\)). Moreover, these polytopes have an induced decomposition into convex cellular faces (with vertices at certain intersections of \(\dim(X)\) many hyperplanes), giving a convex cellular decomposition of \(W\). In particular, \(\partial W\) is a union of convex polytopes of dimension \(\dim(X)-1\), and interior points of these provide the \(w_H\) exhibiting the properly supporting property of Definition \ref{def: polytopal}. This argument shows that \(W\) is polytopal, in the sense of Definition \ref{def: polytopal}, if and only if it is a finite union of convex polytopes (Example \ref{ex: convex polytope}).
\end{remark}

In our setup, we will more generally allow our window to possess a polytopal partition, which produces labelled cut and project sets. More precisely, we may suppose that \(W = \bigcup_{i=1}^\ell W_i\), where each \(W_i\) is polytopal and distinct \(W_i\) intersect on at most their boundaries. We will see in Proposition \ref{prop: unlabel} that adding labels does not present much of a generalisation: up to MLD equivalence \cite[Section 5.2]{AOI}, we may replace the labelled window \(W\) with an unlabelled one. Nevertheless, it seems of value to allow our main theorems to apply directly to them.

In summary, a {\bf polytopal cut and project scheme} is a tuple \(\mathcal{S} = (\tot,\phy,\intl,\Gamma,W)\), where \(W\) is polytopal and may also carry a polytopal partition, as above. We let \(\sH\) denote the supporting hyperplanes and \(\sH_0\) the supporting subspaces of the window (or the union of such over each polytope in the partition, in the labelled case). We make three final standard assumptions:
\begin{enumerate}
	\item \(\pi_\vee\) is injective on \(\Gamma\);
	\item \(\pi_<\) is injective on \(\Gamma\);
	\item \(\Gamma_<\) is dense in \(\intl\).
\end{enumerate}
By Assumption (1), for each \(x \in \Gamma_\vee\) there is a unique lift \(\gamma \in \Gamma\) with \(\gamma_\vee = x\). This allows us to define the star map \(x \mapsto x^\star\), from \(\Gamma_\vee\) to \(\Gamma_<\), by first lifting to the lattice, then projecting to the internal space (by assumption (2) this is an isomorphism from \(\Gamma_\vee\) to \(\Gamma_<\), although a non-continuous one with respect to the subspace topologies). It is well known and a simple exercise to show (even in the absence of (1) and (3)) that (2) is equivalent to non-periodicity of \(\cps\).

We assume (by translating \(W\) if necessary) that the scheme is {\bf non-singular}, meaning that \(\gamma_< \notin \partial W\) (or \(\gamma_< \notin \partial W_i\) for each \(i\), in the labelled case) for all \(\gamma \in \Gamma\). The {\bf cut and project set} determined by \(\mathcal{S}\) is then
\begin{equation} \label{eq: cps}
\cps \coloneqq \{x \in \Gamma_\vee \mid x^\star \in W\},
\end{equation}
where \(x\) has label \(i\) if \(x^\star \in W_i\), in the labelled case (there is precisely one such \(W_i\) by non-singularity). The points of \(\cps\) are those that we \emph{cut} from \(\Gamma\), by keeping only those lattice points in the strip \(W + \phy\), and then \emph{project} to the physical space, see Figure \ref{fig: cps}.

\begin{figure}
  \begin{subfigure}{.45\textwidth}
  \renewcommand\captionlabelfont{}
    \centering\def\svgwidth{\textwidth}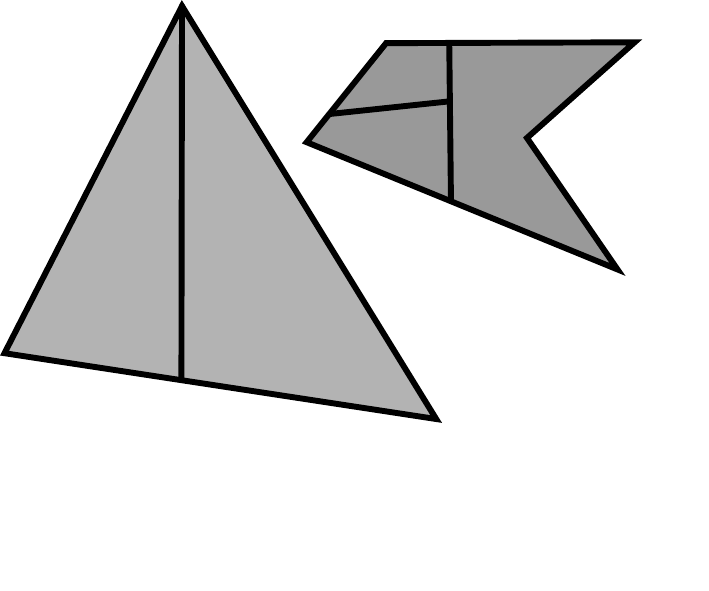%
    \caption{The window does not need to be convex, or even connected, and similarly for the partition labels for labelled windows (here, the window has a partition into 3 labels). Basis vectors of \(\Gamma_<\) are drawn.}
    \label{fig: messy internal}
  \end{subfigure}\hfill
  \begin{subfigure}{.45\textwidth}
  \renewcommand\captionlabelfont{}
    \centering\def\svgwidth{\textwidth}
\begingroup%
  \makeatletter%
  \providecommand\color[2][]{%
    \errmessage{(Inkscape) Color is used for the text in Inkscape, but the package 'color.sty' is not loaded}%
    \renewcommand\color[2][]{}%
  }%
  \providecommand\transparent[1]{%
    \errmessage{(Inkscape) Transparency is used (non-zero) for the text in Inkscape, but the package 'transparent.sty' is not loaded}%
    \renewcommand\transparent[1]{}%
  }%
  \providecommand\rotatebox[2]{#2}%
  \newcommand*\fsize{\dimexpr\f@size pt\relax}%
  \newcommand*\lineheight[1]{\fontsize{\fsize}{#1\fsize}\selectfont}%
  \ifx\svgwidth\undefined%
    \setlength{\unitlength}{340.15748031bp}%
    \ifx\svgscale\undefined%
      \relax%
    \else%
      \setlength{\unitlength}{\unitlength * \real{\svgscale}}%
    \fi%
  \else%
    \setlength{\unitlength}{\svgwidth}%
  \fi%
  \global\let\svgwidth\undefined%
  \global\let\svgscale\undefined%
  \makeatother%
  \begin{picture}(1,0.83333333)%
    \lineheight{1}%
    \setlength\tabcolsep{0pt}%
    \put(0,0){\includegraphics[width=\unitlength,page=1]{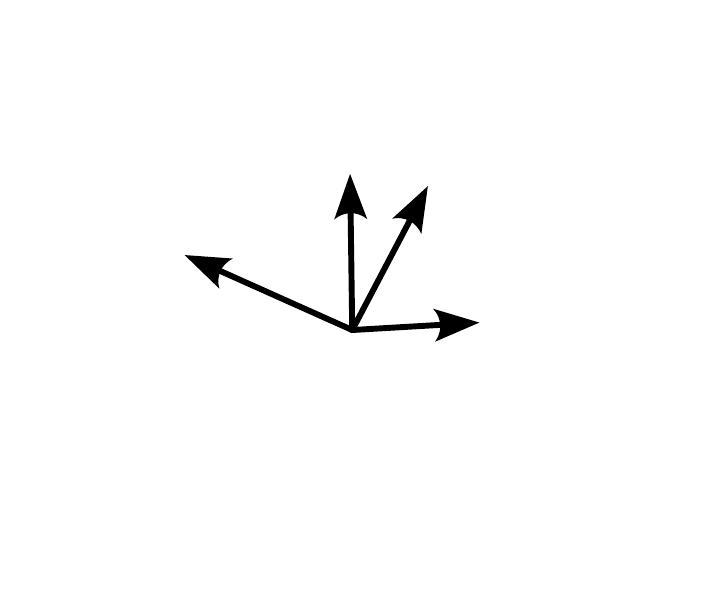}}%
    \put(0.33708437,0.29877744){\color[rgb]{0,0,0}\makebox(0,0)[lt]{\lineheight{1.25}\smash{\begin{tabular}[t]{l}$\Gamma_<$\end{tabular}}}}%
    \put(0,0){\includegraphics[width=\unitlength,page=2]{simple_internal.pdf}}%
    \put(0.0392438,0.73725811){\color[rgb]{0,0,0}\makebox(0,0)[lt]{\lineheight{1.25}\smash{\begin{tabular}[t]{l}$\intl$\end{tabular}}}}%
  \end{picture}%
\endgroup%
    \caption{The same internal space data but with window on the left replaced by its set \(\sH_0\) of supporting subspaces, parallel to those supporting \(W\) (and, in this labelled case, the partition regions). This data is sufficient to determine {\LR} for weakly homogeneous schemes.}
    \label{fig: easy internal}
  \end{subfigure}
  \caption{Internal space data of a \(4\)-to-\(2\) scheme.}
\end{figure}

\subsection{Patches, complexity and repetitivity}
Equip both \(\phy\) and \(\intl\) with norms, and set norm on \(\tot\) as \(\|x\| = \max \{\|x_\vee\|, \|x_<\| \}\). Given \(x \in \cps\), we define the {\bf \(r\)-patch} \(P_r(x)\) to be the intersection of \(\cps\) with the closed \(r\)-ball \(B_r(x)\) centred at \(x\), where \(x\) is additionally recorded as the `centre' of the patch. If \(\cps\) carries a labelling then of course all points of \(P_r(x)\) also carry these labels. Throughout, we will identify two \(r\)-patches whenever they equal, up to a translation taking the centre of one to the centre of the other.

\begin{remark}
One can more generally consider patches of other shapes, such as \(\cps \cap K\) for a bounded set \(K\) (see the definition of a patch, or cluster, in \cite{AOI}). Although not every finite patch in this sense is an \(r\)-patch (such as a two-point patch in \(\Z\)), this does not affect the properties of interest below, since every finite patch is contained in some \(r\)-patch. Similarly, the choice of labelling all patches with their `centre' does not affect properties of interest here, such as the growth rates of the complexity of repetitivity functions, see, for instance, \cite[Proposition 1.11]{Jul10}.
\end{remark}

We then define the {\bf complexity} (or {\bf patch counting}) {\bf function} \(p \colon \R_{\geq 0} \to \N\) by \(p(r) \coloneqq\) number of distinct \(r\)-patches. Note that all cut and project sets here are of {\bf finite local complexity}, that is, \(p(r) < \infty\) for all \(r \geq 0\).

The {\bf repetitivity function} \cite[Section 5.3]{AOI} \(\rho \colon \R_{\geq 0} \to \R_{\geq 0}\) is defined by setting \(\rho(r)\) to be the smallest value of \(R\) for which every \(R\)-patch contains every \(r\)-patch. All cut and project sets here are automatically {\bf repetitive}, that is, \(\rho(r) < \infty\) for all \(r \geq 0\), since by definition all cut and project sets are constructed from non-singular translates of the window.

\begin{notation}
We use Vinogradov Notation: for functions \(f\), \(g \colon \R_{\geq 0} \to \R_{\geq 0}\), we write \(f \ll g\) to mean that there is some \(C > 0\) with \(f(x) \leq Cg(x)\) for sufficiently large \(x\). We write \(f \asymp g\) if \(f \ll g\) and \(g \ll f\).
\end{notation}

The functions \(p\) and \(\rho\) give two measures of orderliness of \(\cps\). If \(p\) grows slowly then \(\cps\) may be considered as less complicated and thus more ordered. For (convex) polytopal cut and project sets, it is known that \(p(r) \asymp r^\alpha\) for \(\alpha \in \N\), \(d \leq \alpha \leq nd\) and \(\alpha\) has a simple description in terms of the defining data; see \cite{Jul10, KoiWalI} or the generalisation of these in Theorem \ref{thm: generalised complexity} later. In particular, \(r^d\) is the slowest polynomial rate of growth in our setting:

\begin{definition} \label{def: low complexity}
We say that \(\cps\) is of {\bf low complexity} (and has property {\Cpx}, for short) if \(p(r) \ll r^d\).
\end{definition}

Slow growth of \(\rho\) means that patches repeat in \(\cps\) with relatively small gaps. The slowest possible growth is linear \cite{LP03}, which is thus singled out as an idealised notion of regularity (for the follow, see also \cite[Section 5.3]{AOI}):

\begin{definition}
We say that \(\cps\) is {\bf linearly repetitive} (and is {\LR}, for short) if \(\rho(r) \ll r\).
\end{definition}

Clearly \LR\ implies \Cpx\ because Delone sets are uniformly discrete, so there are \(\ll r^d\) distinct potential centres of patches in a ball of radius \(r\). The properties {\Cpx} and {\LR} (and more generally, the asymptotic growth rates of \(p\) and \(\rho\)) are MLD invariants which, roughly speaking, means that they are not affected by making cosmetic changes to \(\cps\) governed by local (or `pattern equivariant') rules which are also reversible by local rules. It is not hard to show that any tiling admitted by a primitive FLC substitution rule (see, for instance, \cite{PF08}) is {\LR} and thus so is any Delone set which is MLD to such a substitution tiling, which includes the examples of the Ammann--Beenker and Penrose cut and project sets. See \cite[Section 5.2]{AOI} for more details on substitutions and the notions of MLD and FLC.

All cut and project sets considered here are non-periodic and repetitive. Equation \eqref{eq: cps} defines one cut and project set from a scheme, but in fact a scheme defines an infinite family, by first translating the lattice and then cutting and projecting (at least in the non-singular case; for the singular translates one uses limits of non-singular patterns). Every finite patch in one cut and project set defined by the scheme belongs to any other, that is, each is {\bf locally isomorphic}; see \cite[Section 4.2]{AOI} for the concepts of locally indistinguishable and local isomorphism. In particular, the properties of {\LR} and {\Cpx} do not depend on which cut and project set is chosen, so we may just work with \(\cps\) with no loss of generality.

\subsection{Acceptance domains}
Given any \(r\)-patch \(P\), consider the indicator set \(\cps_P\), given by restricting \(\cps\) to only those \(x \in \cps\) for which \(P_r(x) = P\) (as always, up to translation). Then, as we will see in this subsection, \(\cps_P\) itself is a cut and project set, with the same scheme except with a different window. That window is called an {\bf acceptance domain}. Counting the number of these acceptance domains thus counts \(r\)-patches, and the sizes of the acceptance domains determine the rate of recurrence of patches. Thus, controlling the number and size of acceptance domains allows us to bound the functions \(p\) and \(\rho\).

Here we will just cover the idea of acceptance domains for unlabelled windows, which will be sufficient for our purposes, by Proposition \ref{prop: unlabel} (although the definition in the labelled case is similar). For the following definition, recall that by definition our \(r\)-patches \(P\) have a specified `centre' (which is the centre for the \(r\)-ball used in the intersection to create \(P\)), and that we identify patches that agree up to translation so that we may always choose the unique representative with centre over the origin.

\begin{definition}
For \(r \geq 0\) we let \(\Gamma(r) \coloneqq \Gamma \cap B_r\) (where \(B_r \coloneqq B_r(0)\) is the closed \(r\)-ball centred at the origin). For an \(r\)-patch \(P\), translated with centre over the origin, we consider the following subsets of the lattice:
\begin{align*}
\Gamma^\mathrm{in}_P \coloneqq  & \{x \in \Gamma(r) \mid x_\vee \in P\}, \\
\Gamma^\mathrm{out}_P \coloneqq & \{x \in \Gamma(r) \mid x_\vee \notin P\}.\\
\end{align*}
The {\bf acceptance domain} of \(P\) is the subset \(A_P \subset W\) defined by
\[
A_P \coloneqq \left( \bigcap_{x \in \Gamma^\mathrm{in}_P} (\iW - x_<) \right) \cap \left( \bigcap_{x \in \Gamma^\mathrm{out}_P} (W^\mathrm{c} - x_<) \right).
\]
We denote the set of acceptance domains of \(r\)-patches by \(\sA(r)\).
\end{definition}

It follows from the above that the acceptance domains are (interiors of) polytopes. Note that they need not be convex, even if \(W\) is convex, due to the intersections with translates of \(W^\mathrm{c}\).

\begin{remark} \label{rem: origin of cut regions}
Suppose that \(y \in \cps\), so that \(y^\star \in W\). Then \(y+u \in \cps\) if and only if \(u = x_\vee\) for some \(x \in \Gamma\) with \(y^\star + x_< \in W\), which happens if and only if \(y^\star \in W - x_<\). The proof of Lemma \ref{lem: acceptance domains} below follows quickly from this, see \cite{KoiWalI} for further details. However, we note a minor modification here, by allowing \(x \in \Gamma(r)\) in the definitions of the in/out sets of patches, that is, \(\|x_\vee\|\), \(\|x_<\| \leq r\). In \cite{KoiWalI} the latter is replaced with \(x_< \in W-W\). This yields the same definition once \(B_r \supseteq W-W\). Indeed, in this case, we use just as many lattice translates to define each \(A_P\), but any \(x \in \Gamma(r)\) with \(x_< \notin W-W\) is irrelevant, since automatically \(A_P \subset \iW\) (as \(0 \in \Gamma_P^\mathrm{in}\)) and \((W - x_<) \cap W = \emptyset\) in this case, thus \(x \in \Gamma_P^\mathrm{out}\) with \(\iW \cap (W^\mathrm{c} - x_<) = \iW\).
\end{remark}

\begin{lemma}\label{lem: acceptance domains}
For each \(r \geq 0\),
\[
W = \bigcup_{A \in \sA(r)} \overline{A}
\]
and \(A \cap B = \emptyset\) for distinct \(A\), \(B \in \sA(r)\). For sufficiently large \(r\), for each \(r\)-patch \(P\) and \(x \in \cps\), we have that \(P_r(x)\) and \(P\) are equivalent (up to translation) if and only if \(x^\star \in A_P\). In particular, \(p(r) = \#\sA(r)\).
\end{lemma}

\begin{remark}
We have taken the convention that acceptance domains are open, so only the closures of acceptance domains tile the original window. This has little relevance, since in our non-singular setting \(x^\star\) never hits the boundary of an acceptance domain. Of course, each \(A_P \neq \emptyset\), since otherwise there would be no \(x^\star \in A_P\) for \(x \in \cps\) and thus \(P\) would not occur, contradicting it being an \(r\)-patch.
\end{remark}

\begin{remark}
The acceptance domain containing a point \(w \in W\) (at least for the dense set of \(w \notin \partial W + \Gamma_<\)) is given by
\[
\bigcap_{x \in \Gamma(r)} W(w,x_<),
\]
where \(W(w,x_<) = \iW - x_<\) if \(w \in \iW - x_<\) and \(W(w,x_<) = W^\mathrm{c} - x_<\) otherwise. Since \(p(r) = \# \sA(r)\) for sufficiently large \(r\), this gives an approach to bounding the complexity, by bounding the number of different possible such intersections.
\end{remark}

\subsection{Further properties}

Theorem \ref{thm: main2} simplifies to Theorem \ref{thm: main1} for polytopal windows whose supporting hyperplanes are commensurate with the lattice, formalised as follows:

\begin{definition} \label{def: homogeneous}
We call a polytopal cut and project scheme {\bf weakly homogeneous} if there is some \(N \in \N\) so that, for each \(H \in \sH\), there exists some \(\gamma_H \in \frac{1}{N}\Gamma\) with
\begin{equation} \label{eq: homogeneous}
\bigcap_{H \in \sH} (H - (\gamma_H)_<) \neq \emptyset.
\end{equation}
We call the scheme {\bf homogeneous} if one may take \(N = 1\) in the above.
\end{definition}

\begin{remark} \label{rem: all hyperplanes intersect to point}
To have a compact window, there must be enough hyperplanes so that the intersection of all supporting subspaces in \(\sH_0\) is the origin. Indeed, otherwise, each \(H \in \sH_0\) is parallel to some non-compact subspace \(X < \intl\). Since \(X\) remains parallel with each \(H \in \sH\), which together contain the boundary of \(W\), we would have \(X+w \subset W\) for any \(w \in W\), contradicting \(W\) being compact. So in Equation \eqref{eq: homogeneous} the intersection is a single point belonging to each translated hyperplane. Up to translation, one could assume without loss of generality that this is the origin, although this would make the scheme singular; the above definition is invariant under translations of \(W\).
\end{remark}

\begin{example}
Consider a codimension \(n=1\) cut and project scheme, with internal space \(\intl \cong \R\). Then \(W\) is a (or its partition regions are) finite union(s) of intervals. The supporting hyperplanes are the endpoints of these intervals and the scheme is homogeneous if and only if these all belong to \(\Gamma_<\), perhaps after translating \(W\). The scheme is weakly homogeneous if all vertices belong to some \(\frac{1}{N}\Gamma_<\) (\(N \in \N\)), after an appropriate translate of \(W\).
\end{example}

\begin{example}
For a canonical cut and project set, the window (in a singular position) is the projection of the unit hypercube in \(\tot = \R^k\), with \(\Gamma = \Z^k\). Each supporting hyperplane intersects \(W\) along a face with vertices as elements in \(\Gamma_<\), so these vertices can be used to translate all supporting hyperplanes of \(W\) over the origin. Then the scheme is homogeneous.
\end{example}

To simplify notating going forward, we observe that we may work with unlabelled windows, up to MLD \cite[Section 5.2]{AOI} equivalence:

\begin{proposition} \label{prop: unlabel}
For any polytopal cut and project scheme \(\mathcal{S}\) with labelled window there is another, \(\mathcal{S}'\), which produces MLD equivalent cut and project sets and is equal to \(\mathcal{S}\) except for having a different, unlabelled window. Moreover, the supporting hyperplanes for \(\mathcal{S}'\)  are all \(\Gamma_<\)-translates of those in \(\mathcal{S}\) and vice versa. Thus, in particular, \(\mathcal{S}'\) has the same set of supporting subspaces and is homogeneous, weakly homogeneous or indecomposable (Definition \ref{def: decomposable}) if and only if \(\mathcal{S}\) has those respective properties.
\end{proposition}

\begin{proof}
Take \(\gamma_i \in \Gamma\) so that the translates \(W_i' \coloneqq W_i - (\gamma_i)_<\) do not intersect each other. Consider the scheme \(\mathcal{S}''\) that is equal to \(\mathcal{S}\) except with window \(W' \coloneqq \bigcup_{i=1}^\ell W_i'\), where each \(W_i'\) keeps the label of \(W_i\). Then clearly \(\mathcal{S}''\) defines MLD cut and project sets to \(\mathcal{S}\). Indeed, by sending each point \(x \in \cps\) with label \(i\) to \(x - (\gamma_i)_\vee\), we define a local rule converting a cut and project set \(\cps\) defined by \(\mathcal{S}\) to one from \(\mathcal{S}''\) and analogously for the reverse derivation.

Since the \(W_i'\) are separated we may find \(\epsilon > 0\) so that each \(x \in W_i'\) is at least distance \(\epsilon\) from any other \(W_j'\). Then it is easy to see that the acceptance domains of \(\sA(r)\) for \(\mathcal{S}''\) belong to distinct \(W_i'\) for sufficiently large \(r\), using density of \(\Gamma_<\). So for any \(x \in \Lambda\) we have that \(P = P_r(x)\) determines \(A_P \in \sA(r)\), which can only belong to one partition element of \(W'\), so that \(P\) determines the label of \(x\). This process of adding labels gives an MLD equivalence from cut and project sets of \(\mathcal{S}'\) to \(\mathcal{S}''\) (forgetting labels gives the reverse derivation). By transitivity, \(\mathcal{S}\) and \(\mathcal{S}''\) define MLD cut and project sets.

Since each region \(W_i'\) is a \(\Gamma_<\)-translate of \(W_i\), the sets \(\sH_0\) and \(\sH + \Gamma_<\) are not affected. As indecomposability only depends on \(\sH_0\), and (weak) homogeneity only depends on \(\sH + \Gamma_<\), it follows that the above MLD equivalence preserves these properties.
\end{proof}

\begin{remark}
Many proofs to follow may in principle be converted to ones for labelled cut and project sets, by working with the appropriate definition of the acceptance domains. However, the above result allows us to work with unlabelled schemes from the outset, which simplifies the proofs and notation needed.
\end{remark}

\subsection{Decompositions of cut and project schemes} We finish this section by recalling the notion of a decomposition of the window of a cut and project scheme and some useful associated results. Decompositions in this context were first defined in \cite{FHK02}. An equivalent definition was developed in \cite{KoiWalII}. Most cut and project sets of interest are indecomposable, in which case several of the results and proofs later become simpler, but for our general result they are necessary. The following definition may be found in \cite[Section 3]{KoiWalII}.

\begin{definition} \label{def: decomposable}
For a subset \(S \subseteq \sH_0\), we write
\[
X_S \coloneqq \bigcap_{V \in S} V.
\]
We call a collection \(Q = \{S_1, \ldots, S_\ell\} \) of subsets \(S_i \subset \sH_0\) a {\bf decomposition} of the window if
\[
\intl = X_1 + X_2 + \cdots + X_\ell \text{ for } X_i \coloneqq X_{S_i},
\]
where for every \(H \in \sH_0\) there is exactly one \(S_i \in Q\) with \(H \notin S_i\). If the only decomposition is the trivial one, \(Q = \{\emptyset\}\), then we call the window (and cut and project scheme) {\bf indecomposable} and otherwise it is called {\bf decomposable}.
\end{definition}

We emphasise that each \(S_i\) is a strict subset in the above definition (or else one may include superfluous \(S_i = \sH_0\) elements). Note that, since \(S \subset \sH_0\) rather than \(\sH\) in the definition, the relative positions of the supporting hyperplanes are not relevant, it is their `directions' that are important. The condition that each \(H \in \sH_0\) is in exactly one of \(S_i \in Q\) of course simply means that \(\{S_1^\mathrm{c},\ldots, S_\ell^\mathrm{c}\}\) is a partition of \(\sH_0\). Clearly \(X_i \cap X_j = \{0\}\) for each \(i \neq j\), since for each \(V \in \sH_0\) we have either \(V \in S_i\) or \(V \in S_j\) (or possibly both), so that \(X_i \cap X_j \subseteq \bigcap_{V \in \sH_0} V\), which contains only the origin (see Remark \ref{rem: all hyperplanes intersect to point}). So for a decomposable window, one may split the internal space into two (or more) factors \(X_i\) that are aligned with the supporting hyperplanes, in the sense that each hyperplane is parallel to all but one such factor, for which it intersects to a hyperplane in \(X_i\). Of course, each \(\dim(X_i) > 0\), since for any \(V \notin S_i\) we have \(V \in S_j\), and thus \(X_i \subset V\), for all \(j \neq i\). As noted in \cite[Proposition 3.12]{KoiWalII}, to specify a decomposition it is equivalent to specify the subspaces \(X_i\) that factorise \(\intl\), with the property that each \(V \in \sH_0\) contains all but one \(X_i\).

Given {\Cpx}, an indecomposable scheme has sufficient interaction between the supporting hyperplanes to ensure regularity of \(\Gamma_<\). On the other hand, a decomposable scheme can behave like a sum of cut and project schemes, whose total dimensions sum to \(k\) and codimensions sum to \(n\), with \(\Gamma_<\) splitting into factors lying in subspaces of the internal space with different relative ranks. This is important in defining the correct notion of Diophantine scheme, as we will see in Definition \ref{def: Diophantine scheme}. It is easiest to get a feel for which windows are decomposable by considering some low dimensional examples:

\begin{example}
Every codimension \(n=1\) scheme is indecomposable (there is only one supporting subspace, the origin). In codimension \(n=2\), the scheme is decomposable if and only if there are only two supporting subspaces, meaning that \(W\) is a union of parallelograms, all of which have mutually parallel sides. In higher codimensions the situation is more complicated but, generally, a convex polytopal window is decomposable if and only if it is the Minkowski sum of two lower dimensional polytopes \cite[Theorem 3.19]{KoiWalII} lying in complementary subspaces of dimensions at least \(1\); for example, in codimension \(n=3\) the scheme is decomposable if and only if \(W\) is a prism, given by a product of an interval and a convex \(2\)-dimensional polygon.
\end{example}

\begin{lemma}\label{lem: two element decomposition}
For any decomposition \(Q = \{S_1,\ldots,S_\ell\}\) (with \(\ell \geq 2\)) and any \(S_i \in Q\), we have that \(\{S_i,S_i^\mathrm{c}\}\) is also a decomposition. In particular, a decomposable scheme always has a decomposition of size two.
\end{lemma}

\begin{proof}
Let us denote \(A \coloneqq S_i\), \(B \coloneqq S_i^\mathrm{c}\), \(X_A \coloneqq X_i = \bigcap_{V \in A} V\) and \(X_B \coloneqq \bigcap_{V \in B} V\). Clearly some \(V \in A\) or else \(A = \emptyset\), making each other \(S_j = \sH_0\), contradicting these being strict subsets, so \(B \neq \sH_0\). Since \(B \subseteq S_j\) for each \(j \neq i\) (since any \(V \notin S_i\) must be in all other \(S_j\)) we have \(X_j \subseteq X_B\) and thus \(\intl = (X_1) + (X_2 + \cdots + X_\ell) \subseteq X_A + X_B\) so \(X_A + X_B = \intl\) and \(\{A,B\}\) is a decomposition, as required.
\end{proof}

In the following definition, `flag' is used in the usual sense from linear algebra (rather than the related notion in the context of polytopes), except they would typically be specifically called `complete' or `maximal' flags.

\begin{definition} \label{def: flags}
We call \(f \subseteq \sH_0\) a {\bf flag} if \(f\) consists of precisely \(n = \dim(\intl)\) hyperplanes and
\[
\bigcap_{V \in f} V = \{0\}.
\]
Similarly, we say that \(f \subseteq \sH\) is a flag if \(f\) consists of precisely \(n\) hyperplanes and intersects to a singleton (i.e., when \(V(f) = \{V(H) \mid H \in f\}\) is a flag). We let \(\sF\) denote the set of flags in \(\sH\) and \(\sF_0\) the set of flags in \(\sH_0\).
\end{definition}

A useful equivalent definition of a decomposition, in terms of an associated graph, was given in \cite{KoiWalII} for two-element decompositions. We expand this to arbitrary decompositions in Theorem \ref{thm: graph decomposition} below, whose proof is similar to that of \cite[Theorem 3.7]{KoiWalII}.

\begin{definition} \label{def: graph}
We define the undirected graph \(G(W)\) by setting the vertex set as \(\sH_0\) and connect \(V\) and \(V'\) with a single edge if and only if there is some \(f \subseteq \sH_0\) (with \(n-1\) elements) so that both \(f \cup \{V\}\) and \(f \cup \{V'\}\) are flags.
\end{definition}

\begin{theorem} \label{thm: graph decomposition}
Let \(Q = \{S_1,\ldots,S_\ell\}\) be such that \(\{S_1^\mathrm{c},\ldots,S_\ell^{\mathrm{c}}\}\) is a partition of non-empty subsets of \(\sH_0\). Then \(Q\) is a decomposition if and only if, for every \(S_i \in Q\), \(V \in S_i\) and \(V' \in S_i^\mathrm{c}\), we have that \(V\) and \(V'\) are not connected by an edge in \(G(W)\), that is, \(S_i\) and \(S_i^\mathrm{c}\) are disconnected in \(G(W)\).
\end{theorem}

\begin{proof}
First, assume that \(Q\) is a decomposition and take any \(S_i \in Q\). We denote \(A \coloneqq S_i\) and \(B \coloneqq S_i^\mathrm{c}\), which gives a decomposition \(\{A,B\}\) by Lemma \ref{lem: two element decomposition}. For the proof that \(A\) and \(B\) are disconnected in \(G(W)\), we refer the reader to the proof of \cite[Theorem 3.7]{KoiWalII}.

The converse direction is not directly covered by the proof of the two element case from \cite{KoiWalII}, although can provide a similar argument. Suppose that each pair \(S_i\), \(S_i^\mathrm{c}\) is disconnected from each other in \(G(W)\) and choose any flag \(f = \{V_1,\ldots,V_n\} \in \sF_0\). We denote \(f_i \coloneqq f \cap S_i\). The subspaces \(X_i' \coloneqq \bigcap_{V \in f_i} V\) have sum \(X_1' + \cdots + X_\ell' = \intl\). Indeed, intersecting a subspace with a codimension one subspace decreases its dimension by at most \(1\), so that each \(\dim(X_i) \geq n - \#f_i\) and thus \(\dim(X'_1) + \dim(X'_2) + \cdots + \dim(X'_\ell) \geq (n-\#f_1) + (n-\#f_2) + \cdots + (n-\#f_\ell) = n^2 - n(n-1) = n\), where we note that each \(V \in f\) appears in all but one \(f_i\). But these subspaces are also complementary, in the sense that for any \(V \in f\) we either have \(X_i' \subseteq V\) or \(X_j' \subseteq V\) for every other \(X_j'\). So \(X_i' \cap X \subseteq V\) for all \(V \in f\) and any \(X\) that is a sum of \(X_j'\) for \(j \neq i\). It follows that \(\dim(X_1' + \cdots + X_n') = n\).

We claim that each \(X_i = X_i'\), from which the result follows by the above. Clearly each \(X_i \subseteq X_i'\), since \(S_i \supseteq f_i\). To see that \(X_i \supseteq X_i'\), take any \(V \in S_i\). We wish to show that \(X_i' \subseteq V\), from which the claim follows (since \(X_i\) is the largest subspace containing each \(V \in S_i\)). By definition, if \(V \in f\) then \(V \in f_i\) so that \(X_i' \subseteq V\). Otherwise, let us denote \(A \coloneqq S_i\) and \(B \coloneqq S_i^\mathrm{c}\) and consider the \(n+1\) element set \(f' = f \cup \{V\}\). We list its elements as
\[
f' = \{V_1,V_2,\ldots,V_m,V,V_{m+1},\ldots,V_n\}
\]
where we have ordered the elements of \(f\) so that \(V_1\), \ldots, \(V_m \in A\) appear first, whilst each \(V_{m+1}\), \ldots, \(V_n \in B\). The intersections \(V_1\), \(V_1 \cap V_2\), \(V_1 \cap V_2 \cap V_3\), \ldots either stays the same (if the new hyperplane contains the previous intersection) or drops by \(1\) at each step, since each \(V_j\) and also \(V\) are codimension \(1\). Since there are \(n+1\) elements and the full intersection is the origin, there is exactly one step where the dimension does not drop, and since \(f\) is a flag the dimension drops for each of the first \(m\) intersections. If \(X_i' \nsubseteq V\) then the dimension must fail to drop later, so that \(V_1 \cap \cdots \cap V' = V_1 \cap \cdots \cap V' \cap V_j\), where \(V_j \in B\) and \(V' = V\) if \(j = m+1\) or else \(V' = V_{j-1} \in B\). In either case, we can omit \(V_j\) from the list to obtain a new flag \(f' \setminus \{V_j\}\) with \(V_j \in B\). But \(f\) and \(f'\) differ only by replacing \(V_j \in f\) with \(V \in f'\), so that \(V \in A = S_i\) and \(V_j \in S_i^\mathrm{c}\) are connected by an edge in \(G(W)\), contradicting our assumption on \(G(W)\).
\end{proof}

We conclude this section by noting that decompositions can be made `coarser' by combining elements (generalising Lemma \ref{lem: two element decomposition}) and two different decompositions can be combined to give a decomposition that is `finer' than both. Versions of the results below appear in \cite{KoiWalII} but we provide different, simpler proofs here, using the graph version of the definition of a decomposition provided by Theorem \ref{thm: graph decomposition}.

\begin{lemma}[{\cite[Proposition 3.14]{KoiWalII}}]
Given a decomposition \(Q = \{S_1,\ldots,S_\ell\}\) with \(\ell \geq 2\) elements, define \(Q'\) by replacing any two distinct \(S_i\), \(S_j \in Q\) with \(S_i \cap S_j\). Then \(Q'\) is also a decomposition.
\end{lemma}

\begin{proof}
Clearly the new element \(S_i \cap S_j \neq \sH_0\), since \(S_i\), \(S_j \neq \sH_0\). Given any \(V \in \sH_0\), either \(V \notin S_m\) for some \(m \neq i,j\), in which case \(V \in S_i \cap S_j\) and \(V\) is an all elements of \(Q'\) except \(S_m\), or else \(V \notin S_i \cap S_j\) but is in every other \(S_m \in Q'\), as required for a decomposition. For any \(V \in S_i \cap S_j\) and any \(V' \in (S_i \cap S_j)^\mathrm{c}\) we have that either \(V' \in S_i^\mathrm{c}\) or \(V' \in S_j^\mathrm{c}\). In either case, \(V\) and \(V'\) are not connected by an edge in \(G(W)\) since \(V \in S_i\), \(V \in S_j\) and \(Q\) is a decomposition, so we have verified that \(Q'\) is a decomposition by Theorem \ref{thm: graph decomposition}.
\end{proof}

\begin{lemma}
Suppose that \(Q = \{S_1,\ldots,S_\ell\}\) and \(Q' = \{S_1',\ldots S_m'\}\) are two decompositions. Then there is a decomposition \(Q''\) for which, for all \(S \in Q''\), we have that \(S \supseteq S_i\) and \(S \subseteq S_j'\) for some \(S_i \in Q\) and \(S_j \in Q'\). This decomposition is the coarsest one refining both \(Q\) and \(Q'\), in the sense that for any \(S_i \in Q\) and \(S_j' \in Q'\), with associated subspaces \(X_i = \bigcup_{V \in S_i} V\) and \(X_j' = \bigcup_{V \in S_j} V\) for which \(X_i \cap X_j' \neq \{0\}\), there is some \(S \in Q''\) so that \(X \coloneqq \bigcap_{V \in S} V\) satisfies \(X = X_i \cap X'_j\).
\end{lemma}

\begin{proof}
We define
\[
Q'' \coloneqq \{S_i \cup S_j' \mid S_i \in Q, S_j' \in Q'\} \setminus \{\sH_0\}
\]
i.e., we take arbitrary unions of elements from \(Q\) and \(Q'\) and discard any which are all of \(\sH_0\), so that \(\sH_0 \notin Q''\) by definition. We still have that each \(V \in \sH_0\) is in all but one element of \(Q''\), namely the unique \(S_i \cup S_j'\) for which \(V \notin S_i\) and \(V \notin S_j'\). Take any \(S_i \cup S_j' \in Q''\), and any \(V \in S_i \cup S_j'\), \(V' \in (S_i \cup S_j')^\mathrm{c}\). If \(V \in S_i\) then \(V\) and \(V'\) are not connected by an edge in \(G(W)\) since \(V' \in S_i^\mathrm{c}\), and analogously if \(V \in S_j'\). So we have verified that \(Q''\) is a decomposition by Theorem \ref{thm: graph decomposition}.

Suppose that \(X_i \cap X_j' \neq \{0\}\). Then there exists some \(V \in (S_i \cup S_j')^\mathrm{c}\) since, otherwise, \(X_i \cap X_j' = \left(\bigcap_{V \in S_i} V\right) \cap  \left(\bigcap_{V \in S_j'} V \right) = \bigcap_{V \in \sH_0} V = \{0\}\). It follows that \(S_i \cup S_j' \neq \sH_0\) and thus \(S_i \cup S_j' \in Q''\), as required, since \(X = X_i \cap X_j' = \bigcap_{V \in S_i \cup S_j'} V\).
\end{proof}

The result above shows that there is always a `finest' decomposition, in the sense that it defines subspaces \(X_i\) which are contained in those defined by any other decomposition. This also follows easily from Theorem \ref{thm: graph decomposition}, as this decomposition will be the one consisting of the subsets \(S^\mathrm{c} \subseteq \sH_0\) for which \(S\) is a connected component of \(G(W)\).

\begin{figure}
  \begin{subfigure}{.45\textwidth}
  \renewcommand\captionlabelfont{}
    \centering\def\svgwidth{\textwidth}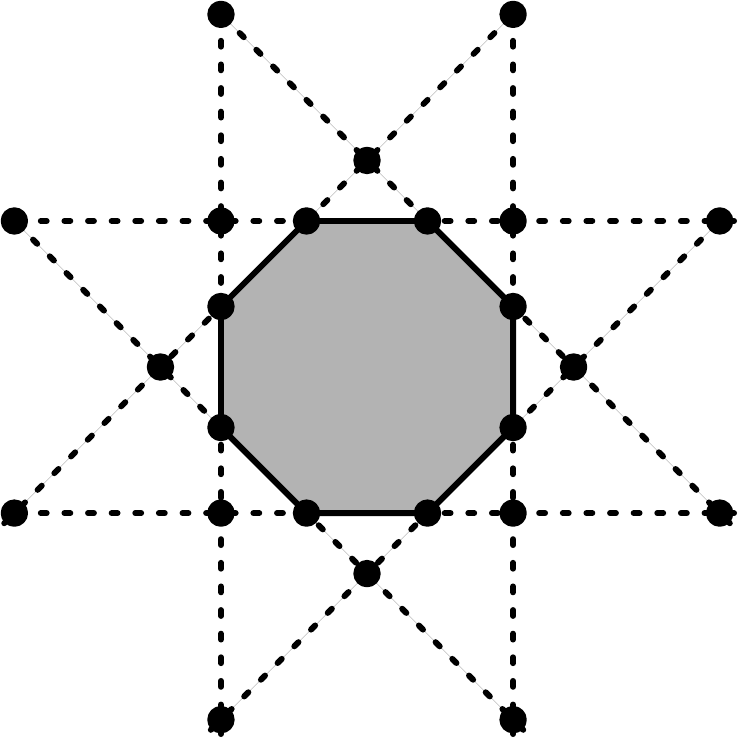%
    \caption{Generalised vertices for the Ammann--Beenker scheme. Consideration of these is not important in applying our results, as the scheme is homogeneous.}
    \label{fig: AB generalised vertices}
  \end{subfigure}\hfill
  \begin{subfigure}{.45\textwidth}
  \renewcommand\captionlabelfont{}
    \centering\def\svgwidth{\textwidth}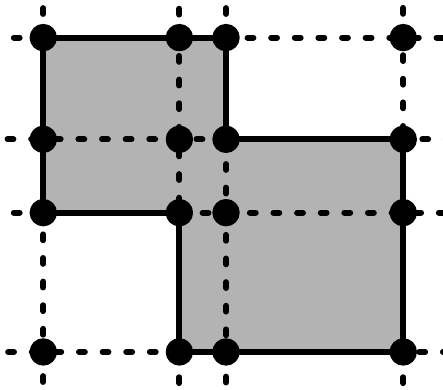%
    \caption{Generalised vertices for a decomposable window. Although the window is not a product of \(1\)-dimensional windows, note that the generalised vertex set is a product with respect to the decomposition of internal space.}
    \label{fig: decomposable}
  \end{subfigure}
  \caption{Generalised vertices in codimension \(n=2\).}\label{fig: generalised vertices}
\end{figure}

\section{Complexity} \label{sec: complexity}
In this section we present a formula for the asymptotic polynomial degree of the complexity function of a polytopal cut and project set. This is the same formula as given in \cite[Theorem 6.1]{KoiWalI} but with hypotheses weakened to allow for non-convex and labelled polytopal windows. Afterwards, we outline some important consequences of having low complexity that follow from this result.

\subsection{Stabiliser subgroups and growth of the complexity function}
Recall the definition of a flag from Definition \ref{def: flags}. We use these to define the generalised vertices, see Figure \ref{fig: generalised vertices}.

\begin{definition} \label{def: vertices}
For a flag \(f \in \sF\), let \(v(f)\) denote the unique intersection point of the flag and
\[
v(\sH) \coloneqq \{v(f) \mid f \in \sF\}
\]
denote the {\bf generalised vertices} of the window. Displacements between generalised vertices are denoted by
\[
F \coloneqq \{v(f) - v(f') \mid f, f' \in \sF\}.
\]
\end{definition}

\begin{definition}
The {\bf stabiliser} of a supporting hyperplane \(H \in \sH\) or subspace \(H \in \sH_0\) is
\[
\Gamma^H \coloneqq \{\gamma \in \Gamma \mid H = H + \gamma_<\}.
\]
Then \(\Gamma^H = \Gamma^{V(H)} = \{\gamma \in \Gamma \mid \gamma_< \in V(H)\}\). For a supporting hyperplane (or subspace) \(H\) we define its {\bf rank} to be \(\rk(H) \coloneqq \rk(\Gamma^H)\), where the latter is the rank of \(\Gamma^H\) considered as a finitely generated free Abelian group.

More generally, for a collection \(S\) of supporting hyperplanes or supporting subspaces, we let
\[
\Gamma^S \coloneqq \bigcap_{V \in S} \Gamma^V.
\]
\end{definition}

\begin{remark}
It is not hard to see that \(\Gamma^S\) is the subgroup of \(\Gamma\) of elements stabilising each hyperplane \(H \in S\), or equivalently the subgroup of \(\Gamma\) which projects into the subspace \(\bigcap_{H \in S} V(H)\) of the internal space. Moreover, the quotient groups \(\Gamma / \Gamma^S\) are all free Abelian. Indeed, we may identify this quotient with its image under projection to \(\intl / \bigcap_{H \in S} V(H)\), which is torsion-free. In the case of a single hyperplane, we see that \(\Gamma / \Gamma^H\) may be identified with the subgroup of \(\intl / V(H) \cong \R\) of distances of points in \(\Gamma_<\) from \(V(H)\).
\end{remark}

The supporting hyperplanes being well-aligned with the lattice, in the sense that the stabilisers have high rank, is a special property of a cut and project scheme. As we will see, this corresponds to the cut and project sets having low complexity. Since the acceptance domains may have shapes that are difficult to control, we introduce the cut regions:

\begin{definition} \label{def: cut regions}
For \(r \geq 0\), an {\bf \(r\)-cut region} is any connected component of
\[
W \setminus \bigcup_{H \in \sH} \left( \bigcup_{\gamma \in \Gamma(r)} (H-\gamma_<) \right).
\]
We let \(\sC(r)\) denote the set of \(r\)-cut regions.
\end{definition}

In other words, we remove all \(\Gamma(r)_<\)-translates of supporting hyperplanes from the internal space, which `cuts' the window into convex polytopal regions. Note that the supporting hyperplanes contain the boundary of the window, so it is clear that the cut regions refine the acceptance domains and thus \(p(r) = \#\sA(r) \leq \#\sC(r)\). One needs further properties for using \(\#\sC(r)\) as an effective lower bound for \(p(r)\), such as the \emph{quasicanonical} property \cite{KoiWalI} for convex polytopal windows. See also Lemma \ref{lem: refining acceptance domains}.

\begin{definition} \label{def: box}
We call \(B \neq \emptyset\) a {\bf box} for a flag \(f\) if \(B\) is an intersection of a finite number of \(\Gamma_<\)-translates of \(\iW\) and \(W^\mathrm{c}\) for which, for each \(H \in f\), there are vectors \(\gamma^0_H\), \(\gamma^1_H \in \Gamma\) with
\begin{equation} \label{eq: box}
B = \left(\bigcap_{H \in f} H^+ - (\gamma^0_H)_<\right) \cap \left(\bigcap_{H \in f} H^- - (\gamma^1_H)_<\right),
\end{equation}
where \(H^+\) is one open half space of \(\intl \setminus H\) and \(H^-\) is the other open half space.
\end{definition}

It follows from density of \(\Gamma_<\) that, for every \(f \in \sF\), we may construct a box of arbitrarily small diameter. It may also be chosen to share a face with any given \(H' \in f\), on either the \((H')^+\) or \((H')^-\) side. Indeed, in a certain neighbourhood of each \(H\), we may choose small \(\gamma_<\) so that the intersections \(\iW \cap (W^\mathrm{c} - (\gamma_H^1)_<)\) look like strips, parallel to \(H\), at least in the neighbourhood \(U_H\) in Definition \ref{def: polytopal}. These can be translated by elements of \(\Gamma_<\) to intersect to a box, at least when restricted to some open neighbourhood \(U\) of the box (and we can translate to \(U_{H'}\) if we want the box to share a face with \(H'\)). Clearly we can also create intersections of \(\iW - \gamma_<\) with non-empty interior but arbitrarily small diameter; \(\Gamma_<\)-translates of the complements of this then easily removes any extra pieces not part of the box (outside of \(U\)) that may have resulted from non-convexity of the window.

The following theorem is the main result of this section, an extension of \cite[Theorem 6.1]{KoiWalI} to possibly non-convex and labelled polytopal cut and project sets. Many parts of the proof remain the same, so in places we only provide a sketch and refer the reader to \cite{KoiWalI} for further details.

\begin{theorem} \label{thm: generalised complexity} For a polytopal cut and project scheme, \(p(r) \asymp r^\alpha\) where
\[
\alpha \coloneqq \max_{f \in \sF_0} \alpha_f, \ \alpha_f \coloneqq \sum_{V \in f} \alpha_V, \ \alpha_H \coloneqq d - \rk(V) + \beta_V \text{ and } \beta_V \coloneqq \dim \langle \Gamma_<^V \rangle_\R.
\]
\end{theorem}

\begin{proof}
We may assume without loss of generality that \(W\) is unlabelled. The result for labelled cut and project schemes then follows from the unlabelled case and Proposition \ref{prop: unlabel}, since for two complexity functions \(p\) and \(p'\) of MLD patterns, we have the bounds \(p(r) \leq p'(r+c)\) and \(p'(r) \leq p(r+c)\) for all \(r \geq 0\) and constant \(c\), where the sets \(\sH_0\), \(\sF_0\) and numbers \(\rk(H)\), \(\beta_H\) are not affected.

\emph{Upper bound}: For this proof only, we make use of the \emph{modified} \(r\)-cut regions \(\sC'(r)\), defined just as in Definition \ref{def: cut regions} except that we only cut with hyperplanes \(H - \gamma_<\) for \(\gamma \in \Gamma\) with \(\|\gamma_\vee\| \leq r\) and \(\gamma_< \in W-W\). It is still clear (see Remark \ref{rem: origin of cut regions}) that, for sufficiently large \(r\), these modified \(r\)-cut regions refine the acceptance domains, that is, for each \(A \in \sA(r)\) there is some \(C \in \sC'(r)\) with \(C \subseteq A\), so \(p(r) \leq \#\sC'(r)\). The proof is now analogous to that of \cite{KoiWalI}, which we outline here.

Each \(C \in \sC'(r)\) is the convex hull of its vertices, where \(v\) is a vertex if, for a flag \(f \in \sF\),
\begin{equation} \label{eq: vertex from flag}
\{v\} = \bigcap_{H \in f} H - (\gamma_H)_<,
\end{equation}
where \(v \in \partial C\), each \(H - (\gamma_H)_<\) intersects \(\partial C\) along a facet and the \(\gamma_H \in \Gamma\) with \(\|(\gamma_H)_\vee\| \leq r\) and \((\gamma_H)_< \in W-W\). Each vertex is shared by a bounded number of such regions, with bound depending only on \(\sH_0\), so we just need to bound the number of vertices.

Set a flag \(f \in \sF\). We have equal vertices \(\bigcap_{H \in f} H - (\gamma_H)_< = \bigcap_{H \in f} H - (\gamma_H')_<\) if and only if each \(H - (\gamma_H)_< = H - (\gamma_H')_<\) (distinct but parallel translated hyperplanes have trivial intersection). However, distinct choices of the \(\gamma_H\), \(\gamma_H'\) can give the same translated hyperplanes, precisely those with \(\gamma_H - \gamma_H' \in \Gamma^H\).

For a given \(H \in f\), choose a basis \(b\) for \(\intl\) of elements from \(\Gamma_<\), where as many as possible (namely, \(\beta_H\)) are chosen from \(\Gamma^H\). Then \(G \coloneqq \langle b \rangle_\Z\) is an \(n\)-dimensional lattice in \(\intl\). Consider the subgroup \(\Gamma^H + G\). This has rank \(\rk(H) + (k-d) -\beta_H\), since \(\Gamma^H\) and \(G\) intersect just on the rank \(\beta_H\) subgroup of elements in \(G\) chosen from \(\Gamma^H\). 

For any \(H \in \sH\) and \([\gamma] \in \Gamma/(\Gamma^H+G)\), there are a uniformly bounded number of representatives \(\alpha \in [\gamma]\) defining distinct cuts \(H - \alpha_<\) of modified cut regions. Indeed, if two representatives differ by an element in \(\Gamma^H\), then they give the same translated hyperplane \(H - \alpha_<\). And for any such choice there are only a bounded number of \(g \in G\) with \((\alpha + g)_< \in W-W\).

The quotient \(\Gamma/(\Gamma^H+G)\) has rank \(k - (\rk(H) + (k-d) -\beta_H) = d - \rk(H) + \beta_H = \alpha_H\), so there are \(\ll r^{\alpha_H}\) distinct translates \(H - (\gamma_H)_<\) to choose from in the intersection of Equation \eqref{eq: vertex from flag}. Thus, there are \(\ll \prod_{H \in f} r^{\alpha_H} = r^{\sum_{H \in f} \alpha_H} = r^{\alpha_f}\) vertices that can come from each flag \(f \in \sF\) (which corresponds to a flag in \(\sF_0\) with the same associated stabiliser ranks/spans). Since every vertex \(v\) is produced by some flag (possibly more, if more than \(n\) hyperplanes happen to intersect at \(v\)), it follows that there are \(\ll \sum_{f \in \sF_0} r^{\alpha_f} \ll r^\alpha\) vertices, as required.

\emph{Lower bound:} Again, potential non-convexity of \(W\) does not affect the main ideas of the proof from \cite{KoiWalI}. We need to show that, for each flag \(f \subseteq \sH\), there are \(\gg r^{\alpha_f}\) \(r\)-acceptance domains in \(\sA(r)\).

We begin by constructing a box \(B \subset W\) (Definition \ref{def: box}) from the flag \(f\), which is `small' relative to the open sets \(U_H\) in Definition \ref{def: polytopal}. Fix some \(H \in f\). We consider intersections of \(B\) with further translates \(\iW - \gamma_<\) and \(W^\mathrm{c} - \gamma_<\) for \(\gamma \in \Gamma(r)\), which are `nice' in the sense that the intersection with \(B\) is identical with a non-trivial intersection of a corresponding translated half-space \(H^\pm - \gamma_<\). By taking \(B\) sufficiently small relative to the open neighbourhoods \(U_H\) of Definition \ref{def: polytopal}, we can ensure a cut is nice by restricting to those \(\gamma\) for which \(w_H - \gamma_<\) is sufficiently close to the centre of \(B\), where \(w_H\) is as defined in Definition \ref{def: polytopal}; this only restricts \(\gamma_<\) to an open set \(U\) (analogously to the bounded set \(W-W\) from the reverse bound). Again, translates differing by an element of \(\Gamma^H\) determine the same cut, but otherwise are distinct. The calculation of how many such cuts are made now proceeds similarly to the lower bound: summarising, the rank \(\rk(H)\) subgroup \(\Gamma^H\) is lost to giving identical cuts, a rank \(k-d\) subgroup \(G < \Gamma\) (giving a lattice \(G_<\) in \(\intl\)) is constrained to position \(\gamma_< \in U\), and \(\Gamma^H \cap G\) has rank \(\beta_H\). We refer the reader to \cite[Section 6.3]{KoiWalI} for further details.

The above shows how we may cut the box \(B\) into \(\gg \prod_{H \in f} r^{\alpha_H} = r^{\alpha_f}\) sub-boxes, using intersections of translates of \(\iW - \gamma_<\) and \(W^\mathrm{c} - \gamma_<\) with \(\gamma \in \Gamma(r)\). The acceptance domains of \(\sA(r)\) are given by such intersections, potentially refining these further, and thus there are \(\gg \max_{f \in \sF_0} r^{\alpha_f} = r^\alpha\) \(r\)-acceptance domains, as required.
\end{proof}

\subsection{Constant stabiliser rank and consequences of low complexity} In the remainder of this section we show that property {\Cpx} induces a certain splitting of the scheme, which is required to correctly define our Diophantine condition in Section \ref{sec: Diophantine schemes}.

\begin{definition}
We say that \(\mathcal{S}\) has {\bf constant stabiliser rank} if \(\rk(V_1) = \rk(V_2)\) for all \(V_1\), \(V_2 \in \sH_0\). 
\end{definition}

\begin{definition}
We say that \(\mathcal{S}\) is {\bf hyperplane spanning} if \(\langle \Gamma_<^V \rangle_\R =V\) for all \(V \in \sH_0\). That is, \(\beta_V = n-1\) in Theorem \ref{thm: generalised complexity}.
\end{definition}

We now collect some useful results from \cite{KoiWalII} on stabiliser ranks, indecomposability and complexity. Note that the properties of constant stabiliser rank, hyperplane spanning, indecomposability and the data required for Theorem \ref{thm: generalised complexity} only depend on \((\sH_0,\Gamma_<)\). Similarly, the omitted proofs of Lemmas \ref{lem: minimum complexity} and \ref{lem: consequences of C} below only depend on this data, so they proceed identically to those in \cite{KoiWalII} and therefore we omit the more technical proofs here.

\begin{lemma}[{\cite[Corollary 2.5]{KoiWalII}}] \label{lem: minimum complexity}
Each \(\alpha_f \geq d\) in Theorem \ref{thm: generalised complexity}, so that \(p(r) \gg r^d\) for any polytopal cut and project scheme.
\end{lemma}

\begin{lemma}[{\cite[Theorem 2.8]{KoiWalII}}] \label{lem: consequences of C}
If \(\mathcal{S}\) has property {\Cpx} then:
\begin{enumerate}
	\item each \(\alpha_f = d\);
	\item \(\mathcal{S}\) is hyperplane spanning;
	\item for a flag \(f = \{V_1,\ldots,V_n\} \in \sF_0\) and \(n \in \{1,\ldots,n\}\), let \(\hat{m} \coloneqq f \setminus \{V_m\}\). Then the \(\Gamma^{\hat{m}}_<\) are contained in complementary one-dimensional subspaces of \(\intl\), and \(\Gamma^{\hat{1}} + \cdots + \Gamma^{\hat{n}}\) is a finite index subgroup of \(\Gamma\).
\end{enumerate}
\end{lemma}

\begin{corollary} \label{cor: sum of ranks for low complexity}
The scheme \(\mathcal{S}\) has property {\Cpx} if and only if \(\sum_{V \in f} \rk(V) = k(n-1)\) for each flag \(f \in \sF_0\).
\end{corollary}

\begin{proof}
Suppose that {\Cpx} holds. By Lemma \ref{lem: consequences of C}, \(\alpha_f = d\) for each \(f \in \sF_0\) in Theorem \ref{thm: generalised complexity}, and each \(\beta_V = n-1\) for each \(V \in f\) by hyperplane spanning, so
\[
\alpha_f = d = \sum_{V \in f} (d - \rk(V) + n-1) =  n(k-1) - \sum_{V \in f} \rk(V) ,
\]
and thus \(\sum_{V \in f} \rk(V) = k(n-1)\), for each \(f \in \sF_0\), as required. Conversely, if \(\sum_{V \in f} \rk(V) = k(n-1)\) holds for each \(f \in \sF_0\) then each
\[
\alpha_f = \sum_{V \in f} (d - \rk(V) + \beta_V) \leq \sum_{V \in f} (d - \rk(V) + (n-1)) = n(k-1) - k(n-1) = k-n = d
\]
and thus \(\alpha \leq d\) in Theorem \ref{thm: generalised complexity} so {\Cpx} holds.
\end{proof}

In the Lemma \ref{lem: consequences of C}, (1) follows from Lemma \ref{lem: minimum complexity} whilst (2) follows quickly from (3). In fact, although not stated in \cite{KoiWalII}, a stronger version of this property follows from it:

\begin{corollary} \label{cor: densely hyperplane spanning}
If \(\mathcal{S}\) has property {\Cpx} then it is {\bf densely hyperplane spanning}, that is, each \(\Gamma_<^H\) is dense in \(V(H)\).
\end{corollary}

\begin{proof}
For any \(V \in \sH_0\), choose some flag \(f \in \sF_0\) containing \(V\). By Conclusion (2) of Lemma \ref{lem: consequences of C}, there are \(n-1\) independent lines in \(V\) containing subgroups \(\Gamma_<^{\hat{m}} \leqslant \Gamma_<^H\). Each \(\Gamma_<^{\hat{m}}\) is dense in the line it spans, or else \(\Gamma_<\) would not be dense in \(\intl\). It follows that \(\Gamma_<^V\) is dense in \(V\).
\end{proof}

The above has the following useful consequence, whose proof is similar to that of \cite[Proposition 5.2]{KoiWalII} using the quasicanonical property.

\begin{lemma}\label{lem: refining acceptance domains}
If \(\mathcal{S}\) is densely hyperplane spanning (in particular, if {\Cpx} holds), then there is some \(\kappa > 0\) for which \(\sA(\kappa r)\) refines \(\sC(r)\), that is, for all \(A \in \sA(\kappa r)\) we have \(A \subseteq C\) for some \(C \in \sC(r)\).
\end{lemma}

\begin{proof}
For each hyperplane \(H \in \sH\) make a small box \(B_H \subset W\) with \(B_H \subset U_H\) (with \(U_H\) as in Definition \ref{def: polytopal}) sharing a face with \(H\). By the dense hyperplane spanning property, we can find a basis \(b_H \in \Gamma_<\) of \(V(H)\) that is dense, relative to the size of this box, in the sense that \(B_H + \langle b_H \rangle_\Z\) exactly covers a half-space \(H^+\) of \(H\) in a neighbourhood of \(H\). By adding further elements \(\gamma_H^i \in \Gamma\), we can use translates \(B_H + \langle b_H \rangle_\Z + (\gamma_H^i)_<\) whose union agrees with the half-space within distance larger than the diameter of \(W\) from \(H\). There are only finitely many elements \(\gamma_H^i\) and lattice elements used to make the boxes \(B_H\), so all can be uniformly bounded with, say, norm \(c\).

Take any cut region \(C \in \sC(r)\). Then it is an intersection of (open) half spaces \(H^+ - \gamma_<\) or \(H^- - \gamma_<\), where \(\gamma \in \Gamma(r)\). By the above, any \(w \in W\) is in \(H^+ - \gamma_<\) if and only if \(w \in B_H + b + (\gamma_H^i)_< - \gamma_<\) for some \(b \in \langle b_H \rangle_\Z\) and \(\gamma_H^i\). Then \(b - \gamma_< \in w - B_H - (\gamma_H^i)_< \subset W-W - (\gamma_H^i)_<\) so that \(\|b-\gamma_<\|\) is uniformly bounded. Moreover, since \(b\) is a basis and \(\|\gamma_<\| \leq r\), we have that \(b \in \Gamma(r')_<\) with \(r' \ll r\). Similarly, \(w \notin H^+ - \gamma_<\) if \(w\) is not in any of these boxes (over \(b \in \Gamma(r')_<\)). Repeating over each hyperplane, we see that \(C\) may be expressed as a union of boxes \(B_H + b + (\gamma_H^i)_< - \gamma_<\), intersected with \(W\) and intersected with the complements of other boxes of the same form. These sets are defined via translates \(g_< \in \Gamma_<\) of \(\iW\) and \(W^\mathrm{c}\) for \(\|g\| \leq r'' = c + r' + r + c \ll r\). It follows that any \(r''\)-acceptance domain \(A\) intersecting \(C\) non-trivially is contained in such a set (possibly refined even further), so that \(A \subseteq C\). As every \(r''\)-acceptance domain is wholly contained in the cut region it intersects with, the result follows.
\end{proof}

For the result below, recall the definition of the graph \(G(W)\) from Definition \ref{def: graph}. The proof of the below is essentially identical to that of \cite{KoiWalII}, so we omit it, but note that it quickly follows from Lemma \ref{lem: consequences of C} and the generalised complexity formula of Theorem \ref{thm: generalised complexity}.

\begin{theorem}[{\cite[Theorem 3.8]{KoiWalII}}] \label{thm: connected => same rank}
Suppose that \(\mathcal{S}\) satisfies {\Cpx} and that \(V_1\) and \(V_2\) are connected by an edge in \(G(W)\). Then \(\rk(V_1) = \rk(V_2)\).
\end{theorem}

The result below follows immediately from the above, along with the fact that \(G(W)\) is a connected graph in the indecomposable case by Theorem \ref{thm: graph decomposition}.

\begin{corollary}[{\cite[Corollary 3.9]{KoiWalII}}] \label{cor: constant stabiliser rank}
If \(\mathcal{S}\) is indecomposable and satisfies {\Cpx} then it is constant stabiliser rank.
\end{corollary}

The result below is essentially \cite[Corollary 3.24]{KoiWalII}, except for omitting statements about writing the window as a Minkowski sum (which does not hold in our current more general setting). As for the other results above (from Lemma \ref{lem: minimum complexity} to here), the proof now follows analogously to \cite{KoiWalII} due to the generalised complexity result of Theorem \ref{thm: generalised complexity} being the same, and we could refer the reader to the proof there. However, we are able to present a more economical argument using Theorem \ref{thm: graph decomposition} and give an explicit description of a decomposition satisfying the result.

\begin{proposition} \label{prop: factorisation to CSR}
Suppose that \(\mathcal{S}\) has property {\Cpx}. For any decomposition \(Q = \{S_1,\ldots,S_\ell\}\), we have that \(\Gamma_1 + \cdots + \Gamma_\ell\) is finite index in \(\Gamma\), where we define \(\Gamma_i \coloneqq \Gamma^{S_i}\). Moreover, there exists a decomposition \(Q\) for which each factor is constant stabiliser rank, in the following sense: given \(S_i \in Q\), the restricted stabilisers
	\[
	\Gamma^V \cap \Gamma_i = \{\gamma \in \Gamma_i \mid \gamma_< \in V\}
	\]
	have constant rank over all \(V \notin S_i\) (of course, by definition, \(\Gamma_i \leqslant \Gamma^V\) for \(V \in S_i\)). Moreover, this common rank is equal to
	\[
	r_i = \rk(\Gamma^V \cap \Gamma_i) = k_i - \delta_i - 1, \text{ where } \delta_i \coloneqq \frac{d_i}{n_i},
	\]
	for \(k_i \coloneqq \rk(\Gamma_i)\), \(n_i \coloneqq \dim(X_i)\) and \(d_i \coloneqq k_i - n_i\). Therefore, \(\rk(V) = k - \delta_i - 1\) depends only on the unique value of \(i\) with \(V \notin S_i\).
\end{proposition}

\begin{proof}
If \(Q = \{\emptyset\}\) is the trivial decomposition then the first part in the above result is trivial (we have that \(\Gamma = \Gamma_1\)) and, if the cut and project scheme is constant stabiliser rank, then \(Q\) satisfies the second part of the result. Indeed, each \(\rk(V)= k - \frac{d}{n} - 1\), by Theorem \ref{thm: generalised complexity} and hyperplane spanning from Lemma \ref{lem: consequences of C}.

More generally, for any given decomposition \(Q = \{S_1,\ldots,S_\ell\}\), take an arbitrary flag \(f = \{V_1,\ldots,V_n\} \in \sF_0\). By Lemma \ref{lem: consequences of C}, this defines subgroups \(\Gamma^{\hat{1}}\), \(\Gamma^{\hat{2}}\), \ldots, \(\Gamma^{\hat{n}} < \Gamma\) whose sum is finite index in \(\Gamma\) and with each \(\Gamma_<^{\hat{m}} \leqslant Y^{\hat{m}} \coloneqq \bigcap_{i \neq m} V_i\). Defining \(X_i'\) to be the intersection of \(V \in f \cap S_i\) (rather than over all of \(S_i\) as for \(X_i\)), it is easily seen that each \(X_i = X_i'\), since trivially \(X_i \subseteq X_i'\) and if this inclusion was strict then \(\dim(X_i' \cap X_j') > 0\) for some distinct \(i \neq j\), which cannot happen (since \(f = S_i \cup S_j\) and has hyperplanes intersecting to the origin). Thus, \((\Gamma_i)_< = \Gamma_< \cap X_i = \Gamma_< \cap X_i' \geqslant \Gamma_< \cap Y^{\hat{m}} = \Gamma_<^{\hat{m}}\) whenever \(V_m \notin S_i\) (since then \(S_i \subseteq \hat{m}\) so that \(X_i \geqslant Y^{\hat{m}}\)). That is, \(\Gamma_i \geqslant \Gamma^{\hat{m}}\). Since each \(\Gamma^{\hat{m}}\) is contained in some \(\Gamma_i\), and the latter have sum that is finite index in \(\Gamma\), we see that \(\Gamma_1 + \ldots + \Gamma_\ell\) is finite index in \(\Gamma\), as required.

We now show that there is a sufficiently fine decomposition to give the constant (restricted) stabiliser rank property. List the distinct possible ranks \(\rk(V)\) over all \(V \in \sH_0\) by \(R_1 < R_2 < \cdots < R_\ell\). Define \(S_i \coloneqq \{V \in \sH_0 \mid \rk(V) \neq R_i\}\). Then it follows immediately from Theorems \ref{thm: graph decomposition} and \ref{thm: connected => same rank} that \(Q = \{S_1,\ldots,S_\ell\}\) is a decomposition, since elements of \(S_i\) and \(S_i^\mathrm{c}\) necessarily have different rank and are thus disconnected in \(G(W)\).

Set a flag \(f \in \sF_0\) and associated notation as above, and take any \(V = V_a \in f\) with \(V \notin S_i\), equivalently \(\rk(V) = R_i\). For each \(m \neq a\) we have \(\Gamma_<^{\hat{m}} \leqslant V\), whilst \(V \cap Y^{\hat{a}} = \{0\}\). Since the \(\Gamma^{\hat{m}}\) have sum that is finite index in \(\Gamma\), it follows that \(\sum_{m \neq a} \Gamma^{\hat{a}}\) is finite index in \(\Gamma^V\). In particular, \(R_i = \rk(V) = k - \rk(\Gamma^{\hat{a}})\) so \(\rk(\Gamma^{\hat{a}}) = k - R_i\) only depends on the index \(i\) for which \(\rk(V_a) = R_i\). Denoting \(\alpha = \rk(\Gamma^{\hat{a}})\), we have that \(k_i = \alpha n_i\), since \(\Gamma_i\) is a sum of \(n_i = \dim(X_i)\) different groups \(\Gamma^{\hat{m}}\). By an identical argument to above (using that the sum of \(\Gamma^{\hat{m}}\) for \(\rk(V_m) = R_i\) are finite index in \(\Gamma_i\)), we have \(\Gamma^V \cap \Gamma_i = \rk(\Gamma_i) - \rk(\Gamma^{\hat{a}}) = k_i - \alpha = k_i - \frac{k_i}{n_i} = k_i - \frac{d_i+n_i}{n_i} = k_i - \frac{d_i}{n_i} - 1\), as required. Since the flag was arbitrary, given any \(V \in \sH_0\) we may find a flag containing \(V\) and, by the same argument, deduce that \(\rk(\Gamma^V \cap \Gamma_i) = k_i - \frac{d_i}{n_i} - 1\) for the value of \(i\) for which \(\rk(V) = R_i\), which is the unique value of \(i\) for which \(V \notin S_i\).
\end{proof}

\begin{remark}
The proof does not just show existence but also provides an explicit construction for the decomposition with constant stabiliser rank property, as well as a description of the \(\Gamma_i\) (up to finite index) in terms of the groups \(\Gamma^{\hat{m}}\) coming from Lemma \ref{lem: consequences of C} for an arbitrary flag. Note that, although a different choice of flag will give different groups \(\Gamma^{\hat{m}}\), with projections in different lines \(Y^{\hat{m}} \leqslant \intl\), the subspaces \(X_i\) that the lines \(Y^{\hat{m}}\) span (taken over all \(m\) for which \(\Gamma^{\hat{m}}\) are equal rank) will not depend on the choice of flag.

It is not hard to show that the decomposition in the above proof is `coarsest', in the sense that for any other decomposition each associated subspace \(X_i\) is contained in one of those from the decomposition defined in the proof. Finer decompositions will, of course, sometimes exist. For example, take the Ammann--Beenker cut and project scheme (see Section \ref{sec: AB}) but replace the window with the square \(W = [0,1]^2 \subset \intl = \R^2\). Then the scheme is already of constant stabiliser rank (each \(\rk(V) = 2\)) and satisfies the above result but can also be decomposed further, where \(X_1\) and \(X_2\) are the two lines along the standard coordinate directions.
\end{remark}

The result below is a corollary of Proposition \ref{prop: factorisation to CSR}, hyperplane spanning from Lemma \ref{lem: consequences of C} and the generalised complexity formula from Theorem \ref{thm: generalised complexity}.

\begin{corollary}
If \(\mathcal{S}\) is indecomposable then it has property {\Cpx} if and only if \(\rk(V) = k - \delta - 1\) for all \(V \in \sH_0\), where \(\delta \coloneqq d/n\). In particular, in this case, the codimension \(n\) divides the dimension \(d\).
\end{corollary}

A decomposable scheme gives a direct sum decomposition of the internal space into complementary subspaces \(X_i\) (with dimensions \(n_i\) summing to \(n\)) which are aligned with the supporting subspaces of the window (each contains all but one \(X_i\)). There is no reason, a priori, for the projected lattice to split compatibly with respect to this decomposition. However, when {\Cpx} holds, \(\Gamma_<\) is well-aligned with the supporting subspaces (the stabiliser ranks are high). In the proof of Proposition \ref{prop: factorisation to CSR}, this can be shown to force the ranks \(k_i = \rk(\Gamma_i)\) to be maximal, in the sense that \(k = k_1 + \cdots + k_\ell\) so, when {\Cpx} holds, the decomposition decomposes the internal space together with the projected lattice (up to finite index). Then a decomposition breaks the scheme into `subsystems', with total dimensions \(k_i\) and codimension \(n_i = \dim(X_i)\) (and thus `physical dimension' \(d_i = k_i - n_i\)).

In fact, in the convex polytopal case, one even has a compatible decomposition of the window \(W = W_1 + \cdots + W_\ell\), where each \(W_i \subset X_i\) is a lower dimensional convex polytope with boundary supported by the intersections of hyperplanes \(H \notin S_i\) with \(X_i\). Of course, we can not expect this for non-convex windows, but this is not important as it is not hard to see that the cut regions split as products, and these may be used to approximate the acceptance domains, given {\Cpx}, by Lemma \ref{lem: refining acceptance domains}. These cut regions will in turn be analysed by their vertices. By considering arbitrary hyperplane intersections, as in the following definition, we obtain the structure of a group for each flag. These will be restricted to the actual vertices in \(W\) needed for \(r\)-cut regions in Definition \ref{def: vertex group subsets}. The second form of \(\Ver(f)\) below follows directly from \(v(f) \in H\) for each \(H \in f\), so that \(H = V(H) + v(f)\).

\begin{definition} \label{def: vertices from flags}
For each flag \(f \in \sF\), we define the vertex intersection points
\[
\Ver(f) \coloneqq \left\{v \in \intl \ \bigg| \ \{v\} = \bigcap_{H \in f} H + (\gamma_H)_<, \text{ for } \gamma_H \in \Gamma\right\} = \left(\bigcap_{V \in V(f)} (V + \Gamma_<)\right) + v(f).
\]
\end{definition}

\begin{proposition} \label{prop: C versus finite index vertices}
For any \(f \in \sF\) we have that \(\Ver(f) - v(f)\) is a free Abelian group of rank \(\sum_{H \in f} (k-\rk(H))\) with \(\Gamma_< \leqslant \Ver(f) - v(f)\). In particular, {\Cpx} holds if and only if \(\Gamma_<\) is finite index in each \(\Ver(f) - v(f)\).
\end{proposition}

\begin{proof}
Since \(\Gamma_<\) and each \(V(H)\) is a group, so is \(\Ver(f) - v(f) = \bigcap_{H \in f} (V(H) + \Gamma_<)\), and the latter also clearly contains \(\Gamma_<\) since \(0 \in V(H)\). It is not hard to show there is an isomorphism
\[
\bigcap_{H \in f} (V(H) + \Gamma_<) \cong \bigoplus_{H \in f} (\Gamma/\Gamma^H),
\]
given by sending the vertex defined by \(\bigcap (V(H) + (\gamma_H)_<)\) to \(([\gamma_H])_{H \in f}\); see the proof of  \cite[Proposition 7.2]{KoiWalII} for further details (where we note that the relevant part of the proof does not require homogeneity). It follows that \(\rk(\Ver(f) - v(f)) = \sum_{H \in f} (k-\rk(H))\), as required. Then \(\Gamma_<\) is finite index in \(\Ver(f) - v(f)\) if and only if the latter has rank \(\sum_{H \in f} (k-\rk(H)) = \sum_{V(H) \in V(f)} (k-\rk(V(H))) = k\) for each flag \(V(f) \in \sF_0\). Equivalently (since each \(f \in \sF_0\) is given by \(V(f')\) for some \(f' \in \sF\)) we have \(nk - \sum_{V \in f} \rk(V) = k\), that is, \(\sum_{V \in f} \rk(V) = k(n-1)\) for each \(f \in \sF_0\), which is equivalent to {\Cpx} by Corollary \ref{cor: sum of ranks for low complexity}.
\end{proof}

\begin{remark} \label{rem: C versus finite index vertices}
The above is relevant to displacements between vertices of cut regions. It follows from Proposition \ref{prop: C versus finite index vertices} that, when {\Cpx} holds, there is some \(N \in \N\) so that
\[
\Ver(f) - \Ver(f') \subseteq \left(\frac{1}{N}\Gamma_< + v(f)\right) - \left(\frac{1}{N}\Gamma_< + v(f')\right) = \frac{1}{N}\Gamma_< + (v(f) - v(f')) \subseteq \frac{1}{N}\Gamma_< + F.
\]
In fact, we can drop the \(F\) term in the weakly homogeneous case:
\end{remark}

\begin{proposition} \label{prop: F for weakly homogeneous}
Given {\Cpx}, weak homogeneity is equivalent to \(F \subset \frac{1}{N}\Gamma_<\) for some \(N \in \N\).
\end{proposition}

\begin{proof}
If the scheme is weakly homogeneous then, for an appropriate translate of \(W\) (one with \(\frac{1}{N} \Gamma_< \cap H \neq \emptyset\) for each \(H \in \sH\)), we have \(V(H) + \frac{1}{N}\Gamma_< = H + \frac{1}{N}\Gamma_<\) for each \(H \in \sH\). Replacing \(\Gamma\) with \(\frac{1}{N}\Gamma\) does not affect {\Cpx}, by Theorem \ref{thm: generalised complexity} (as it does not affect stabiliser ranks or spans). Applying Proposition \ref{prop: C versus finite index vertices} with the finer lattice \(\frac{1}{N}\Gamma_<\), each (finer) \(\Ver(f) - v(f)\) is a finite index subgroup of some \(\frac{1}{M}\Gamma_<\). Moreover, \(\Ver(f) - v(f) = \Ver(f)\), since \(0 \in V(H) \subset H + \frac{1}{N}\Gamma_<\) for each \(H \in f\). In particular, the generalised vertices belong to \(\frac{1}{M}\Gamma_<\), so their difference set \(F\) does too.

Conversely (where we need not assume {\Cpx}), suppose that \(F \subset \frac{1}{N}\Gamma_<\) for some \(N \in \N\), and (by translating \(W\)) assume that some \(v(f) \in \frac{1}{N}\Gamma_<\) for some arbitrary flag \(f \in \sF\). Take any \(H \in \sH\), which is contained in some flag \(f'\), so \(v(f') \in H\). Then \(v(f')-v(f) \in F \subset \frac{1}{N}\Gamma_<\). Thus \(v(f') \in \frac{1}{N}\Gamma_< + v(f) \subset \frac{1}{N}\Gamma_< + \frac{1}{N}\Gamma_< = \frac{1}{N}\Gamma_<\) and thus \(\frac{1}{N}\Gamma_< \cap H \neq \emptyset\). Since \(H \in \sH\) was arbitrary, the scheme is weakly homogeneous.
\end{proof}

We finish this section by defining subgroups of \(\Gamma\) corresponding to the finite index subgroups of Definition \ref{def: vertices from flags}.

\begin{definition} \label{def: lifted vertex groups}
Suppose that {\Cpx} holds, so that for each \(f \in \sF\) we have that \(\Ver(f) - v(f) \leqslant \frac{1}{N}\Gamma_<\) for some \(N \in \N\). Write \(\Gamma[f]\) for the lift of this subgroup to the total space, that is, \(\Ver(f) - v(f) = \Gamma[f]_<\) with \(\Gamma \leqslant \Gamma[f] \leqslant \frac{1}{N}\Gamma\).
\end{definition}

\begin{example}
Consider the Ammann--Beenker cut and project scheme (for more details, see Section \ref{sec: AB}). Two examples of supporting hyperplanes (after translating \(v(f)\) to the origin) are \(H_1 = \langle e_1 \rangle_\R\) and \(H_2 = \langle e_2 \rangle_\R\), for the standard basis vectors \(e_1 = (1,0)\), \(e_2 = (0,1) \in \intl = \R^2\), and \(\Gamma_<\) is given by
\[
\left\langle e_1, e_2, \frac{e_1 + e_2}{\sqrt{2}}, \frac{e_1 - e_2}{\sqrt{2}}  \right\rangle_\Z = \left\{ \left(a_1 + b_1 \frac{\sqrt{2}}{2} , a_2 + b_2 \frac{\sqrt{2}}{2} \right) \middle| a_1, a_2, b_1, b_2 \in \Z, b_1 + b_2 \in 2\Z \right\}.
\]
Taking \(f = \{H_1,H_2\}\) and considering the horizontal and vertical displacements of \(H_1\) and \(H_2\) under translation by \(\Gamma_<\), it is not hard check that they produce intersection points
\[
\Gamma[f]_< =  \left\{ \left(a_1 + b_1\frac{\sqrt{2}}{2} , a_2 + b_2\frac{\sqrt{2}}{2} \right) \middle| a_1, a_2, b_1, b_2 \in \Z\right\}.
\]
Then \(\Gamma\) is index \(2\) in \(\Gamma[f]\) and may be given by replacing the two generators \(a\), \(b \in \Gamma \leqslant \tot \cong \R^4\), mapping to \((e_1\pm e_2)/\sqrt{2} \in \intl\), with the new generators \((a \pm b)/2\) (and leaving the generators mapping to \(e_1\), \(e_2 \in \intl\) the same). A similar analysis can be done for the other supporting hyperplanes.
\end{example}

\section{Diophantine schemes} \label{sec: Diophantine schemes}

In this section we define the notion of a densely embedded lattice being Diophantine. Then, taking any possible decompositions into account, we use this to define the property {\D} of a polytopal cut and project scheme being Diophantine. We finish the section by introducing a new, inhomogeneous version of this that will be relevant for schemes that are not weakly homogeneous.

\subsection{Diophantine numbers, families of vectors and densely embedded lattices}
An irrational number \(\alpha \in \R\) is {\bf badly approximable} if there exists some \(c > 0\) for which
\[
\mathrm{dist}(q\alpha,\Z) \geq \frac{c}{|q|}
\]
for all non-zero \(q \in \Z\), where \(\mathrm{dist}(-,\Z)\) denotes the distance to the nearest integer. We denote the set of badly approximable numbers as \(\Bad\). If \(\alpha \in \Bad\) then \(\mathrm{dist}(q\alpha,\Z) = |q\alpha - n| = |q||\alpha - (n/q)|\), for some \(n \in \Z\), so one may equivalently define \(\alpha \in \Bad\) to mean that
\[
\left| \alpha - \frac{n}{q} \right| \geq \frac{c}{q^2}
\]
for all rational numbers \(n/q\). A classical result is that this is equivalent to \(\alpha\) having bounded entries in its continued fraction expansion (for instance, see \cite[Appendix D]{Bug12}, and \cite{Cass57} for more background on Diophantine Approximation).

More generally, suppose we have a linear map \(L \colon \R^d \to \R^n\), \(L(x) = (x \cdot \alpha_1, x \cdot \alpha_2, \ldots, x \cdot \alpha_n)\); that is, we take a family \(\{\alpha_1,\ldots,\alpha_n\}\) of vectors \( \alpha_i \in \R^d\), generalising from the above case of \(d=n=1\) (where \(L(q) = q\alpha\)). We assume that \(L\) is irrational, in the sense that \(L(n) \notin \Z^n\) for all non-zero \(n \in \Z^d\). Then we call \(L\) {\bf badly approximable} if there exists some \(c > 0\) for which, for all non-zero \(q \in \Z^d\),
\[
\mathrm{dist}(L(q),\Z^n) \geq c |q|^{-\delta},
\]
where \(\delta \coloneqq d/n\) and \(|q| \coloneqq \max |q_i|\) for \(q = (q_1,\ldots,q_d)\). The choice of norm on \(\R^n\) here, for defining the distance \(\mathrm{dist}(-,\Z^n)\) to the integer lattice, is inconsequential (up to changing the constant \(c\)) and similarly for the max norm on \(\R^d\) for \(|q|\). Increasing \(n\) makes it easier to avoid lattice points (more components need to be simultaneously close to integers), whereas increasing \(d\) makes it harder to avoid lattice points (\(L(\Z^d)\) is higher rank), partially explaining the choice of \(\delta \coloneqq d/n\). In fact, \(\delta\) is chosen optimally in the sense that we have the multidimensional Dirichlet's Theorem: for every \(N \in \N\), there exists some non-zero \(q \in \Z^d\) with \(|q| \leq N\) and
\[
\mathrm{dist}(L(q),\Z^n) \leq N^{-\delta}.
\]

A simple extension of the above is to replace \(\Z^n\) with a different lattice, which of course only changes the situation up to linear isomorphism. Instead, to apply the above to the setting of cut and project sets, we wish to remove the reference lattice entirely: we have a \(k = d+n\) densely embedded lattice \(\Gamma_<\) in the internal space but no \emph{natural} splitting of \(\Gamma_< \cong \Z^k\) resembling \(L(\Z^d) + \Z^n\). This motivates the following:

\begin{definition} \label{def: Diophantine}
Let \(G < X\), where \(G \cong \Z^k\) and \(X\) is an \(n\)-dimensional vector space. Choose any isomorphism \(\varphi \colon G \xrightarrow{\cong} \Z^k\) and thus a norm \(\eta(g) = \max_i |g_i|\) on \(G\), where \(\varphi(g) = (g_1,\ldots,g_k) \in \Z^k\). Choose any norm \(\|\cdot\|\) on \(X\). Then we call \(G\) {\bf Diophantine} if there exists some \(c > 0\) so that, for all non-zero \(g \in G\),
\[
\|g\| \geq c\eta(g)^{-\delta}, \text{ where } \delta \coloneqq \frac{d}{n} \text{ for } d \coloneqq k-n.
\]
\end{definition}

\begin{remark}
It is not too hard to see that the Diophantine property only depends on the embedding \(G < X\) up to linear isomorphism of the pair. In particular, it does not depend on the choice of lattice norm \(\eta\) on \(G\), nor on the norm \(\|\cdot\|\) on \(X\), up to changing the constant \(c\). For instance, when applied to cut and project sets, we may take \(\eta(\gamma)\) as the norm of \(\gamma \in \Gamma\) in the total space.
\end{remark}

\begin{remark}
An irrational linear map \(L \colon \R^d \to \R^n\) is badly approximable if and only if \(G = L(\Z^d) + \Z^n < \R^n\) is Diophantine. Conversely, given \(G < X\), we may find \(A\), \(F < G\) with \(G = A+F\), \(A \cong \Z^d\) and \(F \cong \Z^n\) is a lattice in \(X\). Then, up to linear isomorphism, the pair \(F < X\) may be transformed to \(\Z^n < \R^n\). By choosing a basis of \(A \cong \Z^d\), we then obtain a linear map \(L \colon \R^d \to \R^n\) with \(L(\Z^d) = A\) and find that \(L\) is badly approximable if and only if \(G < X\) is Diophantine.

In short, Definition \ref{def: Diophantine} corresponds to the standard definition of a linear map being badly approximable with respect to the integer lattice, but is more natural in our setting, as it is free of a reference lattice: instead of avoiding the integer lattice, points should avoid each other.
\end{remark}

\begin{example} \label{ex: 2-to-1 Diophantine}
By the above, for \(d=n=1\), up to linear isomorphism \(G = \alpha\Z+\Z < \R\) and \(G\) is Diophantine if and only if \(\alpha \in \R \setminus \Q\) is badly approximable. In more detail, for non-zero \(g = n_1 \alpha + n_2 \in G\), we consider the size of \(|n_1 \alpha + n_2|\) relative to, say, \(\eta(g) = \max\{|n_1|,|n_2|\}\). We may as well assume that \(n_1 \neq 0\), or else \(|g| = |n_2| \geq 1\). Otherwise, without loss of generality, we may take \(-n_2\) to be the nearest integer to \(n_1 \alpha\), to minimise \(|n_1 \alpha + n_2| = \mathrm{dist}(n_1 \alpha,\Z)\), so \(\eta(g) \asymp n_1\). Thus, the Diophantine condition equivalently requires \(\mathrm{dist}(n_1 \alpha,\Z) \geq c|n_1|^{-1}\) for all non-zero \(n_1 \in \Z\), which is what it means to be badly approximable.
\end{example}

Diophantine lattices have a useful dual property, that they evenly distribute points in bounded regions without leaving large gaps:

\begin{definition}[{\cite[Definition 5.6]{KoiWalII}}]
We call \(G < X\) {\bf densely distributed} if there is some \(c > 0\) so that, for all \(r > 0\), every point in the unit ball of \(X\) is within distance \(cr^{-\delta}\) of some \(g \in G\) with \(\eta(g) \leq r\).
\end{definition}

Note that the unit ball is arbitrary in the above, it can be replaced by any bounded subset with non-empty interior (see \cite[Lemma 5.7]{KoiWalII}). The following (proved in \cite{KoiWalII} for completeness) is a restatement in our setting of a result usually given in terms of badly approximable numbers or families of vectors, which follows from more general transference results between homogeneous and inhomogeneous problems, see \cite[Chapter V.4]{Cass57}:

\begin{theorem}[{\cite[Theorem 5.8]{KoiWalII}}] \label{thm: Diophantine => densely distributed}
If \(G < X\) is Diophantine then it is densely distributed.
\end{theorem}

Passing between finite index sub- and super-groups does not affect the Diophantine condition; it is not difficult to prove the following:

\begin{lemma}[{\cite[Lemma 5.3]{KoiWalII}}] \label{lem: Diophantine under finite index}
Let \(G \leqslant K\), where \(G\) and \(K\) are free Abelian groups of equal rank \(k\), densely embedded in a vector space \(X\). Then \(G\) is Diophantine if and only if \(K\) is.
\end{lemma}

\subsection{Inhomogeneous approximation}
Next we introduce an inhomogeneous analogue of the Diophantine condition. We recall (see, for instance, \cite{Kle99}) that, given an irrational \(\alpha \in \R\) and \(y \in \R\), we say that \(\alpha\) is {\bf inhomogeneously badly approximable with respect to} \(y\), and write \(\alpha \in \Bad(y)\), if there exists some \(c > 0\) so that
\[
\mathrm{dist}(q\alpha - y,\Z) \geq \frac{c}{|q|}
\]
for all non-zero \(q \in \Z\) with \(q\alpha - y \notin \Z\); note that, since \(\alpha\) is irrational, there can be at most one value of \(q\) with \(q \alpha - y \in \Z\), in which case \(y \in \alpha \Z + \Z\) and the situation reduces to whether or not \(\alpha \in \Bad\). Again, we may introduce higher dimensional versions of this but, for our purposes, we want a form free of any reference lattice:

\begin{definition} \label{def: inhomogeneous Diophantine}
Let \(G\) and \(X\) be as in Definition \ref{def: Diophantine}, and let \(F \subset X\) be non-empty and finite. We say that \(G\) is {\bf (inhomogeneously) Diophantine} with respect to \(F\) if there exists some \(c > 0\) so that, for all \(f \in F\) and \(g \in G\) with \(g \neq f\),
\[
\|g - f\| \geq c\eta(g)^{-\delta}, \text{ where } \delta \coloneqq \frac{d}{n} \text{ for } d \coloneqq k-n.
\]
\end{definition}

In other words, instead of needing to avoid \(0\), we now need \(G\) to avoid elements of \(F\), in the sense that they are distant relative to their lattice norm (except when they exactly coincide with elements from \(F\), which will be a convenient convention here).

\begin{remark} \label{rem: hom versus inhom Diophantine}
Generally, if \(F \subset \Q G\) and \(G\) is (homogeneously) Diophantine then it is Diophantine with respect to \(F\). This is not hard to prove; briefly, if \(g\) is close to some \(f \in F\) and \(Nf \in G\), then \(Ng-Nf\) is close to \(0\), which forces \(\eta(g) \asymp \eta(Ng-Nf)\) to be reasonably large if \(G\) is Diophantine.

The converse is not true: one may have \(F \subset \Q G\) and \(G\) Diophantine with respect to \(F\) but \(G\) not Diophantine, see Example \ref{ex: inhomogeneous bad not integer scaling invariant}. However, we note the following:
\end{remark}

\begin{remark} \label{rem: DF => D}
When applied to cut and project sets, with \(F\) defined as in Definition \ref{def: vertices}, we always have that \(0 \in F\), so that \(G\) being inhomogeneously Diophantine with respect to \(F\) implies that \(G\) is also homogeneously Diophantine.
\end{remark}

It is useful to be able to swap a given \(G\) for a finite index sub or super-group. Clearly replacing \(G\) for a sparser finite index subgroup preserves both homogeneous and inhomogeneous Diophantine properties. The homogeneous Diophantine condition is not affected by passing to finite index super-groups (Lemma \ref{lem: Diophantine under finite index}). However, the question of whether this must hold in the inhomogeneous case appears to the author to be a difficult problem in Diophantine approximation:

\begin{question} \label{q: general question in DA}
Suppose that \(G < X\) is inhomogeneously Diophantine with respect to \(F\), and that \(0 \in F\), i.e., \(G\) is in particular homogeneously Diophantine. Must it hold that \(G' < X\) is inhomogeneously Diophantine with respect to \(F\) whenever \(G\) is a finite index subgroup of \(G'\)?
\end{question}

If the above can be answered in the positive then Theorem \ref{thm: main2} would provide a single necessary and sufficient condition for \LR, that both {\Cpx} and {\DF} hold, valid for arbitrary polytopal cut and project schemes. Otherwise, instead of {\DF}, one may need to check several inhomogeneous Diophantine conditions over different finite index supergroups of \(\Gamma_<\), as we will see in Theorem \ref{thm: iff for LR}. For concreteness, the above question for \(d=n=1\) may be phrased as follows:

\begin{question} \label{q: question in DA}
Suppose that \(\alpha \in \Bad \cap \Bad(y)\). Do we have that \(\alpha/n \in \Bad(y)\) for any \(n \in \N\)?
\end{question}

We assume that \(\alpha \in \Bad\) in the above since otherwise there are counterexamples\footnote{The author thanks Alan Haynes for pointing out this example.}:

\begin{example} \label{ex: inhomogeneous bad not integer scaling invariant}
Let \(d=n=1\), \(F = \{\frac{1}{2}\}\) and \(G = \alpha \Z + \Z < \R\) for \(\alpha > 0\) irrational. Translating into more standard language, we consider the properties of irrationals being in \(\Bad\) and \(\Bad(\frac{1}{2})\). We have \(q | q\alpha - \frac{1}{2} - p | < C\) for \(p \in \Z\), \(q \in \N\) if and only if \(q' |q'\alpha - p'| < 4C\) for \(q' = 2q\) and \(p' = 2p+1\). By the theory of continued fractions, for \(C < \frac{1}{8}\), we have that \(q' |q'\alpha - p'| < 4C\) if and only if \(p'\) and \(q'\) are common integer multiples of, respectively, the numerator and denominator of a principle convergent \(\frac{p_i}{q_i}\) to \(\alpha\). For \(q'\) even, \(q = \frac{q'}{2}\) is even if and only if \(q' \equiv 0 \mod 4\). So consider \(\alpha \notin \Bad\) with \(q_i \not\equiv 0 \mod 4\) for all but finitely many principal convergents \(\frac{p_i}{q_i}\). One example is given by
\[
\alpha = [a_0;a_1,a_2,a_3,\ldots] = [0;2,1,4,8,12,16,20,\ldots] = \frac{1}{2+\frac{1}{1+\frac{1}{4 + \frac{1}{8+\cdots}}}}.
\]
We have the recursive formulae \(p_i = a_np_{i-1} + p_{i-2}\), \(q_i = a_nq_{i-1} + q_{i-2}\), with \(p_0 = 0\), \(q_0 = 1\), \(p_1 = 1\) and \(q_1 = 2\). Since \(a_i \equiv 0 \mod 4\) for \(i \geq 3\), by induction \(q_i \equiv a_iq_{i-1} + q_{i-2} \equiv q_{i-2} \equiv 2 \mod 4\) for all \(i\) odd. Similarly, \(q_i\) is odd for \(i\) even. Finally, \(p_i\) is odd for all \(i \geq 1\).

It follows from the above that \(\alpha \notin \Bad(\frac{1}{2})\) but that \(2\alpha \in \Bad(\frac{1}{2})\). Indeed, since the \(a_i \to \infty\), the principal convergents \(\frac{p_i}{q_i}\) with \(q_i\) even give close inhomogeneous approximations of \(\alpha\) to \(\frac{1}{2}\) (taking \(p = \frac{p_i-1}{2}\) and \(q = \frac{q_i}{2}\)), so \(\alpha \in \Bad(\frac{1}{2})\). However,
\[
\inf q \left| q(2\alpha) - \frac{1}{2} - p \right| = \frac{1}{2} \inf (2q) \left| (2q)\alpha - \frac{1}{2} - p \right| > 0,
\]
since otherwise (using that \(2q\) is even), there would be infinitely many principal convergents of \(\alpha\) with odd numerator and denominator \(\equiv 0 \mod 4\), which is not the case for our \(\alpha\). We see that \(x \in \Bad(d)\) does not in general imply that \(x/n \in \Bad(d)\) for \(n \in \N\).

This example also shows that, for \(F \subset \Q G\), it is possible for \(G\) to be Diophantine with respect to \(F\) but not homogeneously Diophantine. Indeed, consider \(G = 2\alpha \Z + \Z\) with \(\alpha\) as above, and \(F = \{\frac{1}{2}\}\). We have that \(2\alpha \in \Bad(\frac{1}{2})\), equivalently, \(G\) is Diophantine with respect to \(F \subset \frac{1}{2}\Z \subset \Q G\) but \(2\alpha \notin \Bad\), since otherwise \(\alpha \in \Bad\), which we have seen is false.
\end{example}

\subsection{Diophantine schemes}
The Diophantine condition for a densely embedded lattice may be directly applied to cut and project schemes, by asking that \(\Gamma_<\) is Diophantine in \(\intl\). Unfortunately, this definition is not quite the correct one for our main theorems when stabilisers are not of constant rank. This can occur for decomposable schemes, which can satisfy {\Cpx} whilst having factors with varying \(\delta_i = d_i/n_i\). To account for this, we apply the Diophantine condition to factors of constant stabiliser rank, in the sense given in Proposition \ref{prop: factorisation to CSR}.

\begin{definition} \label{def: Diophantine scheme}
Let \(\mathcal{S}\) be a polytopal cut and project scheme with constant stabiliser rank. We call \(\mathcal{S}\) {\bf Diophantine}, and say it has property {\D}, if \(\Gamma_<\) is Diophantine in \(\intl\).

More generally, suppose that \(\mathcal{S}\) is given a decomposition satisfying the conclusions of Proposition \ref{prop: factorisation to CSR}. Then \(\mathcal{S}\) has property {\D} if each \((\Gamma_i)_<\) is Diophantine, considered as a free Abelian subgroup of rank \(k_i\) in the \(n_i\)-dimensional vector space \(X_i\).
\end{definition}

\begin{remark}
In the constant stabiliser rank case, we may take \(\eta(g_<) = \|g\|\) for the lattice norm in Definition \ref{def: Diophantine}, where \(\|\cdot\|\) is any norm on \(\tot\). Then the Diophantine property makes precise the notion that lattice points do not project close to the origin in \(\intl\), relative to their distance from the origin in \(\tot\). Equivalently, lattice points remain distant from the physical space \(\phy\), relative to their norm. This still holds in the decomposable case, but exactly how `close' is defined may depend on a splitting of the internal space.
\end{remark}

\begin{remark} \label{rem: independence of C and D}
If the window is indecomposable then, if it also satisfies {\Cpx}, then it is constant stabiliser rank by Corollary \ref{cor: constant stabiliser rank} and thus we may apply the simpler definition of \D\ that \(\Gamma_<\) is Diophantine (and, indeed, this definition makes sense independently of \Cpx). If \Cpx\ is not satisfied then \LR\ fails regardless of \D.

In this sense \Cpx\ and \D\ may, up to a minor technicality due to the decomposable case, be considered practically independent in nature. However, in the decomposable case the scheme may fail to be constant stabiliser rank and \(\Gamma_<\) being Diophantine is then not the appropriate condition: we must use the slightly more complicated version in this case, which is defined relative to a splitting that can only be guaranteed by \Cpx. However, even in this case, we note that \D\ is equivalent to each \(\Gamma^{\hat{m}}\) (for \(m = 1,\ldots,n\)) being Diophantine as a subgroup of a line (so with \(n=1\) in Definition \ref{def: Diophantine}), where each \(\Gamma^{\hat{m}}\) is defined as in Lemma \ref{lem: consequences of C} with respect to an arbitrary flag \(f \in \sF_0\). This follows from the fact that each \(\Gamma_i\) in a decomposition satisfying Proposition \ref{prop: factorisation to CSR} is given (up to finite index) by the sum of \(\Gamma^{\hat{m}}\) contained in \(X_i\), and all such \(\Gamma^{\hat{m}} < X_i\) are of constant rank. This is not hard to deduce in a manner similar to the details of the proof of Proposition \ref{prop: factorisation to CSR}: the restricted stabiliser ranks of Proposition \ref{prop: factorisation to CSR} will be \(\rk(\Gamma^V \cap \Gamma_i) = k_i - \rk(\Gamma^{\hat{m}})\), so \(\rk(\Gamma^{\hat{m}})\) is independent on the value of \(i\) for which \(\Gamma^{\hat{m}} < X_i\), and it is not hard show that each \(\Gamma^{\hat{m}}\) (of equal rank) is Diophantine if and only if their sum is (see \cite[Lemma 5.5]{KoiWalII}). So property \D\ can always be reduced to the case of checking systems of numbers being simultaneously Diophantine in \(\R\) in the classical sense, where this needs to be checked over \(n\) different such systems, and this characterisation of \D\ does not require reference to property \Cpx\ (although if \Cpx\ fails then the sum of groups \(\Gamma^{\hat{m}}\) will not be finite index in \(\Gamma\) and knowing these ranks is obviously a basic necessity in checking the Diophantine condition). This also shows that \D\ does not depend on the decomposition satisfying the conclusions of \ref{prop: factorisation to CSR} that we choose, which gives a potentially useful flexibility in applying our results.
\end{remark}

The above definition is appropriate in the case that the scheme is weakly homogeneous. Otherwise, we will need an inhomogeneous condition:

\begin{definition} \label{def: DF}
Let \(\mathcal{S}\) be a polytopal cut and project scheme with constant stabiliser rank and let \(F \subset \intl\) be as in Definition \ref{def: vertices}. We call \(\mathcal{S}\) {\bf Diophantine with respect to \(F\)}, and say it has property {\DF}, if \(\Gamma_<\) is Diophantine in \(\intl\) with respect to \(F\).

More generally, suppose that \(\mathcal{S}\) is given a decomposition satisfying the conclusions of Proposition \ref{prop: factorisation to CSR} and write \(F = F_1 + \cdots + F_\ell\) for subsets \(F_i \subset X_i\). Then \(\mathcal{S}\) (with respect to this splitting of \(\Gamma\)) has property {\DF} if each \((\Gamma_i)_<\) is Diophantine with respect to \(F_i\), considered as a free Abelian subgroup of rank \(k_i\) in the \(n_i\)-dimensional vector space \(X_i\).
\end{definition}

We recall (Definition \ref{def: vertices}) that
\[
F \coloneqq \{v(f) - v(f') \mid f, f' \in \sF\},
\]
i.e., the displacements between generalised vertices of \(W\), where \(v(f)\) is the intersection point of a flag \(f \in \sF\). Taking \(f = f'\), it follows that \(0\) is always an element of \(F\) (Remark \ref{rem: DF => D}). Moreover, any decomposition \(Q\) (which induces a splitting \(\intl = X_1 + \cdots + X_\ell\)), induces a corresponding splitting of \(F = F_1 + \cdots + F_\ell\), for finite sets \(F_i \subset X_i\); see Figure \ref{fig: decomposable} for a codimension \(n=2\) example. Indeed, for each \(i\) we may consider the hyperplanes in \(X_i\)
\[
\sH_i = \{H \cap X_i \mid H \in \sH \setminus S_i\}.
\]
By the definition of a decomposition, choosing a flag in \(\sH\) is equivalent to choosing a flag in each \(\sH_i\). Indeed, for any flag \(f = \{H_1,\ldots,H_n\} \in \sF_0\), the dimension of \(H_1 \cap \cdots \cap H_m\) must be \(n-m\) (i.e., decreasing by \(1\) each time we intersect a new subspace), but intersecting \(m > n_i = \dim(X_i)\) hyperplanes from any \(S_i^\mathrm{c}\) will have dimension strictly larger than \(n-m\), since all such subspaces are parallel to the subspaces \(X_j\), \(j \neq i\), whose sum has dimension \(n-n_i\). We have the generalised vertices
\[
v(\sH_i) \coloneqq \{v(f) \mid f \in \sF_i\},
\]
considered as a finite subset of \(X_i\), where \(\sF_i\) is the set of flags (of \(X_i\)) in \(\sH_i\). Displacements between these define sets
\[
F_i = \{v(f) - v(f') \mid f, f' \in \sF_i\}.
\]
Then \(v(\sH) = v(\sH_1) + \cdots + v(\sH_\ell)\) since, by the above, choosing a flag in \(\sH\) is the same as choosing one flag for each \(\sH_i\), each of which determines the \(X_i\) component (by the definition of a decomposition) and thus \(F = F_1 + \cdots + F_\ell\).

\section{Necessary and sufficient conditions for {\LR}} \label{sec: main proof}

We begin this section with some technical notation and results linking vertices of cut regions to the projected lattice. We then quantify how vertices remaining distant corresponds to acceptance domains containing large balls (or, otherwise, having small volumes). These ideas are combined in Theorem \ref{thm: iff for LR} to give a single necessary and sufficient condition for {\LR}, that both {\Cpx} and a certain Diophantine condition holds, in terms of the flag groups. From these, the more easily stated Theorems \ref{thm: main1} and \ref{thm: main2} quickly follow.

\subsection{Cut region vertices and the projected lattice}

Since we will be interested in vertices of \(r\)-cut regions in the window, we introduce the following notation (see also Definition \ref{def: vertices from flags}):

\begin{definition} \label{def: vertex group subsets}
Given \(f \in \sF\) and \(r \geq 0\) we write
\[
\Ver(f,r) \coloneqq \{v \in \intl \mid \{v\} = \bigcap_{H \in f} (H + (\gamma_H)_<), \text{ for } \gamma_H \in \Gamma \cap B_r\}.
\]
We use the following notation for vertices of \(r\)-cut regions in the window:
\[
\Ver_W(f,r) \coloneqq \Ver(f,r) \cap W, \ \ \Ver_W(r) = \bigcup_{f \in \sF} \Ver_W(f,r).
\]
\end{definition}

The \(\gamma_H\) each range over all of \(\Gamma \cap B_r\) in the definition of \(\mathrm{Ver}(f,r)\) (and we have such a choice for each \(H \in f\)). Recall (Definition \ref{def: lifted vertex groups}) the groups \(\Gamma[f] \leqslant \Gamma\) of lifts of flag groups, defined when {\Cpx} holds:
\[
\Ver(f) - v(f) = \{v \in \intl \mid \{v\} = \bigcap_{H \in f} (V(H) - (\gamma_H)_<), \text{ for } \gamma_H \in \Gamma(r)\} = \Gamma[f]_<.
\]
Since these are subgroups of the total space (contained in \(\frac{1}{N}\Gamma\) for some \(N \in \N\)), they naturally inherit its norm, so we may define
\[
\Gamma[f;r] \coloneqq  \Gamma[f] \cap B_r.
\]
We obtain the following close correspondence between projections of points in \(\Gamma[f]\) and the `vertices' (extending out of \(W\)) defined by \(f\):

\begin{lemma} \label{lem: vertices and projected groups}
Given {\Cpx} there exists \(c > 0\) satisfying the following, for \(r \geq c\)
\[
\Gamma[f;r - c]_< \subseteq \Ver(f,r) - v(f) \subseteq \Gamma[f;r+c]_<.
\]
\end{lemma}

\begin{proof}
Let \(\gamma^1\), \ldots, \(\gamma^m \in \Gamma[f] \leqslant \tot\) denote representatives of each coset of \(\Gamma\) in \(\Gamma[f]\), so we may write \(\{(\gamma^i)_<\} = \bigcap_{H \in f} (H + (\gamma^i_H)_<)-v(f) = \bigcap_{H \in f} (V(H) + (\gamma^i_H)_<)\), where each \(\gamma^i_H \in \Gamma\). Since there are only finitely many such \(\gamma^i\) and \(\gamma^i_H\), we may assume that each belongs to \(\Gamma(c_1) = \Gamma \cap B_\alpha\) for some \(c_1 \geq 0\).

Take any \(\gamma \in \Gamma[f;r]\), so that \(\gamma + \gamma^i \in \Gamma(r+c_1)\) for some \(i\). Then
\[
\{(\gamma)_<\} = (\gamma + \gamma^i)_< - \{\gamma^i_<\} = \bigcap_{H \in f} (V(H)+(\gamma + \gamma^i - \gamma^i_H)_<).
\]
Since each \(\gamma + \gamma^i - \gamma^i_H \in \Gamma(r+2c_1)\), it follows from Definition \ref{def: vertex group subsets} that \(\gamma_< \in \Ver(f,r+2c_1) - v(f)\), so the first required inclusion holds for all \(r \geq c = 2c_1\).

For the second inclusion, let \(v \in \Ver(f,r) - v(f)\), i.e., \(\{v\} = \bigcap_{H \in f} (V(H)+(\gamma_H)_<) \) where \(\gamma_H \in \Gamma(r)\). We denote \(v = g_<\) for \(g \in \Gamma[f]\), which we wish to show has norm uniformly bounded away from \(r\). We have \(g + \gamma^i = \gamma\) for some \(i\) and \(\gamma \in \Gamma\), so
\[
\bigcap_{H \in f} (V(H) + \gamma_<) = \{\gamma_<\} = \{(g + \gamma^i)_<\} = \{v + \gamma^i_<\} = \bigcap_{H \in f} (V(H) + (\gamma_H + \gamma^i)_<).
\]
By definition, \(V(H) + \gamma_<^i = V(H) + (\gamma_H^i)_<\) (since \(\gamma_<^i \in V(H) + (\gamma_H^i)_<\)), so by the above \(V(H) + \gamma_< = V(H) + (\gamma_H + \gamma^i)_< = V(H) + (\gamma_H + \gamma_H^i)_< \) for each \(H \in f\). Thus \(\gamma - (\gamma_H + \gamma_H^i) \in \Gamma^H\) for each \(H \in f\). In particular, for any \(H \in f\) we have
\[
\gamma - \gamma_H \in \bigcap_{H \in f} (\Gamma^H - \gamma_H^i).
\]
The above intersection must be a singleton, since each \((\Gamma^H - \gamma_H^i)_< \subset V(H) - (\gamma_H^i)_<\), which intersect to a point as \(H\) runs over a flag (and the projection \(\pi_<\) is injective on \(\Gamma\)). Thus, there is a global bound \(\|\gamma - \gamma_H\| \leq c_2\), depending only on the \(\Gamma^H\) and \(\gamma_H^i\). Since any of the \(\|\gamma_H\| \leq r\), it follows that \(\|\gamma\| \leq r+c_2\) and thus \(\|g\| = \|\gamma - \gamma^i\| \leq r + c_1 + c_2\) and the second inclusion follows by taking any \(c \geq c_1 + c_2\).
\end{proof}

The following result relates the elements of \(\Ver(f,r)\) and the actual vertices \(\Ver_W(f,r)\) of the \(r\)-acceptance domains:

\begin{lemma} \label{lem: general to cut region vertices}
There exist \(\epsilon\), \(C > 0\) so that, for sufficiently large \(r \geq 0\) and all \(f \in \sF\),
\[
(\Ver(f,r) - \Ver(f',r)) \cap B_\epsilon \subseteq \Ver_W(f,Cr) - \Ver_W(f',Cr) \subseteq \Ver(f,Cr) - \Ver(f',Cr)
\]
\end{lemma}

\begin{proof}
The second inclusion is trivial. For the first inclusion, choose any \(2\epsilon\)-ball \(B \subset W\) and a basis \(b\) of \(\intl\) in \(\Gamma_<\) making \(\langle b \rangle_\Z\) sufficiently dense to hit every \(\epsilon\)-ball. Given \(v \in \Ver(f,r)\) and \(v' \in \Ver(f',r)\) with \(\|v-v'\| \leq \epsilon\), take any \(\gamma_< \in \langle b \rangle_\Z\) so that \(\|v-\gamma_<\|\) is within \(\epsilon\) of the centre of \(B\). Then \(v' - \gamma_< \in B\) too, so both \(v - \gamma_<\) and \(v' - \gamma_< \in W\). Since we picked \(\gamma_<\) from a basis, it is not hard to see that \(\|\gamma\| \asymp \|\gamma_<\| \asymp \|v\| \ll r\). Thus, both \(v - \gamma_< \in \Ver_W(f,Cr)\) and \(v' - \gamma_< \in \Ver_W(f',Cr)\) for an appropriate \(C > 0\) and sufficiently large \(r\).
\end{proof}

\begin{notation} \label{not: double flag groups}
For flags \(f\), \(f' \in \sF\), we write \(\Gamma[f,f'] \coloneqq \Gamma[f] - \Gamma[f']\). This is a free Abelian group with \(\Gamma \leqslant \Gamma[f,f'] \leqslant \frac{1}{N}\Gamma\) for some \(N \in \N\) (since the same is true of \(\Gamma[f]\) and \(\Gamma[f']\)) and we may equip \(\Gamma[f,f']\) with the norm from \(\tot\). We denote \(\Gamma[f,f';r] \coloneqq \Gamma[f,f'] \cap B_r\).
\end{notation}

The following result directly relates displacements between vertices with elements of the above (shifted) flag groups:

\begin{proposition} \label{prop: vertices versus flag groups}
Given {\Cpx} there exists some \(\epsilon > 0\) and \(C > 0\) so that, for all \(f\), \(f' \in \sF\) and sufficiently large \(r\),
\[
(\Gamma[f,f';r]_< + (v(f) - v(f'))) \cap B_\epsilon \subseteq \Ver_W(f,Cr) - \Ver_W(f',Cr)
\]
and
\[
\Ver_W(f,r) - \Ver_W(f',r) \subseteq \Gamma[f,f';Cr]_< + (v(f) - v(f'))
\]
\end{proposition}

\begin{proof}
Combining Lemma \ref{lem: vertices and projected groups} and Lemma \ref{lem: general to cut region vertices}, there exist \(c\), \(C\) and \(\epsilon > 0\) with
\begin{equation} \label{eq: ver1}
\left(\Gamma[f;r - c]_< - \Gamma[f;r - c]_< + (v(f)-v(f'))\right) \cap B_\epsilon \subseteq \Ver_W(f,Cr) - \Ver_W(f',Cr).
\end{equation}
It is easy to see that
\begin{equation} \label{eq: ver2}
\Gamma[f;r]_< - \Gamma[f';r]_< = (\Gamma[f;r] - \Gamma[f';r])_< \supseteq ((\Gamma[f] - \Gamma[f']) \cap B_{r-c})_< = \Gamma[f,f';r-c]_<
\end{equation}
for sufficiently large \(r\). Indeed, the first equality is trivial (\(\pi_<\) is a homomorphism) and the last one is by definition. For the inclusion, suppose that \(a-b \in (\Gamma[f] - \Gamma[f']) \cap B_r\), with \(a \in \Gamma[f]\) and \(b \in \Gamma[f']\). Since \(\Gamma\) is finite index in both \(\Gamma[f]\) and \(\Gamma[f']\), there is some \(\gamma \in \Gamma \subseteq \Gamma[f] \cap \Gamma[f']\) with \(\|a-\gamma\| \leq c\) (increasing \(c\), if required, and only depending on the density of \(\Gamma \leqslant \Gamma[f]\)). Then \(a-b = (a-\gamma) - (b-\gamma)\), \(\|a-\gamma\| \leq c\) and \(\|b-\gamma\| = \|(a-\gamma) - (a-b)\| \leq \|c\| + \|a-b\| \leq r + c\), as required (substituting \(r\) with \(r-c\)).

Combining Equations \eqref{eq: ver1} and \eqref{eq: ver2} (after shifting \(r\) by \(-c\)), the first inclusion of Proposition \ref{prop: vertices versus flag groups} follows for sufficiently large \(r\) (and perhaps increasing \(C > 0\)). By the same lemmas as above,
\[
\Ver_W(f,r) - \Ver_W(f',r) \subseteq \Ver(f,r) - \Ver(f',r) \subseteq \Gamma[f;r+c]_< - \Gamma[f';r+c]_< + (v(f)-v(f')).
\]
Clearly \(\Gamma[f;r+c]_< - \Gamma[f';r+c]_< \subseteq (\Gamma[f] - \Gamma[f']) \cap B_{2(r+c)} = \Gamma[f,f';2r+2c]\), so the result follows.
\end{proof}

\subsection{Small acceptance domains versus close vertices}
In this section we relate existence (or not) of close vertices to acceptance domains containing large boxes. First, we note that the flag groups split exactly (not just up to finite index) relative to the decomposition:

\begin{definition}
Given {\Cpx} and a decomposition \(Q = \{S_1, \ldots, S_\ell\}\), we write
\[
\Gamma_i[f] \coloneqq \{\gamma \in \Gamma[f] \mid \gamma_< \in X_i\},
\]
i.e., the unique subgroup of \(\Gamma[f]\) defined by \((\Gamma_i[f])_< = \Gamma[f]_< \cap X_i\). Similarly, we write
\[
\Gamma_i[f,f'] \coloneqq \{\gamma \in \Gamma[f,f'] \mid \gamma_< \in X_i\}.
\]
We write \(\Gamma_i[f;r]\) and \(\Gamma_i[f,f';r]\) for the restrictions of the above to \(B_r\).
\end{definition}

\begin{lemma} \label{lem: exact splitting of decomposition}
Given {\Cpx}, for each \(f \in \sF\) we have
\[
\Gamma[f] = \Gamma_1[f] + \Gamma_2[f] + \cdots + \Gamma_\ell[f].
\]
Thus, for two flags \(f\), \(f' \in \sF\),
\[
\Gamma[f,f'] = \Gamma_1[f,f'] + \Gamma_2[f,f'] + \cdots + \Gamma_\ell[f,f'].
\]
\end{lemma}

\begin{proof}
Let \(\gamma \in \Gamma[f]\) and \(v = \gamma_<\), so that \(\{v\} = \bigcap_{H \in f} (V(H) - (\gamma_H)_<)\) for \(\gamma_H \in \Gamma\). Let
\[
\{v_i\} \coloneqq \bigcap_{H \in f} V(H) - (\gamma_H')_<,
\]
where \(\gamma_H' \coloneqq \gamma_H\) for \(H \in S_i^\mathrm{c}\), and \(\gamma_H' = 0\) otherwise. Then \(v_i \in \Gamma[f]_<\), so we may write \(v_i = (\gamma_i)_<\) for some \(\gamma_i \in \Gamma[f]\). In fact, \(\gamma_i \in \Gamma_i[f]\), since \(v_i \in \bigcap_{H \in f \cap S_i} V(H) = X_i\). On the other hand, we have
\[
v, v_i \in \bigcap_{H \in f \cap S_i^\mathrm{c}} (V(H) + (\gamma_H')_<) = \bigcap_{H \in f \cap S_i^\mathrm{c}} (V(H) + (\gamma_H)_<),
\]
where the above is a common translate of the direct sum of the \(X_j\), \(j \neq i\). It follows that \(v\) projects to \(v_i\), relative to the splitting \(\intl = X_1 + \cdots + X_\ell\), so that \(v = v_1 + \cdots + v_\ell\) and \(\gamma = \gamma_1 + \cdots + \gamma_\ell\), as required. The second claim follows trivially from the first, by taking the difference of elements of two flag groups then splitting each, as above.
\end{proof}

We now wish to show that the existence of small acceptance domains corresponds exactly to elements of the flag groups \(\Gamma_i[f,f';r]\) staying far from \(v(f')_i-v(f)_i\), relative to \(r\). To do this, we need to show how a cut region of small volume can be created, given close vertices of two such flags. To this end, we have the following lemma. Intuitively, it says that a hyperplane from one of the flags always chops off a small corner of a region bounded by the other, or at least chops off a prism of small diameter across one dimension, which thus has small volume when bounded in the remaining dimensions.

\begin{lemma} \label{lem: cutting corners}
Consider two flags \(f\) and \(f'\) in \(X \cong \R^n\) and an arbitrary \(\nu > 0\). Then there is some \(C>0\) satisfying the following. For each \(H \in f\), choose any translate \(H + b_H \neq H\), with the two hyperplanes separated by distance at most distance \(\nu\). If \(0 < \|v(f)-v(f')\| \leq \epsilon\), then there is some \(H' \in f'\) and \(A \subset X\) with \(\vol(A) \leq C\epsilon\) which is a connected component of
\begin{equation} \label{eq: euclidean minus flags}
X \setminus \left(\left(\bigcup_{H \in f} H + \{0,b_H, -b_H\} \right) \cup  H'\right).
\end{equation}
\end{lemma}

\begin{figure}
	\def\svgwidth{.5\textwidth}
	\centering
\begingroup%
  \makeatletter%
  \providecommand\color[2][]{%
    \errmessage{(Inkscape) Color is used for the text in Inkscape, but the package 'color.sty' is not loaded}%
    \renewcommand\color[2][]{}%
  }%
  \providecommand\transparent[1]{%
    \errmessage{(Inkscape) Transparency is used (non-zero) for the text in Inkscape, but the package 'transparent.sty' is not loaded}%
    \renewcommand\transparent[1]{}%
  }%
  \providecommand\rotatebox[2]{#2}%
  \newcommand*\fsize{\dimexpr\f@size pt\relax}%
  \newcommand*\lineheight[1]{\fontsize{\fsize}{#1\fsize}\selectfont}%
  \ifx\svgwidth\undefined%
    \setlength{\unitlength}{340.15748031bp}%
    \ifx\svgscale\undefined%
      \relax%
    \else%
      \setlength{\unitlength}{\unitlength * \real{\svgscale}}%
    \fi%
  \else%
    \setlength{\unitlength}{\svgwidth}%
  \fi%
  \global\let\svgwidth\undefined%
  \global\let\svgscale\undefined%
  \makeatother%
  \begin{picture}(1,1)%
    \lineheight{1}%
    \setlength\tabcolsep{0pt}%
    \put(0,0){\includegraphics[width=\unitlength,page=1]{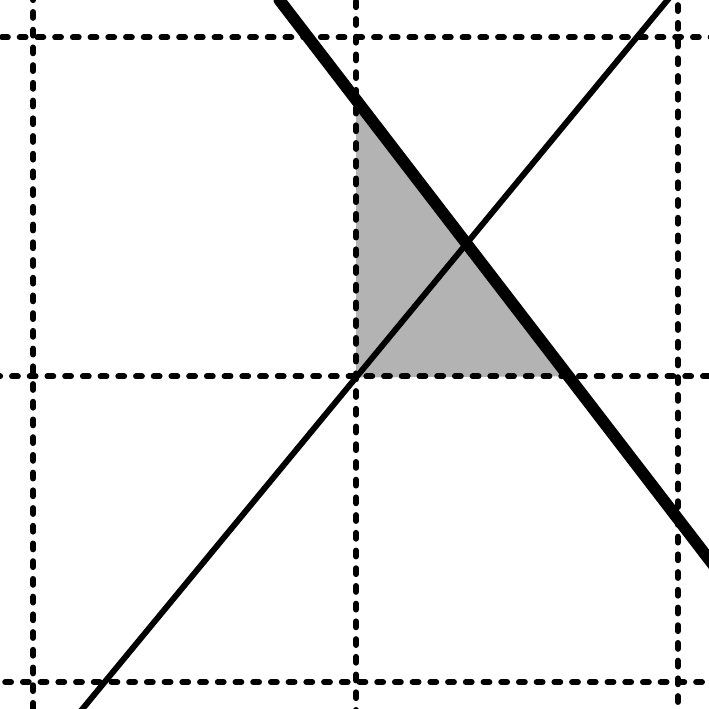}}%
    \put(0.3988437,0.5192186){\color[rgb]{0,0,0}\makebox(0,0)[lt]{\lineheight{1.25}\smash{\begin{tabular}[t]{l}$v(f)$\end{tabular}}}}%
    \put(0.70271066,0.62948701){\color[rgb]{0,0,0}\makebox(0,0)[lt]{\lineheight{1.25}\smash{\begin{tabular}[t]{l}$v(f')$\end{tabular}}}}%
    \put(0,0){\includegraphics[width=\unitlength,page=2]{cutting_lemma.pdf}}%
    \put(0.87369906,0.84292566){\color[rgb]{0,0,0}\makebox(0,0)[lt]{\lineheight{1.25}\smash{\begin{tabular}[t]{l}$P$\end{tabular}}}}%
  \end{picture}%
\endgroup%

	\caption{Lemma \ref{lem: cutting corners} in dimension \(n=2\). For fixed flags \(f\) (dotted lines through \(v(f)\)) and \(f'\) (solid lines), there is a region \(A\) (here in grey) bounded within boxes defined by fixed shifts of elements of \(f\) and a hyperplane of \(f'\) with \(\vol(A) \ll \|v(f)-v(f')\|\). The positive quadrant is denoted by \(P\). In this example, one hyperplane of \(f'\) happens to pass through \(v(f)\), so we use the other one as \(H'\) in Lemma \ref{lem: cutting corners} and its proof.}\label{fig: cutting lemma}
\end{figure}

\begin{proof}
Let \(e_i\) be the \(i\)th standard basis vector of \(\R^n\). By applying a translation we may assume that \(v(f) = 0\). Similarly, without loss of generality (by applying a linear isomorphism), we may assume that \(X = \R^n\), \(f = \{H_i\}_{i=1}^n\), where \(H_i\) is the coordinate aligned hyperplane containing each \(e_j\) except for \(e_i\).

Take any \(H' \in f'\) which does not contain the origin (there is one, since \(f'\) is a flag and \(v(f') \neq 0\)), and let \(\eta = (\eta_1,\ldots,\eta_n)\) be the normal vector of \(H'\), so that \(H' = \{x \in \R^n \mid (x-p) \cdot \eta = 0\}\), for any \(p \in H'\), where we may take \(\|\eta\| = 1\) and \(p = c \eta\) for \(0 < |c| \ll \epsilon\) (since \(v(f') \in H'\) and \(0 < \|v(f) - v(f')\| = \|v(f')\| \ll \epsilon\), with constant only depending on the linear isomorphisms). We may assume that each \(\eta_i \geq 0\), again, by applying a linear isomorphism (if \(\eta_i < 0\) then we may apply the linear isomorphism sending \(e_i \mapsto -e_i\), \(e_j \mapsto e_j\) for \(j \neq i\)) and, by applying \(x \mapsto -x\) if necessary (and replacing \(\eta\) with \(-\eta\)), we may assume that both \(\eta\), \(p\) are in the positive quadrant \(P = \{(x_1,\ldots,x_n) \in \R^n \mid x_i \geq 0\}\).

We note that \(P\) is a connected component of \(X \setminus \bigcup_{H \in f} H\). Consider the subset \(B \subset P\) on the side of \(H'\) containing the origin (so \(B\) is a connected component of Equation \eqref{eq: euclidean minus flags}, except that we do not remove hyperplanes corresponding to \(H \pm b_H\)). That is, we consider
\[
B = \{x \in P \mid (x-p) \cdot \eta < 0\}.
\]
Clearly \(B\) lies above \(V(H') = \{x \in \R^d \mid x \cdot \eta = 0\}\), since \(B \subset P\) and \(\eta \in P\). It follows that \(B\) is contained in the strip of thickness \(\ll \epsilon\) bounded between \(H'\) and \(V(H')\). The basis given by swapping some \(e_j\) for \(\eta\) (using any \(j\) with \(\eta_j \neq 0\)) applies a linear isomorphism distorting volumes in a way bounded by the choices of flags. Then, with respect to this basis, \(B\) is contained in the strip of thickness \(\ll \epsilon\) in the \(\eta\) coordinate direction. By bounding \(B\) to the region \(B'\) that is also below each \(H_i + \nu_i e_i\), where \(\nu_i \ll \nu\), we guarantee that \(x_i \in (0,\nu_i)\) for each \(x \in B'\), so \(B'\) is contained in a box of width at most \(\nu_i \ll \nu\) in each coordinate, and \(\ll \epsilon\) in the \(\eta\) coordinate. Thus, \(\vol(B') \ll \epsilon\) in this basis. So the corresponding region \(A\), before applying the linear isomorphisms, has volume \(\vol(A) \leq C\epsilon\), where \(C\) only depends on the flags \(f\) and \(f'\) and the choice of \(\nu > 0\). It is not hard to see that \(A \neq \emptyset\), for instance, by noting there is a point of \(H'\) in \(P\) and the line segment between this and the origin contains points in \(A\).
\end{proof}

\begin{remark} \label{rem: cutting corners}
In the above proof, if any \(\eta_i = 0\) (such as when \(V(f) = V(f')\), that is, all hyperplanes of \(f'\) are translates of those from \(f\)) then \(A\) is unbounded in the \(e_i\) direction before capping off with \(H \pm b_H\) in Equation \eqref{eq: euclidean minus flags}. On the other hand, in the generic case where each \(\eta_i \neq 0\), we may find a region with volume \(\ll \epsilon^n\).
\end{remark}

The below is the main result of this subsection which will be instrumental in proving our characterisation of \LR\ in Theorem \ref{thm: iff for LR}.

\begin{proposition} \label{prop: close vertices versus small acceptance domains}
Assume that {\Cpx} holds and choose any decomposition satisfying the conclusions of Proposition \ref{prop: factorisation to CSR} (which defines the subspaces \(X_i < \intl\) and \(\delta_i \in \N\)). Suppose that, for sufficiently large \(r\) and any \(v\), \(v' \in \Ver_W(r)\), we have that either \(v_i = v_i'\) or \(\|v_i-v'_i\| \gg r^{-\delta_i}\) for each \(i\), where \(v = v_1 + \cdots + v_\ell\) and \(v' = v'_1 + \cdots + v_\ell'\) for \(v_i\), \(v_i' \in X_i\). Then every acceptance domain \(A \in \sA(r)\) contains a subset \(B_1 + B_2 + \cdots + B_\ell\), where each \(B_i \subseteq X_i\) is a ball of radius \(\gg r^{-\delta_i}\).

Otherwise, given {\Cpx} and a decomposition as above, suppose there exist arbitrarily small \(\epsilon > 0\) and \(r_\epsilon \to \infty\) with \(\|v_i-v_i'\| \leq \epsilon r_\epsilon^{-\delta_i}\) for \(v \neq v' \in \Ver_W(r_\epsilon)\). Then there exist acceptance domains \(A \in \sA(r_\epsilon)\) with volume \(\vol(A) \ll \epsilon r_\epsilon^{-d}\), with constant not depending on \(\epsilon\) or \(r_\epsilon\).
\end{proposition}

\begin{proof}
For any \(A \in \sA(r)\), \(A \supseteq C \in \sC(r)\). We may write \(C = C_1 + \cdots + C_\ell\) where each \(C_i\) is a convex polytope. By assumption, each \(C_i\) has vertices that are distance \(\gg r^{-\delta_i}\) apart, and thus contain a ball of radius \(\gg r^{-\delta_i}\). The latter easily follows from the fact that the bounding subspaces of each convex polytope \(C_i\) belong to a finite set, determining the required constant.

Conversely, suppose we have the \(v\), \(v' \in \Ver_W(r_\epsilon)\) as stated, with flags \(f\), \(f' \in \sF\) and \(v \in \Ver_W(r_\epsilon,f)\), \(v' \in \Ver_W(r_\epsilon,f')\). We write \(\{v\} = \bigcap_{H \in f} (H - (\gamma_H)_<)\).

Without loss of generality (translating \(v\) and \(v'\) by an element in \(\Gamma(c)_<\) for bounded \(c\)), we may assume that \(v\) and \(v'\) are relatively distant from the boundary of \(W\), so all regions constructed below are contained in \(W\). By Proposition \ref{prop: factorisation to CSR}, each \(H \in S_i^\mathrm{c}\) has stabiliser of rank \(k - \delta_i - 1\). Therefore, \(\rk(\Gamma / \Gamma^H) = \delta_i + 1\), so we may find \(g_H \in \Gamma(r_\epsilon)\) with \(H - (\gamma_H)_<\) and \(H - (\gamma_H + g_H)_<\) separated by \(\ll r_\epsilon^{-\delta_i}\).

Consider now only those hyperplanes in \(S_i^\mathrm{c}\) and restrict attention to \(X_i\). To each of the hyperplanes \(H - (\gamma_H)_<\) intersecting to \(v_i\), also add the hyperplanes \(H - (\gamma_H\pm g_H)_<\) to enclose \(v_i\) in a box of hyperplanes with distances \(\ll r_\epsilon^{-\delta_i}\) from \(v_i\). Together with the hyperplanes \(H' - \gamma_H'\) through \(v'_i\) (which, by assumption, has distance \(\ll \epsilon r_\epsilon^{-\delta_i}\) from \(v_i\)) and scaling the picture by \(r_\epsilon^{\delta_i}\), we may directly apply Lemma \ref{lem: cutting corners}; rescaling back (which scales volumes by a factor of \((r_\epsilon^{-\delta_i})^{n_i} = r_\epsilon^{-d_i}\)) we see that the above hyperplanes, together with those through \(v_i\) and \(v_i'\), cut a region \(C_i\) with \(\vol(C_i) \ll \epsilon r_\epsilon^{-d_i}\). Including also the hyperplanes \(H-(\gamma_H)_<\) and \(H-(\gamma_H+g_H)_<\) above for other elements of the decomposition, i.e., with \(H \in S_i\), we cut a region that has cross-section \(C_j\) in the other coordinates with \(\vol(C_j) \ll r_\epsilon^{-d_i}\). Thus, we may enclose the region \(C = C_1 + \cdots + C_\ell \subset W\) with \(\vol(C) \ll \epsilon r_\epsilon^{-d_1} \cdot \cdots \cdot r_\epsilon^{-d_\ell} = \epsilon r_\epsilon^{-d}\). Since each hyperplane used is translated by an element of \(\Gamma(2r_\epsilon)_<\), we have that \(C\) contains a cut region of \(\sC(2r_\epsilon)\), possibly with smaller volume. This in turn contains an \(r'\)-acceptance domain for \(r' \ll r_\epsilon\), by Lemma \ref{lem: refining acceptance domains}, and the result follows.
\end{proof}

\subsection{A necessary and sufficient condition for {\LR}}
We say that \(\cps\) has {\bf positivity of weights} (or {\bf PW}, for short) if \(r\)-patches have frequency \(\gg r^{-d}\), see \cite{BBL} for more details. It is easy to show that \LR\ implies {\bf PW} \cite{BBL} but it will follow from our proofs that the two are equivalent for polytopal cut and project sets. Given {\bf PW}, since the frequencies of \(r\)-patches sum to \(1\), there are \(\ll r^d\) of them, so {\Cpx} holds. See also \cite[Section 6]{KoiWalI} for a brief summary of these notions.

The following well-known result relates volumes of acceptance domains to frequencies. See \cite{Sch98} for further background.

\begin{lemma}\label{lem: PW versus volumes}
A patch \(P\) with acceptance domain \(A_P\) has frequency \(\frac{\vol(A_P)}{\vol(W)}\). Thus, {\bf PW} fails if and only if, for all \(\epsilon > 0\), there exists some \(r\)-patch \(P\) with \(\vol(A_P) \leq \epsilon r^{-d}\).
\end{lemma}

The pieces are now in place to prove our main result characterising {\LR} for all polytopal cut and project sets. Theorems \ref{thm: main1} and \ref{thm: main2} will follow quickly from it.

\begin{theorem} \label{thm: iff for LR}
For a polytopal cut and project set, the following are equivalent:
\begin{enumerate}
	\item {\LR};
	\item {\bf PW};
	\item {\Cpx} holds and, for each pair of flags \(f\), \(f' \in \sF\) and each \(i = 1\), \ldots, \(\ell\), the group \(\Gamma_i[f,f']_<\) is Diophantine in \(X_i\) with respect to the singleton \(\{v(f')_i-v(f)_i\}\).
\end{enumerate}
\end{theorem}

As usual, in the above, the \(\Gamma_i[f,f']_<\) are defined relative to a decomposition satisfying the conclusions of Proposition \ref{prop: factorisation to CSR}.

\begin{proof}
As already stated, {\LR} easily implies {\bf PW}, so suppose that {\bf PW} holds. Then {\Cpx} holds and, by Lemma \ref{lem: PW versus volumes}, for every acceptance domain \(A \in \sA(r)\), we have \(\vol(A) \gg r^{-d}\). Take any pair of flags \(f\), \(f' \in \sF\) and any \(i = 1\), \ldots, \(\ell\). From the converse direction of Proposition \ref{prop: close vertices versus small acceptance domains}, we have that \(v_i = v_i'\) or \(\|v_i - v_i'\| = \|(v-v')_i\| \gg r^{-\delta_i}\) for all \(v\), \(v' \in \Ver_W(r)\). By Proposition \ref{prop: vertices versus flag groups} and Lemma \ref{lem: exact splitting of decomposition}, it follows that \(\|x\| \gg r^{-\delta_i}\) for all non-zero \(x \in \Gamma_i[f,f;r] + (v(f)_i - v(f')_i)\). Equivalently, \(\Gamma_i[f,f]_<\) is Diophantine with respect to \(\{(v(f')_i - v(f)_i\}\) and (3) follows.

So we now assume (3) and will show that {\LR} follows. Take any flag \(f \in \sF\). Then, by (3), \(\Gamma_i[f,f]_<\) is Diophantine with respect to \(\{v(f)_i - v(f)_i = 0\}\), that is, \(\Gamma_i[f,f]_<\) is homogeneously Diophantine. Since \((\Gamma_i)_< \leqslant \Gamma_i[f,f]_<\) (and is of equal rank), it follows that \((\Gamma_i)_<\) is homogeneously Diophantine and thus densely distributed (Theorem \ref{thm: Diophantine => densely distributed}).

By our Diophantine assumption on every flag group, Proposition \ref{prop: vertices versus flag groups} and Lemma \ref{lem: exact splitting of decomposition}, we have that the vertices stay relatively distant in each \(X_i\), that is, they satisfy the hypothesis of the first claim of Proposition \ref{prop: close vertices versus small acceptance domains}. So for any \(A \in \sA(r)\), we have \(A \supseteq B_1 + \cdots + B_\ell\) where each \(B_i \subseteq X_i\) contains a ball of radius \(\gg r^{-\delta_i}\).

Let \(y \in \cps\) and write \(y^\star = y_1 + \cdots + y_\ell \in W\) for \(y_i \in X_i\). By dense distribution, \(y_i + (\Gamma_i(r))_<\) hits all balls of radius \(\gg r^{-\delta_i}\) centred in any given fixed bounded region (say, one large enough to contain the projection of the window to \(X_i\)). Take an arbitrary \(r\)-patch \(P\), with acceptance domain \(A \in \sA(r)\). Since, by the above, \(A\) contains a product of balls of radius \(\gg r^{-\delta_i}\) in each \(X_i\), it follows that there exist \(\gamma_1\), \ldots, \(\gamma_\ell \in \Gamma_i(r')\) with \(y^\star + \gamma_< \in A\) and \(r' \ll r\), where \(\gamma = \gamma_1 + \cdots + \gamma_\ell\). Then \(\|\gamma\|\), \(\|\gamma_\vee\| \ll r\). Since \((y + \gamma_\vee)^\star \in A\), we have that \(P_r(y + \gamma_\vee) = P\) and \(\|(y-\gamma_\vee) - y\| \ll r\). Since \(y \in \cps\) and the \(r\)-patch \(P\) were arbitrary, \(\cps\) is {\LR}.
\end{proof}

By \cite[Theorem 1]{BBL}, a repetitive Delone set \(\cps\) satisfies a subadditive ergodic theorem (SET) if and only if it satisfies positivity of quasiweights. The latter implies {\bf PW} (\cite[Remark 1]{BBL}), which we have seen is equivalent to {\LR}, so SET implies {\LR} for polytopal cut and project sets. Moreover, {\LR} implies SET (by \cite[Theorem 1, Theorem 2]{BBL}, see also \cite{DamLen01}). Thus, Theorem \ref{thm: iff for LR} implies that {\LR} and SET are equivalent for polytopal cut and project sets.

\subsection{Proofs of Theorem \ref{thm: main1} and \ref{thm: main2}}
Theorem \ref{thm: iff for LR} shows that {\LR} is equivalent to {\Cpx} (which may be checked by calculating stabiliser ranks and spans, by Theorem \ref{thm: generalised complexity}), and a Diophantine condition holding: elements of \(\Gamma_i[f,f';r]_<\) stay distant, relative to \(r\), from the point \(f = (v(f')-v(f))_i \in F_i\). It would be preferable to have a Diophantine condition for just a single group, say \(\Gamma_<\) itself (or a splitting of it, in the decomposable case). The below follows easily from Theorem \ref{thm: iff for LR}:

\setcounter{Mainthm}{1}
\begin{Mainthm}
For a polytopal cut and project scheme, if \LR holds then so does {\Cpx} and {\DF}. Conversely, there is some \(N \in \N\) (which may be taken as \(N=1\) in codimension \(n=1\)) so that if {\Cpx} holds and {\DF} still holds after replacing the lattice \(\Gamma\) with \(\frac{1}{N}\Gamma\), then \LR\ holds.
\end{Mainthm}

\begin{remark}\label{rem: value of N}
To be precise in the above, we recall that {\DF} is defined for a specified splitting \(\Gamma_1 + \cdots + \Gamma_\ell\) of \(\Gamma\) up to finite index, in the non-constant stabiliser rank case, and a priori {\DF} may depend on this splitting (subject to the answer to Question \ref{q: general question in DA}). In this case we take the finite index splitting \(\frac{1}{N}\Gamma_1 + \cdots + \frac{1}{N}\Gamma_\ell\) of \(\frac{1}{N}\Gamma\) in the second statement above. We note that a value of \(N \in \N\) satisfying the result can, in theory be extracted from the proof, from the index of \(\Gamma_1 + \cdots + \Gamma_\ell\) in \(\Gamma\) and the indices of \(\Gamma_i\) in each \(\Gamma_i[f,f']\). As such, a value can be determined from just \(\Gamma_<\) and \(\sH_0\) but otherwise independently of \(\sH\), that is, independently of the relative displacements of the \(H \in \sH\).
\end{remark}

\begin{proof}
Condition (3) in Theorem \ref{thm: iff for LR} requires {\Cpx} and we will show it also implies {\DF}. Each \(\Gamma_i \leqslant \Gamma_i[f,f']\), so clearly if \(\Gamma_i[f,f']_<\) is Diophantine relative to \(\{(v(f')-v(f))_i\}\) then the same is true of \((\Gamma_i)_<\). But the scheme satisfies {\DF} if and only if each \((\Gamma_i)_<\) is Diophantine with respect to \(F_i\), which is the case if and only if it is Diophantine with respect to each \(\{(v(f')-v(f))_i\}\), so (3) implies {\DF}. Conversely, we may analogously show that {\DF} for \(\frac{1}{N} \Gamma_<\) implies (3) for some \(N \in \N\), where we note that \(\Gamma_i\) is finite index in each \(\Gamma_i[f,f']\), so each \(\Gamma_i[f,f'] \leqslant \frac{1}{N} \Gamma_i\) for some \(N \in \N\).

In codimension \(n=1\) the scheme is indecomposable. Moreover, it is easy to see that each \(\Gamma[f,f'] = \Gamma\). Indeed, flags just correspond to single window endpoints and thus \(\Ver(f) = \Gamma_< + v(f)\). So each \(\Gamma[f,f'] = \Gamma[f] - \Gamma[f'] = \Gamma - \Gamma = \Gamma\). Thus, {\DF} is precisely the condition that \(\Gamma_<\) is Diophantine with respect to \(F\), which is the same as \(\Gamma[f,f']_< = \Gamma_<\) being Diophantine with respect to \(\{v(f')-v(f)\}\) for each pair \(f\), \(f' \in \sF\).
\end{proof}

In the weakly homogeneous case, the above theorem reduces to checking {\Cpx} and {\D} (the homogeneous Diophantine condition), thus generalising \cite[Theorem A]{KoiWalII} to arbitrary weakly homogeneous polytopal windows:

\setcounter{Mainthm}{0}
\begin{Mainthm}
For a weakly homogeneous polytopal cut and project scheme, the following are equivalent:
\nopagebreak
\begin{enumerate}
	\item \LR;
	\item {\Cpx} and \D.
\end{enumerate}
\end{Mainthm}

\begin{proof}
Suppose that {\LR} holds. By Theorem \ref{thm: main2}, both {\Cpx} and {\DF} hold. Since \(0 \in F\), in particular {\D} must hold. Conversely, suppose that {\Cpx} and {\D} hold. By Lemma \ref{lem: Diophantine under finite index}, {\D} also holds after replacing \(\Gamma\) with any \(\frac{1}{N}\Gamma\), say the value of \(N\) required so as to ensure {\LR}. Since {\D} holds for each \(\frac{1}{N}\Gamma_i\), so does {\DF}, which follows directly from Remark \ref{rem: hom versus inhom Diophantine} and the fact that each \(F_i \subset \Q(\Gamma_i)_<\), since \(F_1 + \cdots + F_\ell = F  \subset \Q\Gamma_< = \Q(\Gamma_1 + \cdots + \Gamma_\ell)_< = \Q(\Gamma_1)_< + \cdots + \Q(\Gamma_\ell)_<\), using Proposition \ref{prop: F for weakly homogeneous}.
\end{proof}

\section{Internal spaces with multiple Euclidean components} \label{sec: multiple components}
A generalisation of our setup is to allow internal space \(\intl \cong \R^n \oplus G \oplus \Z^r\), where \(G\) is a finite Abelian group, that is, the internal space consists of multiple Euclidean components and thus one still has the notion of a window \(W\) being polytopal. In this section we show that such a scheme may always be reduced to one with internal space \(\R^n\) without affecting the MLD class \cite[Section 5.2]{AOI} of the patterns it produces and thus also leaves properties such as \Cpx\ and \LR\ invariant. This process is systematic, so the new window in a single Euclidean component may be described easily in terms of the original window. In this sense, our main theorems characterising \LR\ also cover such examples; see, for instance, the example of the Penrose tilings with internal space \(\R^n \oplus (\Z/5\Z)\) in Section \ref{sec: Penrose}.

We must re-establish the setup for this section. We now take \(\tot \cong \R^k \oplus G \oplus \Z^r\), where \(r \in \N \cup \{0\}\) and \(G\) is a finite Abelian group. We take physical space \(\phy < \tot\) with \(\phy \cong \R^d\) and internal space \(\intl < \tot\) with \(\intl \cong \R^n \oplus G \oplus \Z^r\). These should be complementary, in the sense that \(\phy \cap \intl = \{0\}\) and \(\tot = \phy + \intl\). Of course, up to isomorphism there is no loss of generality in simply taking \(\phy = \R^d\), \(\intl = \R^n \oplus G \oplus \R^n\) and \(\tot = \phy \oplus \intl\) (the agnostic approach of taking merely \(\tot \cong \phy \oplus \intl\), as before, is motivated simply by the fact that it is sometimes natural to fix the total space and vary the other data, such as the lattice or physical space, so we phrase our results in this way which accounts for all such perspectives directly). We fix a lattice \(\Gamma < \tot\), meaning that \(\Gamma\) is discrete and co-compact in \(\tot\). We still assume that \(W \subset \intl\) is compact and equal to the closure of its interior (although, unless otherwise stated, \(W\) does not need to be polytopal in this section). The window may carry a colouring, with each sub-window of the colour partition compact and equal to the closure of its interior and with different sub-windows intersecting on at most their boundaries. We continue to assume the standard properties (1) and (2) from Section \ref{sec: cut and project sets}, that is, the canonical projections \(\pi_\vee\) and \(\pi_<\) from \(\tot\) to \(\phy\) and \(\intl\) are injective on \(\Gamma\). In this case we may assume (3), that \(\Gamma_<\) is dense in \(\intl\), essentially without loss of generality, because we may restrict \(\intl\) to \(\overline{\Gamma_<}\) and the internal space remains of the same form (in contrast to the case where \(\intl \cong \R^n\), where restriction to \(\overline{\Gamma_<}\) can produce \(\Z^r\) factors).

Choosing an identification \(\intl \cong \R^n \oplus G \oplus \Z^r\), we may thus define a canonical projection \(\pi \colon \intl \to \Z^r\). For \(m \in \Z^r\) we let \(W_m \coloneqq W \cap \pi^{-1}(\{m\})\), that is, the part of the window whose \(\Z^r\) component is \(m\). We may choose distinct \(m_i \in \Z^r\), for \(i = 1\), \ldots, \(s\), to write \(W\) as the finite (due to compactness) disjoint union of non-empty sub-windows
\[
W = W_{m_1} \cup W_{m_2} \cup \cdots \cup W_{m_s}.
\]

\begin{lemma} \label{lem: remove Z^r}
For each \(i = 1\), \ldots, \(s\), take any \(\gamma_i \in \Gamma\) with \(\pi^{-1}(\gamma_i)_< \in \pi^{-1}(m_i)\) and with
\[
W' = (W_{m_1} - (\gamma_1)_<) \cup (W_{m_2} - (\gamma_2)_<) \cup \cdots \cup (W_{m_s} - \gamma_s)_< \subset \pi^{-1}(0)
\]
a union of disjoint sets. Then replacing \(W\) with \(W'\) does not affect the MLD class of the cut and project sets produced by the cut and project scheme.
\end{lemma}

\begin{proof}
Let \(\cps\) denote the cut and project set produced from \(W\), and \(\cps'\) the one produced by \(W'\). There is a canonical bijection mapping from \(\cps\) to \(\cps'\), by translating each \(y \in \cps\) with \(\pi(y^\star) = m_i\) to \(y - \gamma_i\). Indeed, \(y \in \cps\) with \(\pi(y^\star) = m_i\) is equivalent to \(y^\star \in W_{m_i}\), or \(y^\star - (\gamma_i)_< \in W_{m_i} - (\gamma_i)_<\) which is equivalent to \(y-(\gamma_i)_\vee \in \cps'\) with \((y-(\gamma_i)_\vee)^\star \in W_{m_i} - (\gamma_i)_<\). For this to be a local derivation we need to check that \(\pi(y^\star) \in \Z^r\) may be determined by the relative displacements of neighbours of \(y \in \cps\) within some globally bounded radius.

Let \(X = \{m_1,m_2,\ldots,m_s\} \subset \Z^r\) be the set of indices in \(\Z^r\) inhabited by the window. For each \(m \in X-X\), take some finite set \(E_m \subset \Gamma_< \cap \pi^{-1}(m)\), large enough so that \(W_{m_i} \subseteq W_{m_j} + E_m\) for all \(i\) and \(j\) with \(m = m_i - m_j\). Take any \(y \in \cps\), say with \(\pi(y^\star) = q \in \Z^r\). For each \(m \in X - X\) we may inspect whether any potential neighbours \(y + \gamma_\vee\) are in \(\cps\) or not, over the finite set of \(\gamma \in \Gamma\) with \(\gamma_< \in E_m\). By construction, this is equivalent to determining if \(\pi(y^\star) + m\) is inhabited by the window or not, equivalently if \(q+m \in X\), that is, \(m \in X-q\). Repeat for each \(m \in X-X\) and define the set \(X'\) by letting \(m \in X'\) if and only if \(m \in X-q\). By construction, we may thus determine \(X' = X-q\) by inspecting occurrence or not in \(\cps\) of potential neighbours of \(y \in \cps\) in a fixed finite set. But since \(X \subset \Z^r\) is finite and fixed, the set \(X-q\) determines \(q = \pi(y^\star)\), as required. The argument for the reverse derivation is analogous.
\end{proof}

The window \(W'\) in the above is contained in the \(0\) component (of the \(\Z^r\) factor), so restricting the scheme to this component does not affect the cut and project set. We have thus shown how to replace a cut and project set with \(\intl \cong \R^n \oplus G \oplus \Z^r\) with one of the form \(\intl \cong \R^n \oplus G\), up to MLD equivalence, and how to define the window for this new cut and project set. So we henceforth restrict to internal space \(\intl \cong \R^n \oplus G\) and aim to show that we may also remove the \(G\) component.

Let \(e \in G\) denote the identity element. To simplify notation, we will henceforth identify \(\tot = \R^k \oplus G\) and \(\intl = \R^n \oplus G\). Explicitly, we really fix an isomorphism \(\varphi \colon \R^k \oplus G \to \tot\) for which \(\varphi\{(x_1,\ldots,x_d,0,\ldots,0,e) \mid x_i \in \R\} = \phy\) and \(\varphi\{(0,\ldots,0,y_1,\ldots,y_n,g) \mid y_i \in \R\} = \phy\) but allow ourselves to write \(\R^k \oplus G\) instead of \(\varphi \R^k \oplus G = \tot\), \((x,g) \in \tot\) rather than \(\varphi (x,g) \in \tot\), etc. For \(g \in G\) we call \(\R^k \times \{g\}\) a {\bf component} of \(\tot\), specifically the {\bf \(g\) component}. We similarly define components of the internal space, denoted
\[
C_g = \{(y,g) \mid y \in \R^n\} \subset \intl\ .
\]
Note that the restriction \(\Gamma_e \coloneqq \Gamma \cap (\R^k \times \{e\})\) to the identity component is a lattice in \(\R^k \times \{e\} \cong \R^k\), with \(\Gamma\) a disjoint union of translates \(\Gamma = \bigcup_{\gamma \in S} (\Gamma_e + \gamma)\), where \(S \subset \Gamma\) has \(|S| = |G|\) and is given by selecting any one element of \(\Gamma\) for each component (which we may do, by the assumption that \(\Gamma_<\) is dense in \(\intl\)). From injectivity of \(\pi_\vee\) on \(\Gamma\), we have that \(\Gamma\) is torsion free and thus \(\Gamma \cong \Z^k\). In particular, the restriction of \(\Gamma_<\) to the \(g\) component, for \(g\) non-trivial, does not contain the Euclidean origin \((0,g) \in \R^n \times \{g\}\), as otherwise this element's pre-image in the lattice is a torsion element.

To remove \(G\) we begin with the following simple observation:

\begin{lemma}
Suppose that \(W\) is coloured so that each element \(W_i \subseteq W\) of the colour partition is wholly contained in a unique component, say component \(C_{g_i}\) for \(g_i \in G\). For each \(i\), take any \(\gamma_i \in \Gamma\) with \((\gamma_i)_< \in C_{g_i}\) in a way so that each \(W_i - (\gamma_i)_<\) disjoint. Then replacing \(W\) with \(W'\), which is defined by replacing each \(W_i\) with \(W_i - (\gamma_i)_<\), does not affect the MLD class of the cut and project sets.
\end{lemma}

\begin{proof}
We may construct a bijection between the cut and project sets \(\cps\) to \(\cps'\) defined, respectively, with windows \(W\) and \(W'\) analogously to the construction at the start of the proof of Lemma \ref{lem: remove Z^r}. Both this and its inverse automatically define local derivations, since by hypothesis we may determine the displacement applied to each \(y \in \cps\) from its colour.
\end{proof}

For cut and project sets where we may always deduce the component of a point using local rules, it thus follows from the above that we may always replace the cut and project scheme with another, up to MLD, whose window is wholly contained in the identity component, which we may then restrict to so as to remove the \(G\) component.

\begin{example} \label{ex: Penrose labelling}
Consider the Penrose cut and project scheme (for further details, see Example \ref{sec: Penrose}). Here, \(G = \Z / 5\Z\), with each component except the identity component inhabited by a pentagonal window. Giving each pentagon a unique colour does not affect the cut and project scheme, up to MLD, since this information can be deduced from local information in the resulting cut and project sets. Indeed, consider the pentagons \(P_i\) contained in component \([i] \in \Z / 5\Z\), for \(i = 1\), \ldots, \(4\). Take a set \(S \subset \Gamma\) so that each \(\gamma \in S\) is contained in component \([1]\), and \(P_{i+1} \subseteq P_i + S_<\) for \(i = 1\), \(2\) and \(3\). Given a Penrose cut and project set \(\cps\), the points \(y \in \cps\) that are projections from component \([4]\) are recognised, through a local rule, as those with no neighbours \(y + s_\vee\) for \(s \in S\). Similarly, those \(y \in \cps\) coming from component \([3]\) can be locally recognised as those who have a neighbour \(y + s_\vee\) in component \([4]\), and so on.

It follows that attaching the information of the component of each point as a colour does not change its MLD class. Then, by the previous lemma, we may replace the Penrose cut and project scheme, up to MLD, with one whose internal space is Euclidean with a single component, and has window a disjoint union of lattice translates of the \(P_i\) into the identity component (we can then also discard their colours, up to MLD, by the argument in Proposition \ref{prop: unlabel}). Of course, this will not preserve special symmetry properties of the Penrose patterns, such as 10-fold statistical rotational symmetry, but that is not important for their (translational) MLD class. Also note that the restriction \((\Gamma_e)_<\) of \(\Gamma_<\) to the identity component, and the collection \(\sH + (\Gamma_e)_<\) of supporting hyperplanes, up to lattice translations, retain 10-fold symmetry.
\end{example}

It would be desirable to be able to systematically determine the internal space components of points in a cut and project scheme using local rules, in a manner similar to the above example. Unfortunately this is not as simple as the case of \(\Z^r\) components, as dealt with in Lemma \ref{lem: remove Z^r}. Consider, for example, a window \(W\) that is identical in each \(G\) component \(C_g = \R^n \times \{g\} \cong \R^n\). Then the acceptance domains remain identical in each component too (lattice translates of the window and its complement, which are intersected to the acceptance domains, are translated together), so there is no way to determine the component of \(y^\star\) from rules local to \(y \in \cps\). However, in this example one may simplify the cut and project scheme by noting that the \(G\) factor is superfluous: we may collapse it, which does not change the cut and project sets produced. We are thus motivated to introduce the following:

\begin{definition}
We define \(\Aut(W) = \{v \in \intl \mid W = W+v\}\).
\end{definition}

It is easy to see that \(\Aut(W)\) is a group under addition in \(\intl\). For \(v = (x,g) \in \R^n \oplus G = \intl\), if \(v \in \Aut(W)\) then necessarily \(x = 0\) (since \(W\) is bounded), so we henceforth identify \(\Aut(W)\) with a subgroup of \(G\). When \(\Aut(W)\) is non-trivial, we may reduce our scheme:

\begin{lemma} \label{lem: mod out Aut(W)}
Given a cut and project scheme \(\mathcal{S}\) define a new one, \(\mathcal{S}'\), by replacing \(\intl = \R^n \oplus G\) with the quotient \(\intl \oplus (G/\Aut(W))\). More precisely, let \(q\) denote the projection \(q \colon G \to G/\Aut(W)\) and also the derived projections \(q \colon \R^n \oplus G \to \R^n \oplus (G/\Aut(W))\) and \(q \colon \tot \to \tot' \coloneqq \phy \oplus (\R^n \oplus (G/\Aut(W)))\). We take the window of \(\mathcal{S}'\) to be \(W' \coloneqq q(W)\) (with obvious meaning if \(W\) is coloured), and lattice \(\Gamma' \coloneqq q(\Gamma)\). Then \(\mathcal{S}\) and \(\mathcal{S}'\) produce identical cut and project sets.
\end{lemma}

\begin{proof}
As usual, we assume that \(\mathcal{S}\) is defined with \(W\) in non-singular position. This implies the same for \(\mathcal{S}'\). Indeed, suppose that \(g_< \in \partial W'\) for some \(g \in \Gamma'\). Let \(\gamma \in \Gamma\) with \([\gamma] = g\). Then \(\gamma_< \in \partial (W+v)\) for some \(v \in \Aut(W)\). By definition of \(\Aut(W)\), we have \(W = W+v\) and thus \(\gamma_< \in \partial W\), a contradiction, so \(W'\) is also in non-singular position.

It follows directly from our definitions that \(\pi_\vee\) remains injective on the lattice. Similarly, if \(\pi_<(q(\gamma)) = 0\), for \(\gamma \in \Gamma\), then \(\gamma\) is a torsion element in \(\Gamma_< \cong \Gamma \cong \Z^k\) so \(\gamma = 0\). Clearly \(\Gamma'_<\) is still dense in the internal space. It is trivial to check that all other standard conditions we are assuming for cut and project sets still hold for \(\mathcal{S}'\).

Let \(\cps\) and \(\cps'\) be the cut and project sets produced by \(\mathcal{S}\) and \(\mathcal{S}'\), respectively. Let \(y \in \cps\), equivalently \(\gamma_< \in W\) for \(\gamma \in \Gamma\) and \(y = \gamma_\vee\). Then \(q(\gamma_<) = (q(\gamma))_< \in q(W) = W'\), so that \(y = q(\gamma)_\vee \in \cps'\) and \(\cps \subseteq \cps'\). Conversely, let \(y \in \cps'\), so \(y = q(\gamma)_\vee\) and \((q(\gamma))_< \in q(W)\) for some \(q(\gamma) \in q(\Gamma) = \Gamma'\). Then \(\gamma_< \in W + v\) for some \(v \in \Aut(W)\). But \(W = W+v\), so \(\gamma_< \in W\) and thus \(y = \gamma_\vee \in \cps\). It follows that \(\cps \subseteq \cps'\) and thus \(\cps = \cps'\). Analogous arguments hold if the window is coloured, by replacing \(W\) with \(W_i\) in the above arguments, showing that the Delone sets \(\cps_i = \cps_i'\) of colour \(i\) agree for each \(i\).
\end{proof}

So, without loss of generality, we may assume that \(\Aut(W)\) is trivial, as otherwise we may reduce \(G\) using the above lemma. In this case \(W \neq W+(0,g)\) for all non-trivial \(g \in G\), that is, there are at least two components of the internal space, related by the translation \(g\), that contain different sub-windows. We wish to exploit this, in a similar way to Example \ref{ex: Penrose labelling}, to be able to colour points \(y \in \cps\) using local rules (i.e., checking existence or not of potential neighbours with a finite radius) which determine the \(G\) components of each \(y^\star\). It seems necessarily rather technical to derive an approach systematically covering the general case, but one may be given as described in the proof below.

\begin{proposition} \label{prop: trivial Aut(W)}
Suppose that \(\Aut(W)\) is trivial. Then colouring the restrictions of \(W\) to each component with a distinct colour each does not affect the MLD class of the cut and project sets.
\end{proposition}

\begin{proof}
We let \(W_a \coloneqq W \cap C_a\) (the restriction of \(W_a\) to the \(C_a\) component). Take any non-trivial \(g \in G\) and set \(v_g \coloneqq (0,g) \in \intl\). We wish to create an acceptance domain \(X_g\) (that is, a subset given by a finite intersection of \(\iW\) with \(\Gamma_<\) translates of \(\iW\) and \(W^\mathrm{c}\)) so that exactly one of \(X_g \cap C_a\) and \(X_g \cap C_{a+g}\) is non-empty, for some \(a \in G\).

Since \(\Aut(W)\) is trivial, there is some \(a \in G\) with \(W_a + v_g \neq W_b\), where \(b = a+g\). If the window is coloured, this holds for the partition of the window in some colour \(i\) (and we simplify notation by simply letting \(W\) denote this element of the colour partition).

Since \(W_a + v_g \neq W_b\), there must either be an interior point \(v \in W_a\) with \(v + v_g \notin W_b\), or an interior point \(v \in W_b\) with \(v-v_g \notin W_a\) (in either case we may also assume that \(v\) is non-singular, in the sense that \((v + \Gamma_<) \cap \partial W = \emptyset\)). We assume the former case; otherwise, we may replace the roles of \(C_a\) and \(C_b\) below, and replace each occurrence of \(g\) with \(-g\).

By density of \(\Gamma_<\), clearly we may take a finite subset of \(S \subset \Gamma_e \coloneqq \Gamma \cap (\R^k \times \{e\})\) so that the intersection
\[
X_1 = \iW \cap \bigcap_{\gamma \in S} (\iW+\gamma_<)
\]
defines an open neighbourhood about \(v = (x,a) \in W_a\) with arbitrarily small diameter \(\epsilon > 0\) when restricted to \(C_a\). If \(X_1 \cap C_b = \emptyset\) then, since \(X_1 \cap C_a \neq \emptyset\), we are done, by setting \(X_g \coloneqq X_1\). Otherwise, since \(v+v_g \notin W_b\), we may take \(\epsilon\) sufficiently small so that \((X_1 \cap C_a) + v_g\) has trivial intersection with \(W\) (since an open neighbourhood of \(v+v_g \in C_b\) has trivial intersection with \(W\)). We show that a nested sequence \(X_1 \supseteq X_2 \supseteq \cdots \supseteq X_m\) of acceptance domains can be constructed, terminating with an \(X_m\) satisfying \(X_m \cap C_a = \emptyset\) and \(X_m \cap C_b \neq \emptyset\).

To construct \(X_2\), first take \(z_1 = (x_1,g) \in \Gamma_<\) with \((X_1 \cap (X_1 + z_1)) \cap C_b \neq \emptyset\). Since, by assumption, \(X_1 \cap C_a \neq \emptyset\) and \(X_1 \cap C_b \neq \emptyset\), this can be achieved, and we may also assume that \(\|x_1\|\) is much larger than \(\epsilon\), since \(X_1 \cap C_a\) has diameter arbitrarily small \(\epsilon > 0\) with \((X_1 \cap C_a) + (0,g)\) exterior to \(W\) and thus also exterior to \(X_1 \cap C_b \subset W\). We define \(X_2 \coloneqq X_1 \cap (X_1 + z_1)\). By construction, \(X_1 \supseteq X_2\) and \(X_2 \cap C_b\) is contained in an \(\epsilon\)-neighbourhood of \(v + z_1 = (x+x_1,g)\). If \(X_2 \cap C_a = \emptyset\), then we are done by taking \(X_g \coloneqq X_2\)

Otherwise, we continue inductively. Suppose we have constructed \(X_1 \supseteq X_2 \supseteq \cdots \supseteq X_j\) where each \(X_{i+1} = X_i \cap (X_i + z_i)\) for \(z_i = (x_i,g) \in \Gamma_<\) with \(\|x_1 - x_i\| \leq 2\epsilon\). Suppose also that each \(X_i \cap C_a \neq \emptyset\) and \(X_i \cap C_b \neq \emptyset\). Then we may find \(z_j = (x_j,g) \in \Gamma_<\) with \((X_j \cap (X_j + z_j)) \cap C_b \neq \emptyset\). Note that, since \(X_j \cap C_a \subseteq X_1 \cap C_a\), with the latter in an \(\epsilon\)-neighbourhood of \(v\), and \(X_j \cap C_b \subseteq X_2 \cap C_b\) is contained in a \(\epsilon\)-neighbourhood about \(v+z_1\), we must take \(\|x_j-x_1\| \leq 2\epsilon\). We may now take \(X_{j+1} \coloneqq X_j \cap (X_j + z_j)\). By construction, \(X_{j+1} \cap C_b \neq \emptyset\). If \(X_{j+1} \cap C_a = \emptyset\) we are done, otherwise we continue inductively.

We claim that, eventually, some \(X_m \cap C_a = \emptyset\) (but, by construction, \(X_m \cap C_b \neq \emptyset\) at each step). Indeed, consider taking \(m\) equal to the order of \(g\). By construction,
\begin{align*}
X_{m+1} \cap C_{a+mg} & \subseteq (X_m \cap C_{a+(m-1)g}) + z_m,\\
X_m \cap C_{a+(m-1)g} & \subseteq (X_{m-1} \cap C_{a+(m-2)g}) + z_{m-1},\\
                      & \ldots,
\end{align*}
\begin{align*}
X_3 \cap C_{a+2g}     & \subseteq (X_2 \cap C_{a+g}) + z_2,\\
X_2 \cap C_{a+g}      &\subseteq (X_1 \cap C_a) + z_1,
\end{align*}
where \(\|z_i - z_1\| \leq 2\epsilon\) for each \(i\). Combining, we see that
\[
X_{m+1} \cap C_a = X_{m+1} \cap C_{a+mg} \subseteq (X_1 \cap C_a) + z_1 + z_2 + \cdots + z_m.
\]
But \(z_1 + \cdots + z_m = (x_1,g) + \cdots + (x_m,g) = (x_1 + \cdots + x_m,e)\) and \(\|(x_1 + \cdots + x_m) - m x_1\| \leq 2(m-1)\epsilon\) (since each \(\|x_i-x_1\| \leq 2\epsilon\)). So the Euclidean part of \(z_1 + \cdots + z_m\) is approximately equal to \(mx_1\), where by assumption \(\|x_1\|\) was taken much larger than \(\epsilon > 0\). But \(X_1 \cap C_a\) has diameter at most \(\epsilon\), so clearly \((X_1 \cap C_a) \cap ((X_1 \cap C_a) + z_1 + \cdots +z_m) = \emptyset\) . Since \(X_{m+1} \subseteq X_1\), we see that \(X_{m+1} \cap C_a = \emptyset\), as required. It follows that we may find an acceptance domain \(X_g\) with exactly one of \(X_g \cap C_a\) or \(X_g \cap C_{a+g} = \emptyset\).

We may now easily create acceptance domains \(X'_g\) so that at exactly one of \(X'_g \cap C_a\) or \(X'_g \cap C_{a+g} = \emptyset\), for some \(a \in G\), and \(X'_g \cap W_a = W_a\) or \(X'_g \cap W_a = \emptyset\) for each \(a \in G\). Indeed, take a finite subset \(S \subset \Gamma_e\) large enough so that \(W_a \subseteq (X_g \cap C_a) + S_<\) whenever \(X_g \cap C_a \neq \emptyset\) and define \(X'_g \coloneqq (X_g + S_<) \cap W\).

It can be determined whether any given \(y \in \cps\) is such that \(y^\star\) is in a component inhabited by \(X'_g\) or not, by the fact that we constructed \(X'_g\) using unions and intersections of \(\Gamma_<\) translates of \(\iW\) (since a given \(y \in \cps\) being in any given \(W-\gamma_<\) or not is equivalent to \(y + \gamma_\vee \in \cps\) or not). Moreover, this distinguishes points \(y_1\), \(y_2 \in \cps\) whenever \(y_1^\star \in C_a\) and \(y_2^\star \in C_{a+g}\) for some \(a \in G\) (with \(a\) depending on \(g\)). In other words, we may add further colours to our window which realise these distinctions without changing the MLD class of \(\cps\), by extending any current labelling scheme with additional label either \(g\) or \(\overline{g}\) on each component \(W_a\) of the window, according to whether, respectively, \(W_a \cap X'_g = W_a\) or \(W_a \cap X'_g = \emptyset\).

We may repeat the above for each non-trivial \(g \in G\). For each \(a \in G\) define a finite subset \(S_a \subset \Gamma\), contained in the \(a\) component of \(\Gamma\) and large enough so that \(W_b + (S_a)_< \supseteq W_{a+b}\) whenever \(W_b \neq \emptyset\). Add further colours to the window, by giving \(W_b\) colour \((a,g)\) whenever \(b+a \in X'_g\), and colour \((a,\overline{g})\) if \(b+a \notin X'_g\). At the level of cut and project sets, this may still be derived by local rules, since it just involves checking if neighbours \(y+s\) have colour \(g\) or \(\overline{g}\), over all \(s \in (S_a)_\vee\), which in turn is local by the constructions above. Now, consider two distinct window components \(W_p\) and \(W_q\). Let \(g = q-p \neq e\). There are components \(C_u\) and \(C_{u+g}\) with one inhabited by \(X'_g\), the other not. Take \(a = u-p\). Then the components of \(p+a = u\) and \(q+a = u+g\) have distinct colours \(g\) and \(\overline{g}\), and thus one of \(W_p\) and \(W_q\) is coloured \((a,g)\) and the other \((a,\overline{g})\). Since \(p \neq q\) were arbitrary and these components have now been assigned distinct colours without affecting MLD classes, the proof is complete.
\end{proof}

Combining the above results in this section, we have:

\begin{corollary} \label{cor: reduce to Euc}
For a cut and project scheme with internal space \(\intl \cong \R^n \oplus G \oplus \Z^r\), for \(G\) finite and Abelian, we may construct another, with internal space \(\mathbb{R}^n\), which produces MLD cut and project sets.
\end{corollary}

We note that the proofs did not require the window to be polytopal, merely that \(W\) is topologically regular. The previous results are of course more constructive than the summary given by Corollary \ref{cor: reduce to Euc}: we may explicitly describe the new cut and project scheme. Firstly we remove the \(\Z^r\) component by using \(\Gamma_<\) translates to move each \(\Z_r\) component of the window into disjoint subwindows of the \(\Z^r\) identity component, which does not affect MLD classes by Lemma \ref{lem: remove Z^r}, and restrict to this component. Next we take the quotient \(q\) by \(\Aut(W)\), when non-trivial, which does not affect MLD classes by Lemma \ref{lem: mod out Aut(W)}. By Proposition \ref{prop: trivial Aut(W)} the required cut and project scheme then has lattice the restriction of \(q(\Gamma)\) to the identity component, and the new window is given by taking disjoint \(\Gamma_<\) translates of each window component into the identity component. In particular, we may directly describe the data \(\Gamma_<\) and \(\sH\) required for applying Theorems \ref{thm: main1} and \ref{thm: main2}, and the more general Theorem \ref{thm: iff for LR}.

\section{Examples} \label{sec: examples}

\subsection{Codimension one}
In codimension \(n=1\) the window \(W\) is a finite union of intervals. As \(n=1\), the scheme is indecomposable and, by Theorem \ref{thm: generalised complexity}, {\Cpx} always holds. We have the set \(F\) of displacements between pairs of window (or label partition) endpoints. The scheme is weakly homogeneous if there is some \(N \in \N\) with \(F \subset \frac{1}{N}\Gamma_<\) (e.g., by Proposition \ref{prop: F for weakly homogeneous}). In this case Theorem \ref{thm: main1} applies and thus {\LR} is equivalent to {\D}, that is, that \(\Gamma_<\) is Diophantine as a rank \(k\) densely embedded lattice in \(\intl \cong \R\). In the non-weakly homogeneous case, by Theorem \ref{thm: main2}, {\LR} is equivalent to {\DF}.

For \(k = 2\), we have \(\Gamma_< < \intl\) is linearly isomorphic to \(\alpha\Z + \Z < \R\), where \(\alpha\) may be taken as the ratio of lengths of any two generators of \(\Gamma_<\). This satisfies {\D} if and only if \(\alpha \in \R\) is badly approximable (Example \ref{ex: 2-to-1 Diophantine}), which is necessary for {\LR} (and sufficient in the weakly homogeneous case). This includes all quadratic irrationals, as their continued fractions are eventually periodic and thus bounded. Generally, for {\LR} it is necessary and sufficient that {\DF} holds, i.e., \(\alpha \in \Bad(y)\) for each \(y \in F\). This is much more restrictive than simply \(\alpha \in \Bad = \Bad(0)\), which already limits \(\alpha\) to a set of measure \(0\), since \(\alpha \notin \Bad(x)\) for almost every \(x\) when \(\alpha \in \Bad(0)\), by \cite{Kur55}. However, the set of \(\alpha \in \bigcap_{y \in F} \Bad(y)\) still has full Hausdorff dimension (in particular, it is uncountable), since each \(\Bad(y)\) is winning, see \cite[Remark 7.8]{BRV16}.

\subsection{The Ammann--Beenker tilings} \label{sec: AB}
We refer the reader to \cite[Example 7.8]{AOI} for the construction of the \(4\)-to-\(2\) Ammann--Beenker cut and project scheme, although we recall the main details here. It is already known that the Ammann--Beenker tilings satisfy {\LR}, since they are generated by substitution. Nonetheless, they are useful examples to demonstrate how to apply our main theorem.

We define \(\Gamma_\vee = \Z[\xi_8]\), where \(\xi_8 = \exp(2\pi i/8) \in \C \cong \R^2\) is the standard primitive eighth root of unity. Then \(\Gamma_\vee = \{n_0 + n_1 \xi_8 + n_2 \xi_8^2 + n_3 \xi_8^3 \mid n_i \in \Z\}\). The star map \(x \mapsto x^\star\) into internal space \(\intl \cong \C \cong \R^2\) is defined by
\[
n_0 + n_1 \xi_8 + n_2 \xi_8^2 + n_3 \xi_8^3 \mapsto n_0 + n_1 \xi_8^3 + n_2\xi_8^2 + n_3\xi_8,
\]
in short, \(\xi_8 \mapsto \xi_8^3\) and extended as a ring automorphism. This defines the Minkowski embedding \(x \mapsto (x,x^\star) \in \C \times \C \cong \R^4\), which defines our lattice \(\Gamma \cong \Z^4\) in the total space, with projections \((x,x^\star)_\vee \coloneqq x\), \((x,x^\star)_< \coloneqq x^\star\). Finally, the window \(W\) (see Figure \ref{fig: AB generalised vertices}) is defined as a centred regular octagon of edge length \(1\) (in this position, the window is in fact non-singular). Alternatively, \(\Gamma\) may be taken as a rotation of \(\Z^4 \leqslant \R^4 = \tot\), with \(W\) (up to translation) a projection of the associated unit hypercube in \(\R^4\) with vertices in \(\Gamma\). From this it follows that each supporting hyperplane (after an appropriate, singular translation of \(W\)) intersects \(\Gamma_<\), so the scheme is homogeneous and Theorem \ref{thm: main1} applies.

We first check {\Cpx} by considering the set \(\sH_0\) of supporting subspaces. This is the 8-fold symmetric configuration of 4 lines through the origin, of the form \(\ell_n \coloneqq \langle \xi_8^n \rangle_\R\). We consider their stabilisers, whose projections to internal space are the intersections \(\ell_n \cap \Gamma_<\). For example, \(\xi_8^0 = 1\), \(\xi_8 + \xi_8^{-1} = \sqrt{2}\in \ell_0\) give two linearly independent elements in \(\Gamma_<\) generating the (projection of) the stabiliser of \(\ell_0\), which we see is rank \(2\) (clearly not higher, or \(\Gamma_<\) could not be dense). By rotational symmetry the same is true of the other \(\ell_n\) and \(\rk(H) = 2\), \(\beta_H = 1\) for each \(H \in \sH\). By Theorem \ref{thm: generalised complexity}, \(p(r) \asymp r^2\) so {\Cpx} holds.

Note that the projected stabiliser considered above, \(\langle 1,\sqrt{2} \rangle_\Z\), is Diophantine as a densely embedded subgroup of \(\ell_0 \cong \R\), since \(\sqrt{2}\) is badly approximable (see Example \ref{ex: 2-to-1 Diophantine}). Similarly, the other stabilisers are Diophantine in their supporting subspaces, for example \(G_1 \coloneqq \langle \xi_8,\sqrt{2}\xi_8 \rangle_\Z \leq \ell_1\) is Diophantine as a subgroup of \(\ell_1\). Sums of Diophantine lattices of equal rank and embedding dimension (into complementary spaces) are still Diophantine (more generally, see \cite[Lemma 5.5]{KoiWalII}), so that \(G_0 + G_1\) is Diophantine in \(\ell_0 + \ell_1 = \intl\). Since \(\rk(G_0 + G_1) = 4\) it is finite index in \(\Gamma_<\), which is thus Diophantine by Lemma \ref{lem: Diophantine under finite index}. Since the scheme is constant stabiliser rank, by Definition \ref{def: Diophantine scheme} property {\D} holds. Since {\Cpx} also holds, {\LR} follows from Theorem \ref{thm: main1}. Of course, the same would be true if we replace the standard Ammann--Beenker window with any other which has the same set \(\sH_0\) of supporting hyperplanes and retains weak homogeneity, which could be non-convex or disconnected.

\subsection{The decorated Ammann--Beenker tilings}
It is known that the Ammann--Beenker tilings do not admit local rules \cite{Bur88, BF15} but they may be decorated (in a non-local way) to define tilings that do; see \cite{Le97, Katz95, Gah93} for more details. Such tilings are MLD to the labelled cut and project sets with scheme as above, but where \(W\) is partitioned into eight triangles by cutting the octagon with straight lines between opposite vertices. Although {\LR} is known to follow from a substitution rule, again this is a good example to demonstrate our main theorem.

We now have supporting subspaces \(\sH = \{\langle \xi_{16}^n \rangle_\R \mid n = 0, 1, \ldots, 7\}\). Notice that, since the new supporting hyperplanes (\(\langle \xi_{16}^n \rangle_\R\) for \(n\) odd) still pass through vertices of \(W\), the scheme is still homogeneous. So we just need to consider the stabiliser ranks of these new supporting subspaces, for example \(V = \langle \xi_{16} \rangle_\R\). Clearly this line contains the midpoint \(\frac{1}{2}(1 + \xi_8)\) between \(1\) and \(\xi_8\), so \(1 + \xi_8 \in V \cap \Gamma_<\). In a similar way, \(\sqrt{2} + \sqrt{2}\xi_8 = (\xi_8 + \xi_8^{-1}) + (\xi_8^2 + 1) \in V \cap \Gamma_<\), so we have two independent elements and \(\rk(V) = 2\). The same holds for the other new hyperplanes, by \(8\)-fold rotational symmetry. Thus, {\Cpx} holds by Theorem \ref{thm: generalised complexity}; we recall that all theorems apply to windows with polytopal partitions by Proposition \ref{prop: unlabel} and taking \(\sH\) as the hyperplanes supporting the partition regions (in particular, in applying Theorem \ref{thm: generalised complexity} to determine the complexity, and checking constant stabiliser rank in applying Definition \ref{def: Diophantine scheme}). Since the scheme is constant stabiliser rank, {\D} still holds (as \(\Gamma_<\) remains unchanged), so the decorated Ammann--Beenker tilings are {\LR} by Theorem \ref{thm: main1}.

\subsection{Penrose tilings} \label{sec: Penrose}
A comprehensive description of the Penrose cut and project scheme may be found in \cite[Example 7.11]{AOI}. This scheme generates Delone sets MLD to Penrose's famous tilings \cite{Pen79}. We now take \(\Gamma_\vee = \Z[\xi_5] < \C \cong \R^2\), with \(\xi_5 = \exp(2\pi i/5)\) the standard primitive fifth root of unity. This is rank \(4\); note that \(\sum_{i=0}^4 \xi_5^i = 0\). One may define a star map into \(\C \cong \R^2\) by extending \(\xi_5 \mapsto \xi_5^2\) to a ring homomorphism, and indeed this is used for the cut and project scheme of the T\"{u}bingen triangle tilings \cite[Example 7.10]{AOI}. However, for the Penrose cut and project scheme a cyclic component is introduced. We have the homomorphism \(\kappa \colon \Z[\xi_5] \to \Z/5\Z\) defined by
\[
\kappa \left( \sum_j m_j \xi_5^j\right) \coloneqq \sum_j m_j \mod 5,
\]
and the star map \(\star \colon \Z[\xi_5] \to \Z[\xi_5] \times (\Z/5\Z)\) is induced from \(\xi_5 \mapsto \xi_5^2\) and tracking \(\kappa\), that is (noting \(\xi_5^4 = -1 - \xi_5 - \xi_5^2 - \xi_5^3\)):
\[
n_0 + n_1 \xi_5 + n_2 \xi_5^2 + n_3 \xi_5^3 \mapsto ((n_0-n_2) + (n_1-n_2) \xi_5^2 - n_2 \xi_5^3 + (n_3-n_2) \xi_5,\kappa(x)).
\]
Then the Minkowski embedding defines a \(4\)-to-\(2\) cut and project scheme with lattice \(\Gamma = \{(x,x^\star) \mid x \in \Z[\xi_5]\} \cong \Z^4\) and internal space \(\C \times (\Z/5\Z) \cong \R^2 \times (\Z/5\Z)\). The window is given by \(W = W^{(1)} \cup W^{(2)} \cup W^{(3)} \cup W^{(4)}\), where each \(W^{(i)} \subset \R^2 \times \{[i]\}\). More precisely (see \cite[Example 7.11]{AOI}),
\[
W^{(1)} = P, \ W^{(4)} = -P, \ W^{(3)} = \varphi P \ \text{ and } \ W^{(2)} = -\varphi P,
\]
where \(\varphi = 1 + \xi_5 + \xi_5^4\) is the golden ratio (with \(\kappa(\varphi) =3\)) and \(P \subset \R^2 \times \{[1]\}\) is the regular pentagon, given as the convex hull of \(\{\xi_5^j\}_{j=0}^4\). Alternatively, one may imagine an internal space \(\intl \cong \R^2\) but where, to determine if \(\gamma_\vee \in \cps\), one checks if \(\gamma_<\) belongs to a window depending on \(\kappa(\gamma)\). This viewpoint falls out of the alternative dualisation method, see \cite[Section 5]{BKSZ90}.

By the results of Section \ref{sec: multiple components}, our machinery is able to deal with internal spaces with cyclic components, by shifting all window components into \(\R^2 \times \{[0]\}\) using elements of \(\Gamma_<\). This produces MLD cut and project sets to the Penrose cut and project scheme, with internal space \(\R^2\) and window \(W = \bigcup_{i=1}^4 W^{(i)} + (\gamma_i)_<\) a disjoint union of \(4\) pentagons, where each \(\gamma_i \in \Gamma\) and lattice restricted to the index 5 subgroup \(G = \{\gamma \in \Gamma \mid \kappa(\gamma) = 0\} < \Gamma\). Since each vertex of each \(W^{(i)} + \gamma_i\) is contained in \(G_<\), we see in particular that the scheme is homogeneous, so that Theorem \ref{thm: main1} applies. So we just need to check {\Cpx} and, since the scheme is indecomposable, that \(\Gamma_<\) is Diophantine.

Both of these checks are similar to the Ammann--Beenker tiling. For example, take the supporting subspace \(V = \langle 1 \rangle_\R < \C\). This intersects the \(\Z\)-independent elements \(5 \in G_<\) and \(5\varphi = 5 + 5\xi_5 + 5\xi_5^4 \in G_<\), so \(\rk(V) = 2\). By rotation, we similarly find \(\rk(V) = 2\) for all \(V \in \sH_0\) and thus {\Cpx} holds by Theorem \ref{thm: generalised complexity}. We have that \(G^V_< < V\) is linearly isomorphic to \(\langle 1, \varphi \rangle_\Z < \R\), which is Diophantine since \(\varphi\) is badly approximable (Example \ref{ex: 2-to-1 Diophantine}), and similarly for other stabilisers; the sum of two such, like for the Ammann--Beenker, thus gives a rank 4 Diophantine sublattice of \(G_<\), which thus itself is Diophantine and so {\D} holds.

\subsection{Generalised Penrose Tilings} \label{sec: generalise Penrose}
The Penrose tilings are often described as \(5\)-to-\(2\) cut and project sets, with window the projection of the \(5\)-dimensional hypercube, originally observed through De Bruijn's grid method \cite{dBr81}, see also \cite[Section 7.5.2]{AOI}. This description can be inconvenient for analysis, as \(\Gamma_<\) is not dense in \(\intl \cong \R^3\) and instead \(\overline{\Gamma_<} \cong \R^2 \times \Z\). The projection \(\pi_\vee\) is also not injective on \(\Gamma\); the description above, as a \(4\)-to-\(2\) cut and project set with cyclic internal space component is derived from a quotient of this one, restricted to internal space \(\intl \cong \R^2 \times \Z\), which slices the window into \(4\) pentagons (and isolated vertices in two components, which may be discarded, since one considers non-singular cut and project sets and their limits), see \cite[Figure 7.8]{AOI}. However, the \(5\)-to-\(2\) description naturally defines an interesting family of cut and project sets, parametrised by \(\gamma \in \R\), given by translating the window in the discrete direction, defining (up to MLD) the \emph{generalised} Penrose tilings. See \cite{PK87}, \cite[Remark 7.8]{AOI} and \cite[Section 7.5.2]{AOI} for further details.

The generalised Penrose tilings may still be described by cut and project schemes with internal space \(\intl \cong \R^2 \times (\Z/5\Z)\) and, by the results of Section \ref{sec: multiple components}, can then be taken up to MLD as having internal space \(\intl \cong \R^2\) and lattice \(G\) as above for the standard Penrose tilings. Then \(W\) is a union of pentagons or decagons, with supporting subspaces \(\sH_0\) unchanged, so {\Cpx} still holds by Theorem \ref{thm: generalised complexity} and, since \(G_<\) remains the same and the scheme is indecomposable, {\D} also holds.

It follows that the generalised Penrose tilings that satisfy weak homogeneity are {\LR}, by Theorem \ref{thm: main1}. Presumably weak homogeneity holds precisely when \(\gamma \in \Q\). Otherwise, one needs to calculate the displacement set \(F\) between the generalised vertices, which depends on \(\gamma\). Then, by Theorem \ref{thm: main2}, {\LR} fails if \(G_<\) is not Diophantine with respect to \(F\), and there is some \(N \in \N\) so that {\LR} holds if \(\frac{1}{N}G_<\) is Diophantine with respect to \(F\). By Remark \ref{rem: value of N}, the value \(N\) may be taken independently of the particular displacements \(\sH\) of the supporting subspaces, and thus independently of the parameter \(\gamma\). Of course, it is tempting to hope that Question \ref{q: general question in DA} is answered in the positive, in which case this \(N\) can be ignored and {\DF} is both necessary and sufficient for {\LR}, given {\Cpx}.

\printbibliography

\end{document}